\documentclass{amsart}
\usepackage{amsmath,amssymb} % use AMS symbol
\usepackage{amscd}
\usepackage[dvipdfm]{graphicx}

\theoremstyle{plain}
\newtheorem{thm}{Theorem}[section]
\newtheorem{lem}{Lemma}[section]

\newtheorem{prop}{Proposition}[section]
\newtheorem{clm}{Claim}[section]
\theoremstyle{definition}
\newtheorem{defn}{Definition}[section]
\theoremstyle{remark}
\newtheorem*{rem}{Remark}
\newtheorem*{prf}{Proof}

\newtheorem*{prf7}{Proof of Theorem \ref{classify}}



\title{Heegaard Floer homology, L-spaces and smoothing order on links II}
\author{TAKUYA USUI}

\begin{document}

\begin{abstract}
We focus on L-spaces for which the boundary maps of the Heegaard Floer chain complexes vanish. In previous paper \cite{Usui}, we collect such manifolds systematically by using the smoothing order on links. In this paper, we classify such L-spaces under appropreate constraint. 
\end{abstract}

\address{Graduate School of Mathematical Science, University of Tokyo, 3-8-1 Komaba Meguroku Tokyo 153-8914, Japan}
\email{usui@ms.u-tokyo.ac.jp}
\maketitle

%
%一節
\section{Introduction}%セクション%

\ In \cite{OS2} and \cite{OS1}, Ozsv\'{a}th and Szab\'{o} introduced the \textit{Heegaard-Floer homology} $\widehat{HF}(Y)$ for a closed oriented three manifold $Y$. The Heegaard Floer homology $\widehat{HF}(Y)$ is defined by using a pointed Heegaard diagram representing $Y$ and a certain version of Lagrangian Floer theory. The boundary map of the chain complex \textit{counts} the number of pseudo-holomorphic Whitney disks. Of course, the boundary map depends on the pointed Heegaard diagram. In this paper, the coefficient of homology is $\mathbb{Z}_{2}$. A rational homology three-sphere $Y$ is called an L-space when its Heegaard Floer homology $\widehat{HF}(Y)$ is a $\mathbb{Z}_{2}$-vector space with dimension $|H_{1}(Y;\mathbb{Z})|$, where $|H_{1}(Y;\mathbb{Z})|$ is the number of elements in $H_{1}(Y;\mathbb{Z})$. 

%Note that any lens-space is L-space. Actually we can draw a genus one Heegaard diagram representing $L(p,q)$ for which the two circles $alpha$ and $\beta$ meet transversely in $p$ points.

In this paper, we consider a special class of L-spaces. 
\begin{defn}
An L-space $Y$ is \textit{strong} if there is a pointed Heegaard diagram representing $Y$ such that the boundary map vanishes.
\end{defn}

Strong L-spaces are originally defined in \cite{LL} in another way (see Proposition \ref{equi}), and discussed in \cite{BGW} and \cite{Greene}.

Now, We prepare some notations to state the main theorems.

\begin{figure}[h]

\includegraphics[width=7cm,clip]{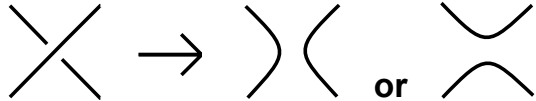} \\
\caption{smoothing}\label{smoothing}
\end{figure}

For a link $L$ in $S^3$, we can get a link diagram $D_{L}$ in $S^2$ by projecting $L$ to $S^2 \subset S^3$. To make other link diagrams from $D_{L}$, we can smooth a crossing point in different two ways (see Figure \ref{smoothing}.)

In \cite{EIT} and \cite{Taniyama}, the following ordering on links is defined.
\begin{defn}

Let $D_{L_{1}}$ and $D_{L_{2}}$ be alternating link diagrams in $S^2$. 
We say $D_{L_{1}} \subseteq D_{L_{2}}$ if $D_{L_{2}}$ contains $D_{L_{1}}$ as a connected component after smoothing some crossing points of $D_{L_{2}}$.

Let $L_{1}$ and $L_{2}$ be alternating links in $S^3$. Then, we say $L_{1} \leq L_{2}$ if for any minimal crossing alternating link diagram $D_{L_{2}}$ of $L_{2}$, there is a minimal crossing alternating link diagram $D_{L_{1}}$ of $L_{1}$ such that $D_{L_{1}} \subseteq D_{L_{2}}$. 
\end{defn}

\begin{figure}[h]

\includegraphics[width=4cm,clip]{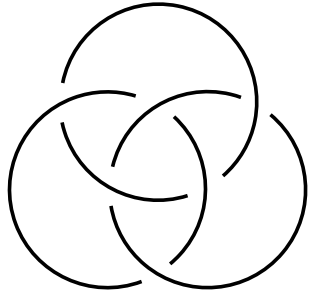} \\
\caption{The Borromean rings}\label{borromean}
\end{figure}

These orderings on links and diagrams are called \textit{smoothing orders} in \cite{EIT}. Note that smoothing orders become partial orderings. Let us denote the minimal crossing number of $L$ by $c(L)$. If $L_{1} \leq L_{2}$, then $c(L_{1}) \leq c(L_{2})$. We can check the well-definedness by using this observation. Actually, if $L_{1} \leq L_{2}$ and $L_{2} \leq L_{1}$, then $c(L_{1}) = c(L_{2})$ and there is no smoothed crossing point. So $L_{1} = L_{2}$. Next, if $L_{1} \leq L_{2}$ and $L_{2} \leq L_{3}$, then $L_{1} \leq L_{3}$ by defintion.
Note that we can define $\leq$ for any two links by ignoring alternating conditions. But in this paper we consider only alternating links and alternating link diagrams. The Borromean rings $\rm{Brm}$ are an alternating link in $S^3$ whose diagram looks  as in Figure \ref{borromean}. We fix this diagram and denote it by $\rm{Brm}$ too.

\begin{defn}
$\mathcal{L}_{\overline{\rm{Brm}}} = \{$ an alternating link $L$ in $S^3$ such that $\rm{Brm} \nleq L\}$, where $\rm{Brm}$ is the Borromean rings.
\end{defn}

Denote $\Sigma(L)$ a double branched covering of $S^3$ branched along a link $L$.
The following theorem is proved in \cite{Usui}:

\begin{thm} \label{thmMain}
Let $L$ be a link in $S^3$. If $L$ satisfies the following conditions:
\begin{itemize}
\item $L \in \mathcal{L}_{\overline{\rm{Brm}}}$, 
\item $\Sigma(L)$ is a rational homology three-sphere,
\end{itemize}
then $\Sigma(L)$ is a strong L-space and a graph manifold (or a connected sum of graphmanifolds).
\end{thm}

A graph manifold is defined as follows.
\begin{defn}
A closed oriented three manifold $Y$ is a \textit{graph manifold} if $Y$ can be decomposed along embedded tori into finitely many Seifert manifolds.
\end{defn}

In \cite{Usui}, we also defined the following class of manifolds.

\begin{defn}
$(T,\sigma,w)$ is an \textit{alternatingly weighted tree} when the following three conditions hold.
\begin{itemize}
\item $T$ is a disjoint union of trees (i.e a disjoint union of simply connected, connected graphs). 
Let $V(T)$ denote the set of all vertices of $T$.
\item $\sigma:V(T) \rightarrow \{\pm 1\}$ is a map such that if two vertices $v_{1}$, $v_{2}$ are connected by an edge, then $\sigma(v_{1})=-\sigma(v_{2})$. 
\item $w:V(T) \rightarrow \mathbb{Q}_{\ge 0} \cup \{\infty\}$ is a map.
\end{itemize}
\end{defn}  

\begin{figure}[h]

\includegraphics[width=9cm,clip]{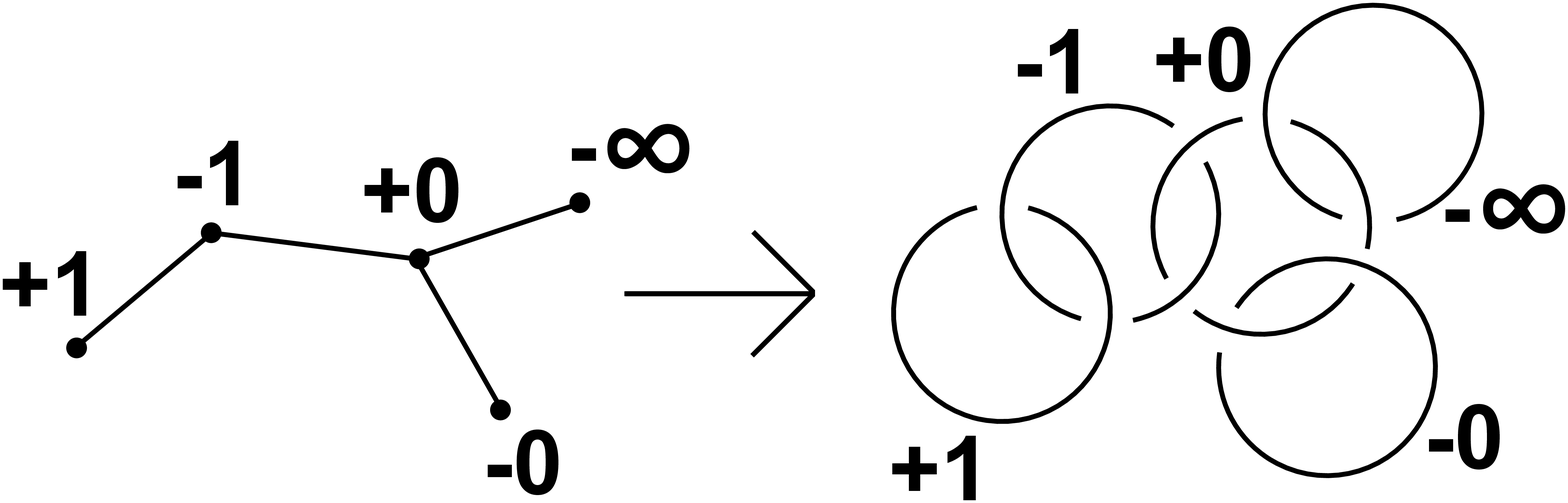} \\
\caption{}\label{tree1}
\end{figure}

Denote $\mathcal{T}$ the set of all alternatingly weighted trees.
For an alternatingly weighted tree $(T,\sigma,w)$, shortly $T$, we  define a three manifold $Y_{T}$ as follows.
First, we can take a realization of the tree $T$ in $\mathbb{R}^2 \subset S^3$. For each vertex $v$, we introduce the unknot in $S^3$. Next if two verteces in $T$ are connected by an edge, we link the corresponding two unknots with linking number $\pm 1$. Thus, we get a link $L_{T}$ in $S^3$. Then, we can get a new closed oriented three-manifold $Y_{T}$ by the surgery of $S^3$ along every unknot component of $L_{T}$ with the surgery coefficients $\sigma(v) w(v)$ (see Figure \ref{tree1})

This process gives a natural map $\mathcal{T} \rightarrow \mathcal{M}_{\mathcal{T}} = \{Y_{T} ; T \in \mathcal{T}\}/\text{homeo}$.

%\begin{thm}
%$\mathcal{M}_{\mathcal{AA}}=\mathcal{M}_{\mathcal{T}}$.

%\end{thm}
However, we can prove that these classes are equivalent as follows (see \cite{Usui}). 

\begin{thm} \label{thmcon}
The set of the three manifolds $Y_{T}$ induced from alternatingly weighted trees $(T,\sigma,w)$ is equal to the set of the branched double coverings $\Sigma(L)$ of $S^{3}$ branced along $L$ with $\rm{Brm} \nleq L\}$. That is,
$\mathcal{M}_{\mathcal{T}} = \mathcal{L}_{\overline{\rm{Brm}}}$.
\end{thm}

\begin{defn}
A Heegaard diagram is called \textit{strong} if the induced boundary map of Heegaard Floer chain complex vanishes. 
\end{defn}

Then, we prove the following classification theorem.

\begin{thm}\label{classify}
Let $(\Sigma,\alpha,\beta,z)$ be a strong Heegaard diagram representing a strong $L$-space $Y$. If the genus of $\Sigma$ is at most three, then $Y$ is in $\mathcal{L}_{\overline{\rm{Brm}}} = \mathcal{M}_{\mathcal{T}}$.
\end{thm}

\section{Heegaard-Floer homology and L-spaces}

\ The Heegaard Floer homology of a closed oriented three manifold $Y$ is defined from a pointed Heegaard diagram representing $Y$.
Let $f$ be a self-indexing Morse function on $Y$ with $1$ index zero critical point and $1$ index three critical point. Then, $f$ gives a Heegaard splitting of $Y$. That is, $Y$ is given by glueing two handlebodies $f^{-1}([0,3/2])$ and $f^{-1}([3/2, 3])$ along their boundaries. If the number of index one critical points or the number of index two critical points of $f$ is $g$, then $\Sigma = f^{-1}(3/2)$ is a closed oriented genus $g$ surface. 
We fix a gradient flow on $Y$ corresponding to $f$. We get a collection $\alpha = \{ \alpha_{1}, \cdots, \alpha_{g} \}$ of $\alpha$ curves on $\Sigma$ which flow down to the index one critical points, and another collection $\beta = \{ \beta_{1}, \cdots, \beta_{g} \}$ of $\beta$ curves on $\Sigma$ which flow up to the index two critical points.
Let $z$ be a point in $\Sigma \setminus (\alpha \cup \beta)$.
The tuple $(\Sigma,\alpha,\beta,z)$ is called a \textit{pointed Heegaard diagram} for $Y$.
Note that $\alpha$ and $\beta$ curves are characterized as pairwise disjoint, homologically linearly independent, simple closed curves on $\Sigma$. We can assume $\alpha$-curves intersect $\beta$-curves transversaly.%A Heegaard diagram specifies a three-manifold $Y$ as follows. (書け). 

%Symg(Sigma)の定義とWhitney-diskの定義
Next, we review the definition of the Heegaard Floer chain complex.

Let $(\Sigma,\alpha,\beta,z)$ be a pointed Heegaard diagram for $Y$.
The $g$-fold \textit{symmetric product} of the closed oriented surface $\Sigma$ is defined by $\rm{Sym}^{g}(\Sigma) = \Sigma^{\times g}/S_{g}$. That is, the quotient of $\Sigma^{\times g}$ by the natural action of the symmetric group on $g$ letters.

Let us define $\mathbb{T}_{\alpha } = \alpha_{1} \times \cdots \times \alpha_{g} /S_{g}$ and $\mathbb{T}_{\beta} = \beta_{1} \times \cdots \times \beta_{g}/S_{g}$.

%Let $(\Sigma,\alpha,\beta,z)$ be a pointed Heegaard diagram representing $Y$.
Then, the chain complex $\widehat{CF}(\Sigma,\alpha,\beta,z)$ is defined as a free $\mathbb{Z}_{2}$-module generated by the elements of 
$$\mathbb{T}_{\alpha } \cap \mathbb{T}_{\beta }= \{ x=(x_{1 \sigma (1)},x_{2 \sigma (2)},\cdots , x_{g \sigma (g)}) | x_{i \sigma (i)} \in \alpha_{i} \cap \beta_{\sigma(i)}, \sigma \in S_{g} \}.$$ 
%where $S_{g}$ is the symmetric group on $g$ letters. 

Then, the boundary map $\widehat{\partial}$ is given by%T\alphaは書かなくてもいいかも

\begin{equation}\label{def}
\widehat{\partial}x = \sum_{y \in \mathbb{T}_{\alpha} \cap \mathbb{T}_{\beta}} c(x,y) \cdot y,
\end{equation}
where $c(x,y) \in \mathbb{Z}_{2}$ is defined by \textit{counting} the number of pseudo-holomorphic Whitney disks. For more details, see \cite{OS2}.

\begin{defn}\cite{OS2}
The homology of the chain complex $(\widehat{CF}(\Sigma,\alpha,\beta,z),\widehat{\partial})$ is called
the Heegaard Floer homology of a pointed Heegaard diagram. We denote it by $\widehat{HF}(\Sigma,\alpha,\beta,z)$.
\end{defn}

\begin{rem}
%There are some more remarks.
%\begin{itemize}

%\item 
For appropriate pointed Heegaard diagrams representing $Y$, their Heegaard Floer homologies become isomorphic. So we can define the Heegaard Floer homology of $Y$. Denote it by $\widehat{HF}(Y)$. (For more details, see \cite{OS2}).

\end{rem}
%次を随所に埋め込む

In this paper, we consider only L-spaces, in particular strong L-spaces. The following proposition enables us to define strong L-spaces in another way. The second condition comes from \cite{LL}.

\begin{prop}\label{equi}
Let $(\Sigma,\alpha,\beta,z)$ be a pointed Heegaard diagram representing a rational homology sphere $Y$. Then, the following two conditions (1) and (2) are equivalent. 
\begin{enumerate}
\item the boundary map $\widehat{\partial}$ is the zero map, and $Y$ is an $L$-apace.
\item $|\mathbb{T}_{\alpha} \cap \mathbb{T}_{\beta}|=|H_{1}(Y;\mathbb{Z})|$.
\end{enumerate}
\end{prop}

For example, any lens-spaces are strong L-spaces. Actually, we can draw a genus one Heegaard diagram representing $L(p,q)$ for which the two circles $\alpha$ and $\beta$ meet transversely in $p$ points. That is, $|\mathbb{T}_{\alpha} \cap \mathbb{T}_{\beta}|=|H_{1}(L(p,q);\mathbb{Z})| = p$.

%\subsection{Proof of Proposition \ref{equi}}\label{prfofprop}
%In this subsection, we prove Proposition \ref{equi}.

%lemの証明
To prove this proposition, we recall that the Heegaard Floer homology $\widehat{HF}(Y)$ admits a relative $\mathbb{Z}/ 2 \mathbb{Z}$ grading(\cite{OS1}) By using this grading, the Euler characteristic satisfies the following equation.
%Ozsv\'{a}th-Szab\'{o} proved the next lemmma.%next
$$\chi(\widehat{HF}(Y,s)) = |H_{1}(Y;\mathbb{Z})|.$$

%propの証明
\begin{prf}
The first condition tells us that $\widehat{CF}(Y)$ becomes a $\mathbb{Z}_{2}$-vector space with dimension $|H_{1}(Y;\mathbb{Z})|$. By definition of $\widehat{CF}(Y)$, we get that $|\mathbb{T}_{\alpha} \cap \mathbb{T}_{\beta}|=|H_{1}(Y;\mathbb{Z})|$.
Conversely, the second condition and the above equation tell us that both $\widehat{CF}(Y)$ and $\widehat{HF}(Y)$
become $\mathbb{Z}_{2}$-vector spaces with dimension $|H_{1}(Y;\mathbb{Z})|$. Therefore, the first condition follows.
\qed
\end{prf}

\section{Strong diagram and induced matrix}
\subsection{Characterization of strong diagram}

\begin{figure}[h]

\includegraphics[width=4cm,clip]{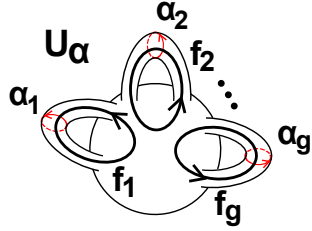} \\
\caption{A generator of $H_{1}(Y;\mathbb{Z})$}\label{fs}
\end{figure}

Let $(\Sigma,\alpha,\beta)$ be a Heegaard diagram. Fix orientations of $\alpha$-  and $\beta$-curves. 
Let $\{f_{1}, \cdots , f_{g}\}$ be pairwise disjoint simple closed oriented curves with $\# (f_{i} \cap \alpha_{j}) = + \delta_{ij}$ (see Figure \ref{fs}). Then, $\{f_{1}, \cdots , f_{g}\}$ generates $H_{1}(Y; \mathbb{Z})$. So $\beta_{j}$ can be written as linear combinations
$$\beta_{j} = \sum^{g}_{i = 1} a_{ij} f_{i}, \text { where } a_{ij} \in \mathbb{Z}.$$
In other words, $a_{ij}$ is defined as the algebraic intersection number of $\alpha_{i}$ and $\beta_{j}$, i.e., $a_{ij} = \#(\alpha_{i} \cap \beta_{j})$.
Let us define $A = (a_{ij})$.

Next, let $b_{ij}$ be the number of the elements of the set $\alpha_{i} \cap \beta_{j}$. $b_{ij}$ become non negative integers. Let $B = (b_{ij})$. By the definition of $a_{ij}$ and $b_{ij}$, it is easy to see that $|a_{ij}| \le b_{ij}$. Moreover, we get that $|\rm{det}(A)| = |H_{1}(Y; \mathbb{Z})|$ and $\overline{\rm{det}}(B) = \mathbb{T}_{\alpha} \cap \mathbb{T}_{\beta}$, where $$\overline{\rm{det}} (B) = \sum_{\sigma \in S_{g}} b_{1\sigma(1)}\cdots b_{g\sigma(g)}.$$ 

%For example, $\overline{det} (C)$ for 
%\[ C = \left (\begin{array}{@{\,}cccc@{\,}}
%-4 & -3 \\
%-1 & -2 \end{array} \right) \]
%is $11$.
If $(\Sigma,\alpha,\beta)$ is a strong Heegaard diagram, then $|\rm{det}(A)| = \overline{\rm{det}}(B)$ by its definition. Now we describe a characterization of strong Heegaard diagrams.

\begin{lem}\label{effec}
A Heegaard diagram $(\Sigma,\alpha,\beta)$ is strong if and only if $A$ is effective and $|a_{ij}| = b_{ij}$ for all $(i,j)$.
\end{lem}

\begin{prf}
Let $(\Sigma,\alpha,\beta)$ be a Heegaard diagram and define $A$ and $B$ as above.
Then, $|\rm{det}(A)|$ and $\overline{\rm{det}}(B)$ are expanded as follows.

$$|\rm{det}(A)| = \bigl| \sum_{\rm{sgn}(\sigma) = +1} a_{1\sigma(1)}\cdots a_{g\sigma(g)} - \sum_{\rm{sgn}(\sigma) = -1} a_{1\sigma(1)}\cdots a_{g\sigma(g)}\bigr|, \text{ and }$$
$$\overline{\rm{det}}(B) = \sum_{\rm{sgn}(\sigma) = +1} b_{1\sigma(1)}\cdots b_{g\sigma(g)} + \sum_{\rm{sgn}(\sigma) = -1} b_{1\sigma(1)}\cdots b_{g\sigma(g)}.$$

By using the triangle inequality, 
$$\overline{\rm{det}}(A) = \bigl| \sum_{\rm{sgn}(\sigma) = +1} a_{1\sigma(1)}\cdots a_{g\sigma(g)} - \sum_{\rm{sgn}(\sigma) = -1} a_{1\sigma(1)}\cdots a_{g\sigma(g)} \bigr| $$
$$\le \sum_{\rm{sgn}(\sigma) = +1} |a_{1\sigma(1)}\cdots a_{g\sigma(g)}| + \sum_{\rm{sgn}(\sigma) = -1} |a_{1\sigma(1)}\cdots a_{g\sigma(g)}|.$$

This inequation becomes the equality if and only if $A$ is effective.

On the other hand, the inequality $|a_{ij}| \le b_{ij}$ implies that 
$$|\rm{det}(B)| = \sum_{\rm{sgn}(\sigma) = +1} b_{1\sigma(1)}\cdots b_{g\sigma(g)} + \sum_{\rm{sgn}(\sigma) = -1} b_{1\sigma(1)}\cdots b_{g\sigma(g)} $$
$$\ge \sum_{\rm{sgn}(\sigma) = +1} |a_{1\sigma(1)}\cdots a_{g\sigma(g)}| + \sum_{\rm{sgn}(\sigma) = -1} |a_{1\sigma(1)}\cdots a_{g\sigma(g)}|.$$

This inequation becomes the equality if and only if $|a_{ij}| = b_{ij}$ for all $(i,j)$.

Thus, this lemma follows.

\qed
\end{prf}

\begin{defn}
Let $(\Sigma,\alpha,\beta)$ be a strong Heegaard diagram and orient $\alpha$- and $\beta$- curves. \textit{The induced matrix } $A_{(\Sigma,\alpha,\beta)}$ is defined by $A_{(\Sigma,\alpha,\beta)} = (\#(\alpha_{i} \cap \beta_{j}))_{ij}$.
\end{defn}

\subsection{Equivalence class of matrices}

\begin{defn}
Let $A$ be a $g \times g $ effective matrix. 
$A$ is \textit{ not maximal} if $A$ can be changed into a new effective matrix $A'$ by replacing some $0$-component of $A$ with $+1$ or $-1$.
\end{defn}

\begin{defn}
Let $A = (a_{ij})$ and $B = (b_{ij})$ be $g \times g$ matrices. 
We say $A \sim B$ if $A$ can be changed into a new matrix $A'$ so that each $(i,j)$-component of $A'$ has the same signature $\{+,-,0\}$ as the $(i,j)$-component of $B$ by some of the following operations.
\begin{itemize}
\item to permutate two rows or two columns,
\item to multiple a row or a column by $(-1)$, 
\item to transpose the matrix.
\end{itemize}
\end{defn}

Note that if an effective (or maximal) matrix $A$ is equivalent to $B$, then $B$ is also effective (or maximal).
Thus, we can define 
$$E_{g} = \{A : \text{ effective }\} / \sim, \text{ and }$$
$$ME_{g} = \{A : \text{ maximal and effective }\} / \sim \subset E_{g}.$$

\begin{defn}
Let $[A]$ and $[B]$ be two elements of $E_{g}$. 
We say $[A] \le [B]$ if we can change $A$ into $A'$ by replacing some zero components of $A$ by $+1$ or $-1$ so that the new class $[A']$ is equal to $[B]$.
\end{defn}

Note that this ordering tells us that an element of $ME_{g}$ is \textit{maximal} in $E_{g}$.  

For $g = 2 $, we can prove easily that there exists only one maximal effective matrix $A_{2}$, i.e., $ME_{g} = \{[A_{2}]\}$;
\[ A_{2} = \left (\begin{array}{@{\,}cccc@{\,}}
+ & + \\
- & + \end{array} \right), \]
where we consider only signatures of components of matrices.

\begin{lem}
For $g = 3$, there exists only one maximal effective matrix $A_{3}$, where 
\[ A_{3} = \left (\begin{array}{@{\,}cccc@{\,}}
+ & + & + \\
- & + & + \\
0 & - & + \end{array} \right), \] i.e., $ME_{g} = \{[A_{3}]\}$
\end{lem}

\begin{prf}
It is easy to see that $A_{3}$ is a maximal effective matrix. 
Let $A$ be a maximal effective $3 \times 3$ matrix. We prove $A \sim A_{3}$.
First, put
\[ A = \left (\begin{array}{@{\,}cccc@{\,}}
x_{1} & x_{2} & x_{3} \\
x_{4} & x_{5} & x_{6} \\
x_{7} & x_{8} & x_{9} \end{array} \right). \]

Since $A$ is maximal, we can assume that the diagonal components of $A$ are all positive, i.e., $x_{1} = x_{5} = x_{9} = +$.
We can also assume that there exists at least one $0$ component, and we can put $x_{7} = 0$ (by permutating rows and columns.)
Next, we consider the signature of $x_{4}$.

\begin{itemize}
\item If $x_{4} = 0$, we can put $x_{6} = +$ and $x_{8} = -$ because $A$ is maximal. (by permutating rows and columns). Then, we can assign arbitrary signatures to $x_{2}$ and $x_{3}$. ($A$ is still effective.) By multipling second or third columns by $(-1)$, we can assume $x_{2} = x_{3} = +$. This operations may change the signature of $x_{5}$, $x_{6}$,$x_{8}$ and $x_{9}$. But by multipling second or third rows by $(-1)$ or permutating rows and columns, we can assume $x_{5} = x_{6} = x_{9} = +$ and $x_{8} = -$. Thus, all signatures of $A$ are decided, however $A$ is not maximal because we can take $x_{4}$ to be positive. This is contradiction.

\[ A \to \left (\begin{array}{@{\,}cccc@{\,}}
+     & x_{2} & x_{3} \\
x_{4} & +     & x_{6} \\
0     & x_{8} & +     \end{array} \right)
\to
\left (\begin{array}{@{\,}cccc@{\,}}
+     & x_{2} & x_{3} \\
0     & +     & x_{6} \\
0     & x_{8} & +     \end{array} \right)
\to
\left (\begin{array}{@{\,}cccc@{\,}}
+     & x_{2} & x_{3} \\
0     & +     & +     \\
0     & -     & +     \end{array} \right)
\to
\left (\begin{array}{@{\,}cccc@{\,}}
+     & +     & +     \\
0     & +     & +     \\
0     & -     & +     \end{array} \right) \].

\item If $x_{4} = +$, by multipling second row and column by $(-1)$, we can change so that $x_{4} = -$.
\item If $x_{4} = -$, then $x_{2} = +$ because $A$ is maximal. If $x_{8} = 0$, we can transpose $A$ and we return to the case when $x_{4} = 0$. So $x_{8} \neq 0$. Moreover, we can put $x_{8} = -$ by multipling third rows and columns by $(-1)$. Lastly, $x_{3}$ and $x_{6}$ must be positive because $A$ is maximal. Thus, $A \sim A_{3}$.
\[ A \to \left (\begin{array}{@{\,}cccc@{\,}}
+     & x_{2} & x_{3} \\
-     & +     & x_{6} \\
0     & x_{8} & +     \end{array} \right)
\to \left (\begin{array}{@{\,}cccc@{\,}}
+     & +     & x_{3} \\
-     & +     & x_{6} \\
0     & x_{8} & +     \end{array} \right)
\to
\left (\begin{array}{@{\,}cccc@{\,}}
+     & +     & x_{3} \\
-     & +     & x_{6} \\
0     & -     & +     \end{array} \right)
\to
\left (\begin{array}{@{\,}cccc@{\,}}
+     & +     & +     \\
-     & +     & +     \\
0     & -     & +     \end{array} \right)
= A_{3} \].
\end{itemize}
\qed
\end{prf}

Now we define $A'_{3}$ as follows.
\[ A'_{3} = \left (\begin{array}{@{\,}cccc@{\,}}
+ & 0 & + \\
- & + & 0 \\
0 & - & + \end{array} \right). \]

Note that $[A'_{3}] \le [A]$.

\section{In the case $g = 2$}\label{g2case}

\subsection{Types of strong diagrams with $g = 2$}

\begin{figure}[h]

\includegraphics[width=11cm,clip]{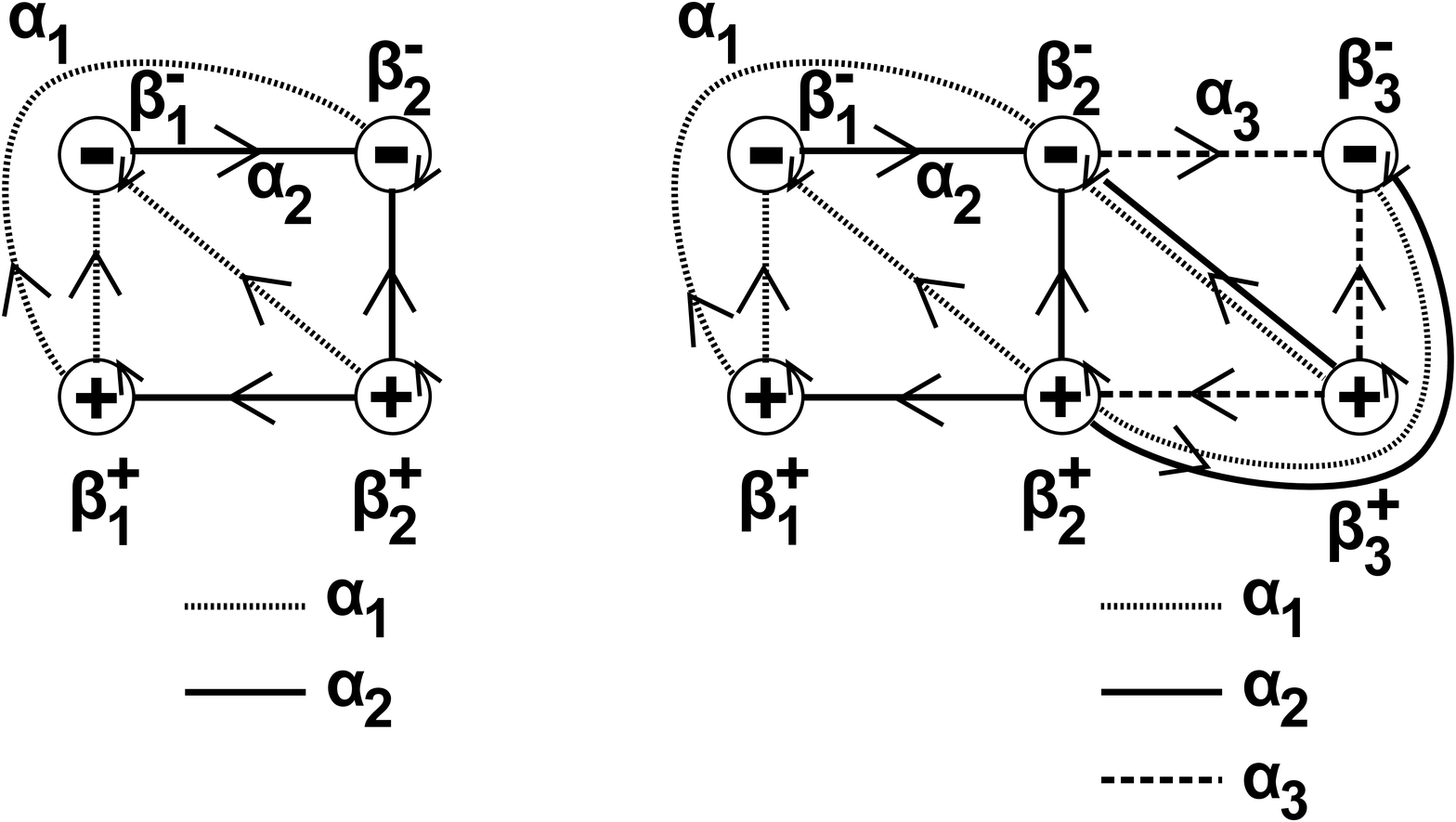} \\
\caption{examples of Heegaard diagrams ($g = 2, 3$)}\label{exhdiag}
\end{figure}

First, recall that we can describe a Heegaard diagram $(\Sigma,\alpha,\beta)$ in $\mathbb{R}^{2}$ as in Figure \ref{exhdiag}, where $2g$ signed oriented circles are $\beta$-curves and oriented arcs are $\alpha$ curves. By attaching the corresponding $\beta$ circles, we recover the Heegaard diagram $(\Sigma,\alpha,\beta)$.
Denote each $\beta$-circle by $\beta_{j}^{+}$ or $\beta_{j}^{-}$. 

Each $\alpha$-curve is divided into arcs by $\beta$-circles, and each arc is oriented and has its endpoints at $\beta_{j}^{\pm}$. So let us denote the set of $\alpha_{i}$-arcs from $\beta_{j}^{\pm}$ to $\beta_{j'}^{\pm}$ by $\Gamma_{i}(\pm j, \pm j')$. We put $\Gamma_{i}(\pm j', \pm j) = \Gamma_{i}(\pm j, \pm j')$ and $\Gamma(\pm j, \pm j') = \Gamma(\pm j', \pm j) = \amalg_{i} \Gamma_{i}(\pm j, \pm j')$.
Moreover, let us denote the union of all elements in $\Gamma(\pm j, \pm j')$ by $\Gamma'(\pm j, \pm j') \subset S^{2}$.
For example, the following $\alpha$-arc is an element in $\Gamma_{1}(+1, -2)$.

\begin{figure}[h]

\includegraphics[width=4cm,clip]{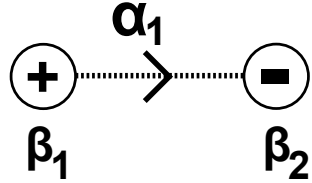} \\
\caption{an element in $\Gamma_{1}(+1, -2)$}\label{gamma}
\end{figure}

Let us consider the $g = 2$ case.

Let $(\Sigma,\alpha,\beta)$ be a strong diagram representing $Y$. By Lemma \ref{effec}, the induced matrix $A_{(\Sigma,\alpha,\beta)}$ is effective. Now we first consider the case when $A_{(\Sigma,\alpha,\beta)}$ is $A_{2}$. We call such a diagram \textit{an } $A_{2}$\textit{-strong diagram}. There exist just $8$ $\Gamma$-sets which may have some elements as follows.

$$\Gamma_{1}(+ \{1,2\}, - \{1,2\}),$$
$$\Gamma_{2}(- 1, + 1),\Gamma_{2}(+ 2, - 2),
\Gamma_{2}(+ 2,+ 1),
\Gamma_{2}(- 1, -2),$$

where $\{1,2\}$ means $1$ or $2$.

Now we prepare a notation.
Let $E$ a subset of $S^2$. For two elements $\gamma_{1}$ and $\gamma_{2}$ in $\Gamma(\pm j, \pm j')$, we say $\gamma_{1} \sim \gamma_{2}$ \textit{ out of } $E$ if $\gamma_{1}$ and $\gamma_{2}$ are isotopic on $(S^2 \setminus E, \beta_{j}^{\pm} \cup \beta_{j'}^{\pm})$.

\begin{figure}[h]

\includegraphics[width=4cm,clip]{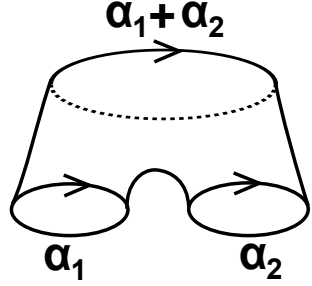} \\
\caption{Handle-slide}\label{hslide}
\end{figure}

Recall that we can transform a Heegaard diagram by isotopies and handle-slides and stabilizations.
For two $\alpha$- or $\beta$-curves, (for example $(\alpha_{1}, \alpha_{2})$, ) we can get a new pair of curves by adding one curve to the other curve (for example $(\alpha_{1}, \alpha_{2} \pm \alpha_{1})$) (see Figure \ref{hslide}). Denote this handle-slide by $\pm \alpha_{1} \leadsto \alpha_{2}$. In particular, for $\beta$-curves, we denote a handle-slide by $\beta_{j}^{\pm} \leadsto \beta_{j'}^{\pm}$.

\begin{prop}\label{2ste}
Let $(\Sigma,\alpha,\beta)$ be an $A_{2}$-strong diagram. Then,  $(\Sigma,\alpha,\beta)$ can be transformed by handle-slides and isotopies (if it is necessary) so that the new diagram is strong and $\Gamma(+ j, - j)$ has at least one element for each $j = 1,2$. (However, The new diagram may not be $A_{2}$-strong.)

%the one of the following conditions satisfies.
%\begin{itemize}
%\item the new diagram is also $A_{2}$-strong and $\Gamma(+ j, - j)$ has at least one element for each $j = 1,2$. 
%\item the new diagram is strong but not $A_{2}$-strong. 
%\end{itemize}
\end{prop}

\begin{prf}
We prove this proposition in two steps. 
\begin{enumerate}
\item We can transform $(\Sigma,\alpha,\beta)$ so that $\Gamma(+ 2, - 2) \neq \emptyset$.  
\item If $\Gamma(+ 2, - 2) \neq \emptyset$, we can transform the diagram so that $\Gamma(+ 1, - 1) \neq \emptyset$ (while keeping the condition $\Gamma(+ 2, - 2) \neq \emptyset$).
\end{enumerate}

In each step, we must take a strong diagram.

\underline{Step 1}
Let $(\Sigma,\alpha,\beta)$ be an $A_{2}$-strong diagram. 
If $\Gamma(+ 2, - 2) \neq \emptyset$, there is nothing to do. 

Suppose that $\Gamma(+ 2, - 2) = \emptyset$.
Since $\# (\alpha_{2} \cap \beta_{2}) \neq 0$, we get that $\Gamma(+ 2, + 1) \neq \emptyset$ and $\Gamma(- 1, - 2) \neq \emptyset$.
For any two element $\gamma_{1}$ and $\gamma_{2}$ in $\Gamma(+ 2, + 1)$, we get that $\gamma_{1} \sim \gamma_{2}$ out of $(\beta_{2}^{\mp} \cup \beta_{1}^{\mp} \cup \Gamma'(-1, -2))$ because $S^{2} \setminus (\beta_{2}^{\mp} \cup \beta_{1}^{\mp} \cup \Gamma'(-1, -2))$ consists of disjoint disks. 
Then, we can transform the diagram by a handle-slide $\beta_{2}^{+} \leadsto \beta_{1}^{+}$. Note that the new diagram is strong but not $A_{2}$-strong because $\Gamma_{2}(-1,-2) = \emptyset$.

\begin{figure}[h]

\includegraphics[width=7cm,clip]{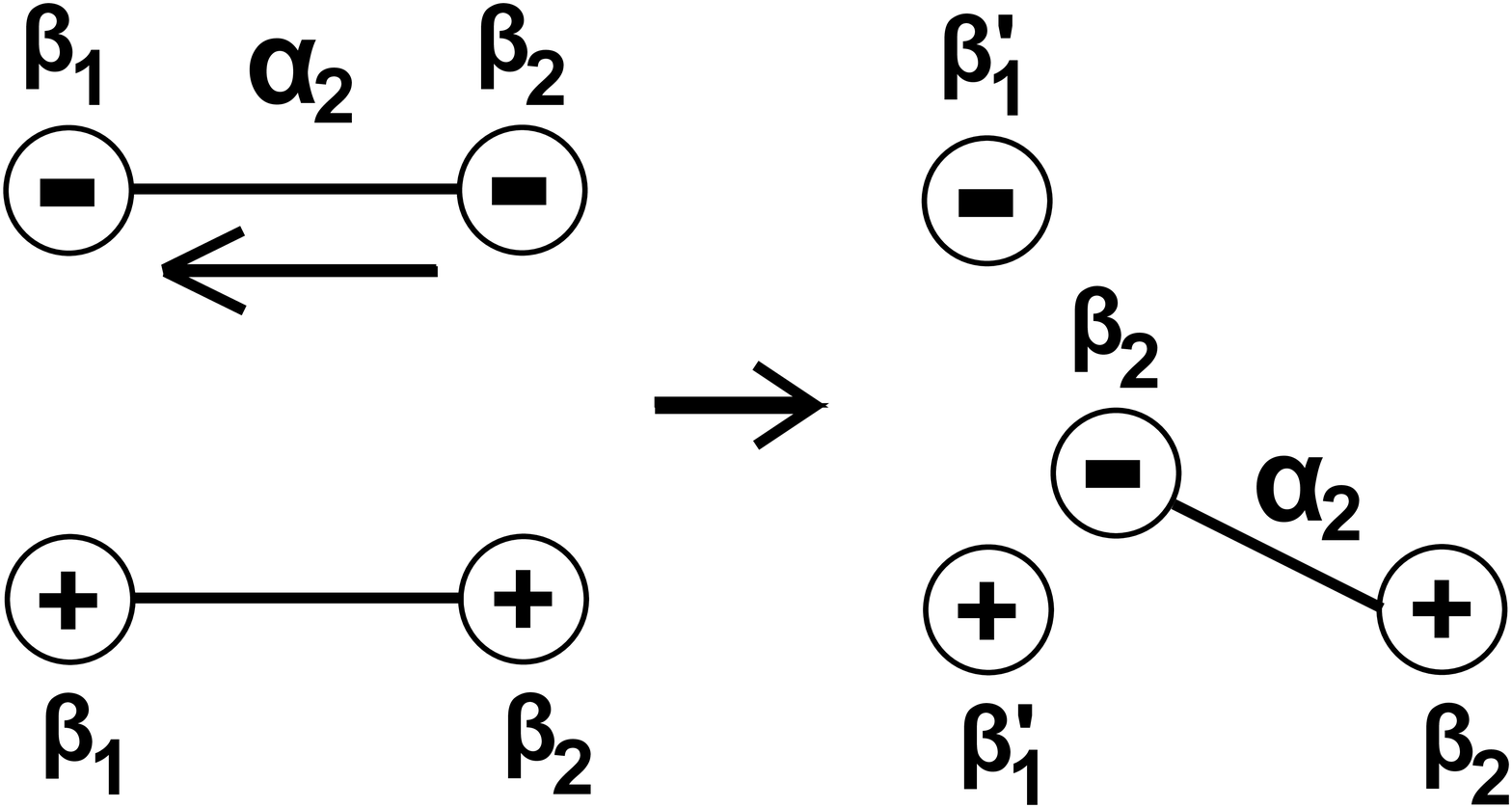} \\
\caption{}\label{selfcurve5-}
\end{figure}

\underline{Step 2}
Let $(\Sigma,\alpha,\beta)$ be an $A_{2}$-strong diagram with $\Gamma(+ 2, - 2) \neq \emptyset$.
If $\Gamma(+ 1, - 1) \neq \emptyset$, there is nothing to do. 

Suppose that $\Gamma(+ 1, - 1) = \emptyset$.
Since $\# (\alpha_{2} \cap \beta_{1}) \neq 0$, we get that $\Gamma(+ 2, + 1) \neq \emptyset$ and  $\Gamma(- 1, - 2) \neq \emptyset$.
For any two element $\gamma_{1}$ and $\gamma_{2}$ in $\Gamma(+ 2, + 1)$, $\gamma_{1} \sim \gamma_{2}$ out of $(\beta_{2}^{\mp} \cup \beta_{1}^{\mp} \cup \Gamma'(- 1, - 2))$, because $S^{2} \setminus (\beta_{2}^{\mp} \cup \beta_{1}^{\mp} \cup \Gamma'(- 1, - 2) \cup \Gamma'(+ 2, - 2))$ also consists of disjoint disks.
Then, we can transform the diagram by handle-slides $\beta_{1}^{+} \leadsto \beta_{2}^{+}$ finitely many times (see Figure \ref{selfcurve6-}). In finitely many steps, we will get an $A_{2}$-strong Heegaard diagram where $\Gamma(+ 1, - 1) \neq \emptyset$. We also get $\Gamma(+ 2, - 2) \neq \emptyset$ because the set $\alpha_{1} \cap \beta_{2}$ becomes non empty after these handle-slides.

\begin{figure}[h]

\includegraphics[width=9cm,clip]{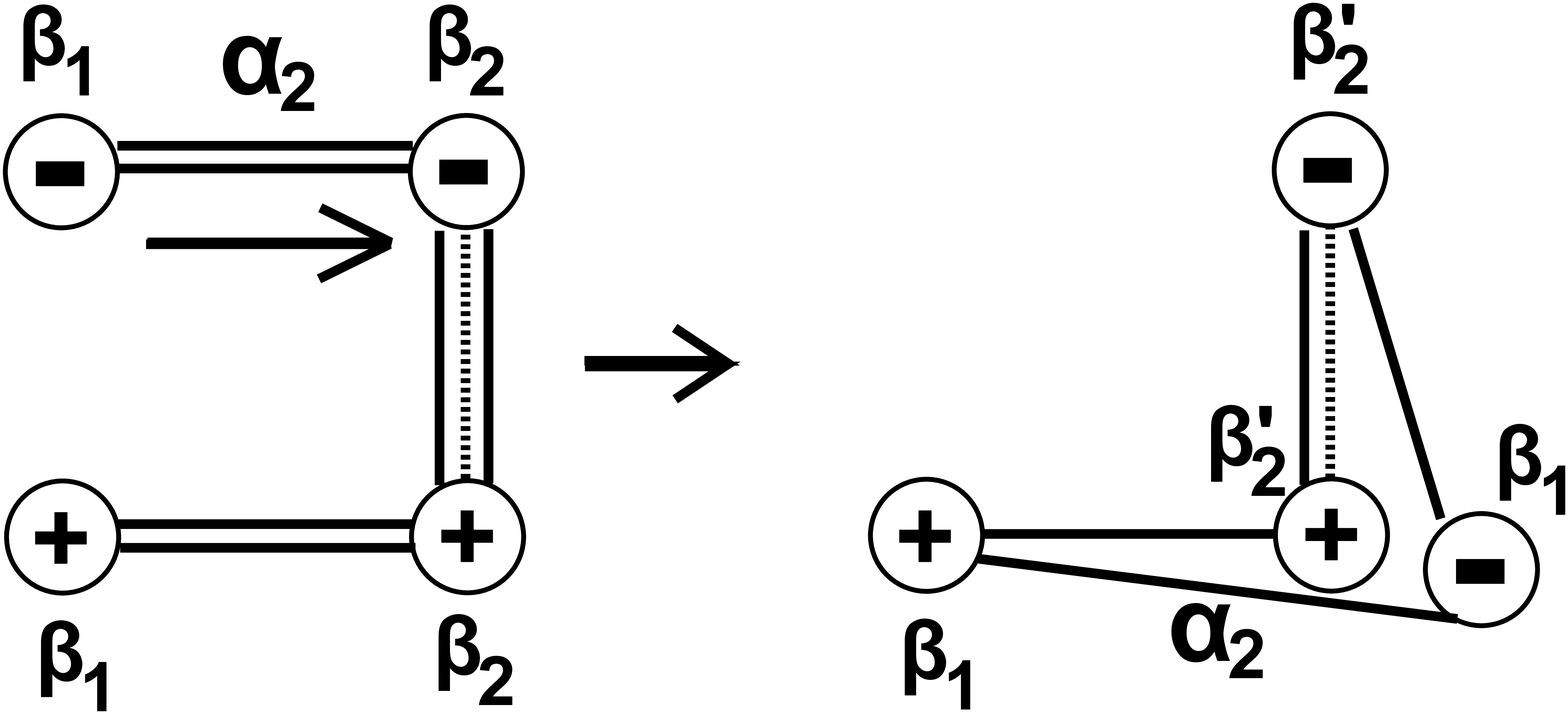} \\
\caption{}\label{selfcurve6-}
\end{figure}
\qed
\end{prf}

\begin{figure}[h]

\includegraphics[width=9cm,clip]{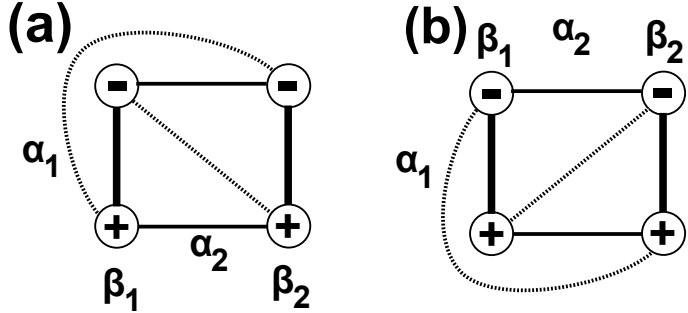} \\
\caption{type (a) and type (b)}\label{22case}
\end{figure}

\begin{prop}\label{aorb}
Let $(\Sigma,\alpha,\beta)$ be an $A_{2}$-strong diagram. Suppose that $\Gamma(+ j, - j) \neq \emptyset$ for each $j = 1,2$. 
Then, $(\Sigma,\alpha,\beta)$ is of type (a) or (b) shown in Figure \ref{22case}.
\end{prop}

\begin{prf}
It is easy to see that there are at most two possible types of diagrams under these assumptions. They are type (a) and type (b).
\qed
\end{prf}

Note that these two types of diagrams are equivalent under permutations of $\alpha$-curves and changes of the orientations of $\beta$-circles. Thus, we consider only the type (a).

\subsection{Surgery representations ($g = 2$)}

In this subsection, we give surgery representations of the diagrams.
To do this, we first study about $\Gamma(+j,-j)$ for $j = 1,2$ more precisely.
After that, we consider auxiliary attaching circles $\gamma = (\gamma_{1}, \gamma_{2})$.

Let $(\Sigma,\alpha,\beta)$ be an $A_{2}$-strong diagram of type (a) in Proposition \ref{aorb}.

\begin{figure}[h]

\includegraphics[width=6cm,clip]{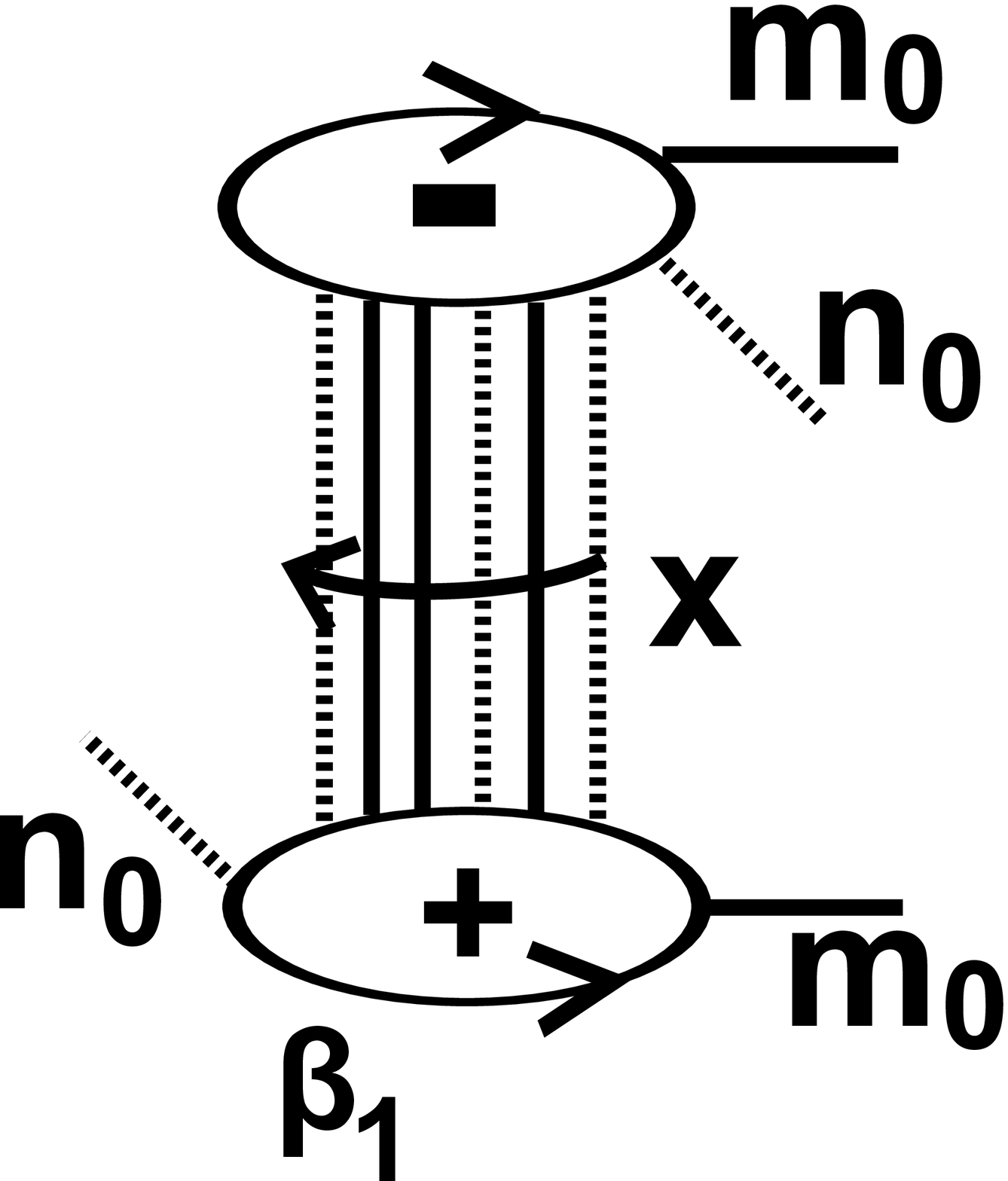} \\
\caption{}\label{nbd1}
\end{figure}

\underline{$\Gamma(+1,-1)$}:

The neighborhood of $\beta_{1}$-circle looks as in Figure \ref{nbd1}. 

Let $n_{0} = \#\Gamma(+1,-2) = \#\Gamma(+2,-1) > 0$ and $m_{0} = \#\Gamma(+2,+1) = \#\Gamma(-1,-2) > 0$.
Since $\Gamma(+1,-1) = \Gamma_{1}(+1,-1) \amalg \Gamma_{2}(-1,+1)$, we can define $x = (n_{1}, m_{1},n_{2},m_{2},\cdots, n_{k},m_{k})$ to be the sequence of the number of intersection points induced from $\Gamma(+1,-1)$. $n$ comes from  $\beta_{1} \cap \alpha_{1}$ and $m$ comes from $\beta_{1} \cap \alpha_{2}$. For example, $x = (n_{1}, m_{1},n_{2},m_{2},n_{3}) = (1,1,1,2,1)$ for Figure \ref{nbd1}. (By the way, this diagram in Figure \ref{nbd1} never become a Heegaard diagram because $m_{1} \neq m_{2}$.)

Since $\beta_{1}^{+}$ and $\beta_{1}^{-}$ are attached to be $\Sigma$, $x$ has the form $(n_{1}, m_{1},\cdots,m_{k-1}, n_{k})$, $(n_{1} \neq 0, n_{k} \neq 0)$ or $(m_{1},n_{1},\cdots, n_{k-1},m_{k})$, $(m_{1} \neq 0, m_{k} \neq 0)$. Moreover, one of the following two cases happen.

\begin{itemize}
\item $x = (n_{1}, m_{1},\cdots,m_{k-1}, n_{k})$, where $m_{i} = m_{0}$ for all $i$ and $n_{1} = n_{k}$,
\item $x = (m_{1},n_{1},\cdots, n_{k-1},m_{k})$, where $n_{i} = n_{0}$ for all $i$ and $m_{1} = m_{k}$.
\end{itemize}

\begin{figure}[h]

\includegraphics[width=7cm,clip]{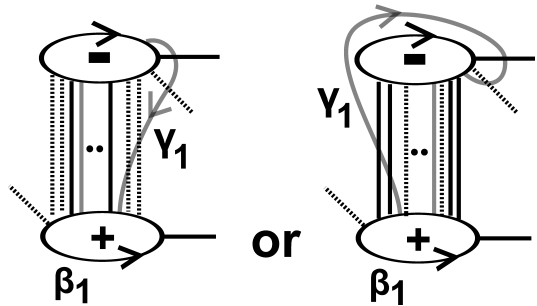} \\
\caption{$\gamma_{1}$-curve}\label{nbd4}
\end{figure}

In each case, take a simple closed oriented curve $\gamma_{1}$ on $\Sigma$ so that the following conditions hold (see Figure \ref{nbd4}). If $\Gamma(+1, -1)$ consists of only $\alpha_{1}$-arcs (or only $\alpha_{2}$-arcs), then we take $\gamma_{1} = \beta_{1}$. 
\begin{itemize}
\item $\gamma_{1}$ intersects each arc in $\Gamma(+2,-1)$ and $\Gamma(-1,-2)$ at a point. 
\item $\gamma_{1}$ intersects $\beta_{1}$ at some points, but does not intersect $\beta_{2}$.

\item If $x = (n_{1}, m_{1},\cdots,m_{k-1}, n_{k})$, $\gamma_{1}$ intersects only $\alpha_{1}$-arcs in $\Gamma(+1,-1)$.
\item If $x = (m_{1},n_{1},\cdots, n_{k-1},m_{k})$, $\gamma_{1}$ intersects only $\alpha_{2}$-arcs in $\Gamma(+1,-1)$.
\end{itemize}

\begin{figure}[h]

\includegraphics[width=7cm,clip]{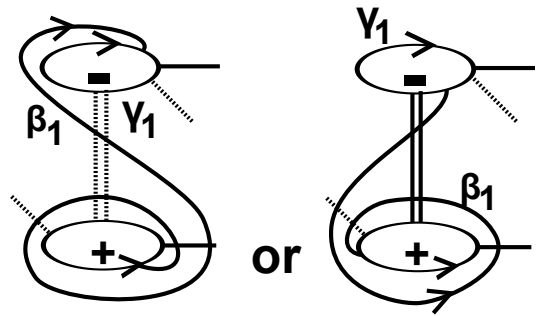} \\
\caption{new diagram $(\Sigma, \alpha, (\gamma_{1}, \beta_{2}))$}\label{nbd5}
\end{figure}

We consider the new $A_{2}$-strong Heegaard diagram $(\Sigma, \alpha, (\gamma_{1}, \beta_{2}))$. Then, the neighborhood of $\gamma_{1}$-circles looks as in Figure \ref{nbd5}. (In general, $\beta_{1}$ becomes arcs.) 
We describe $\beta_{1}$-curve more precisely later.


\begin{figure}[h]

\includegraphics[width=4cm,clip]{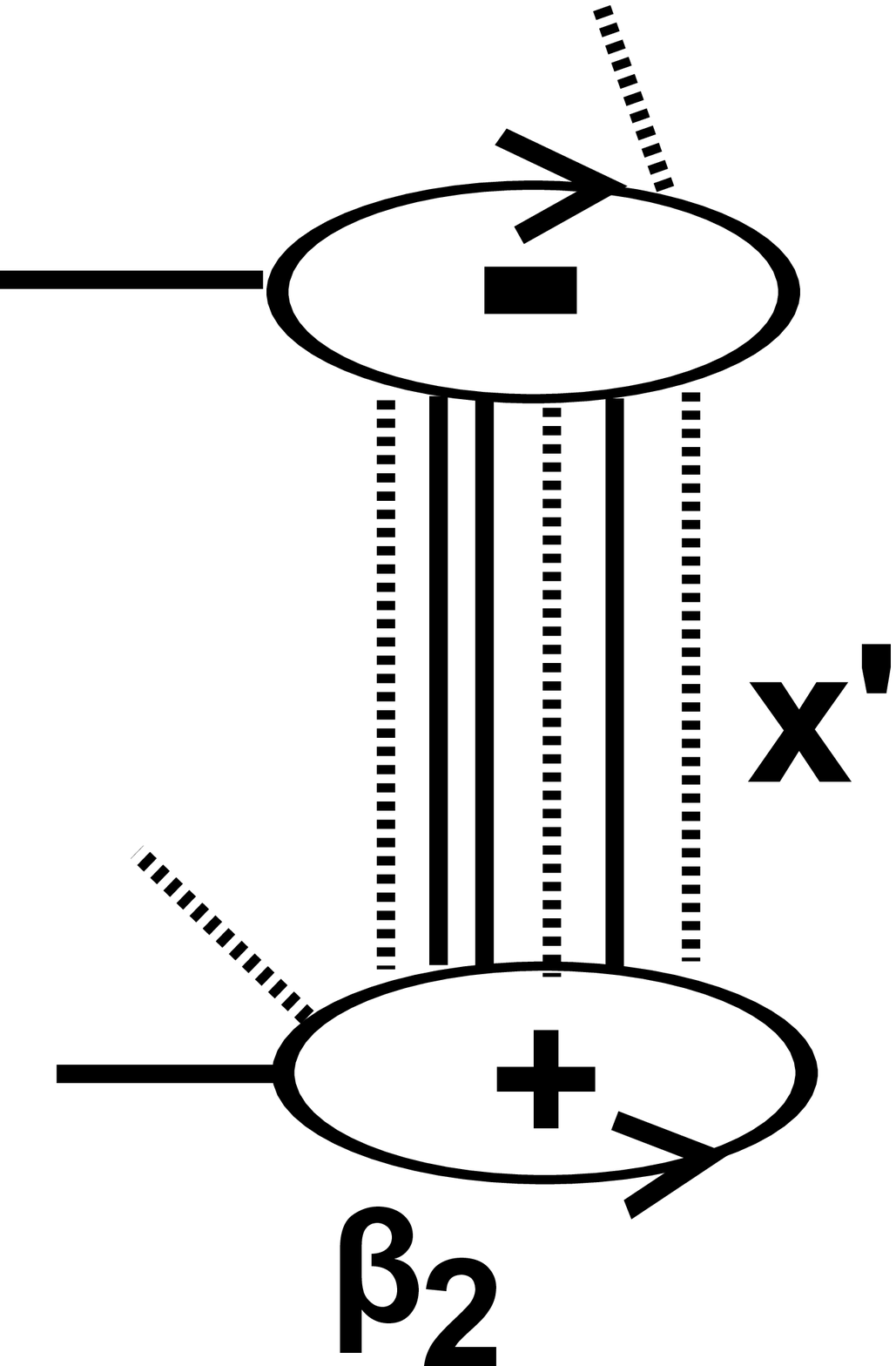} \\
\caption{}\label{nbd2-}
\end{figure}

\underline{$\Gamma(+2,-2)$}:

Next, consider the neighborhood of $\beta_{2}$ as in Figure \ref{nbd2-}.

Similarly, we can define $x'$ to be the the sequence of the number of intersection points induced from $\Gamma(+2,-2)$. There exist two cases, where $n'$ comes from $\beta_{2} \cap \alpha_{1}$ and $m$ comes from $\beta_{2} \cap \alpha_{2}$.

\begin{itemize}
\item $x' = (n'_{1}, m'_{1},\cdots,m'_{k-1}, n'_{k})$, where $m'_{i} = m'_{0}$ for all $i$ and $n'_{1} = n'_{k}$,
\item $x' = (m'_{1},n'_{1},\cdots, n'_{k-1},m'_{k})$, where $n'_{i} = n'_{0}$ for all $i$ and $m'_{1} = m'_{k}$,
\end{itemize}

\begin{figure}[h]

\includegraphics[width=7cm,clip]{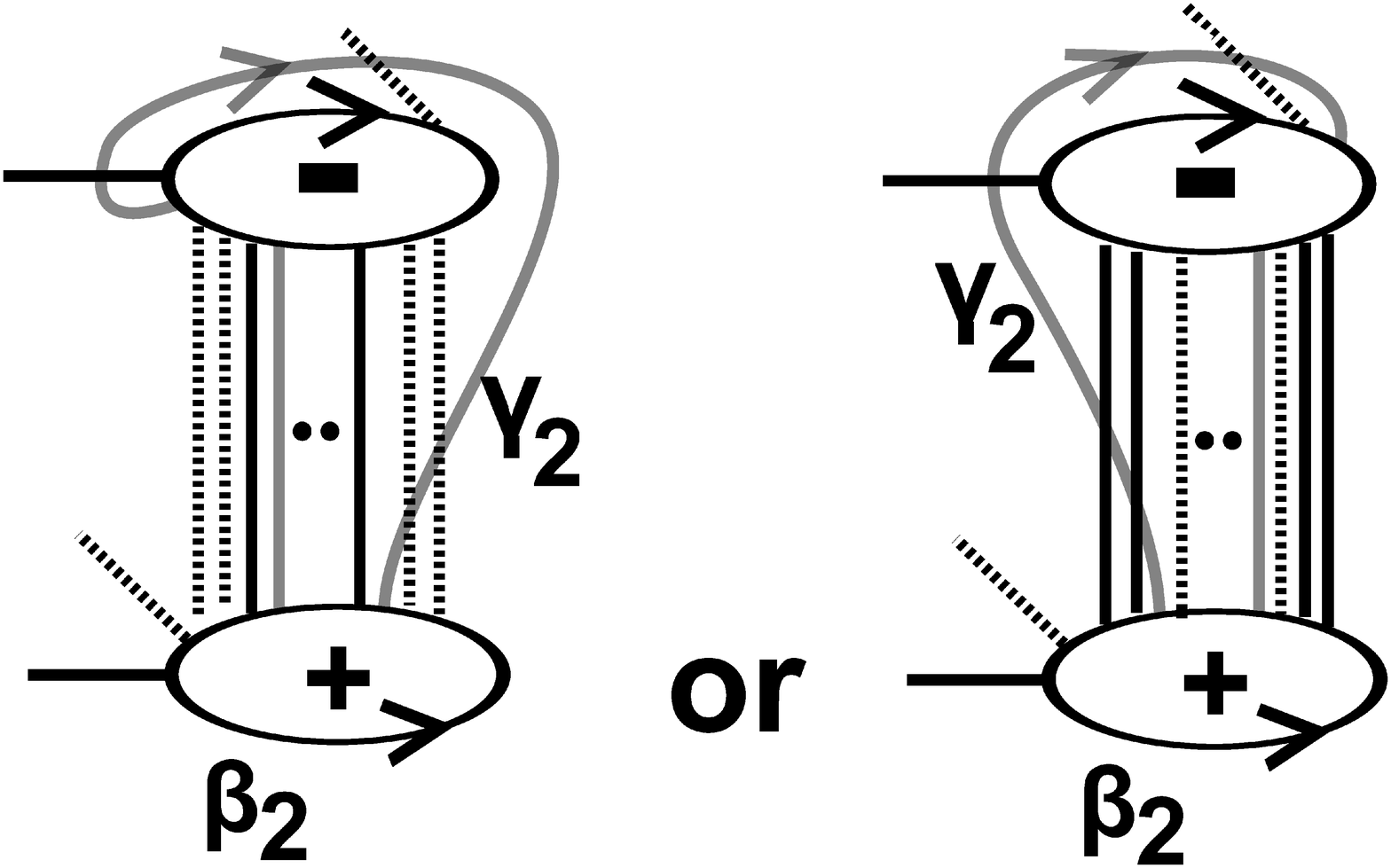} \\
\caption{$\gamma_{2}$-curve}\label{nbd4-}
\end{figure}

In each case, take a simple closed oriented curve $\gamma_{2}$ on $\Sigma$ similarly so that the following conditions hold (see Figure \ref{nbd4-}). If $\Gamma(+2, -2)$ consists of only $\alpha_{1}$-arcs (or only $\alpha_{2}$-arcs), then we take $\gamma_{2} = \beta_{2}$. 
\begin{itemize}
\item $\gamma_{2}$ intersects each arc in $\Gamma(+2,-1)$ and $\Gamma(-1,-2)$ at a point. 
\item $\gamma_{2}$ intersects $\beta_{2}$ at some points, but does not intersect $\beta_{1}$.

\item If $x' = (n'_{1}, m'_{1},\cdots,m'_{k-1}, n'_{k})$, $\gamma_{2}$ intersects only $\alpha_{1}$-arcs in $\Gamma(+2,-2)$.
\item If $x' = (m'_{1},n'_{1},\cdots, n'_{k-1},m'_{k})$, $\gamma_{2}$ intersects only $\alpha_{2}$-arcs in $\Gamma(+2,-2)$.
\end{itemize}

\begin{figure}[h]

\includegraphics[width=7cm,clip]{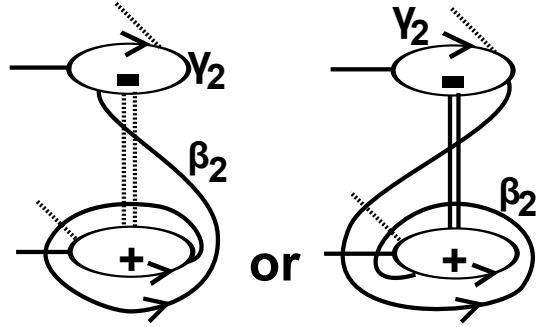} \\
\caption{new diagram $(\Sigma, \alpha, \gamma)$}\label{nbd5-}
\end{figure}

We consider the new $A_{2}$-strong Heegaard diagram $(\Sigma, \alpha, \gamma=(\gamma_{1}, \gamma_{2}))$. Then, the neighborhood of $\gamma_{2}$-circles looks as in Figure \ref{nbd5-}. (In general, $\beta_{2}$ becomes arcs.) 

We describe $\beta_{2}$-curve more precisely later.

\begin{figure}[h]

\includegraphics[width=10cm,clip]{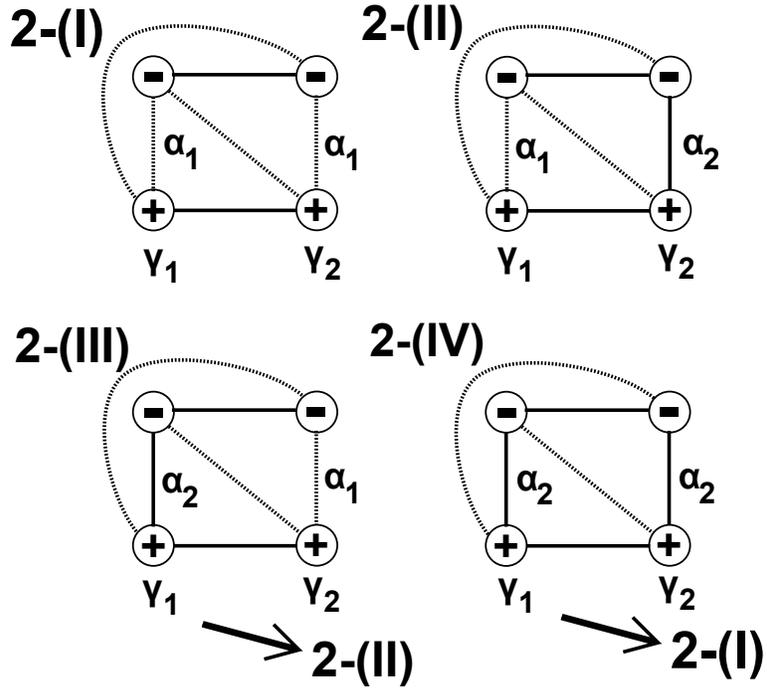} \\
\caption{possible types}\label{4case1-}
\end{figure}

Since the new diagram $(\Sigma, \alpha, \gamma)$ is easier than the old diagram, we classify them into four types (see Figure \ref{4case1-}). Actually, the new $\Gamma(+j,-j)$ consists of only $\alpha_{1}$-arcs or $\alpha_{2}$-arcs.

But, the type 2-(III) (resp. 2-(IV)) are equivalent to the type 2-(II) (resp. 2-(I)) under permutations of $\alpha$-curves and changes of the orientations of $\beta$-curves (before taking $\gamma$-curves).

\underline{type 2-(I)}
In this case, we get that $\#\Gamma(-1,-2) = \#\Gamma(+2,+1) = 1$. 
Thus, we can take another attaching circles $\delta = (\delta_{1},\delta_{2})$ in this diagram such that
\begin{itemize}
\item $\#(\delta_{i} \cap \gamma_{j}) = \delta_{ij}$ for any $(i,j)$ (thus, $(\Sigma,\delta,\gamma)$ represents $S^3$),
\item $\delta_{1}$ and  $\delta_{2}$  intersects $\alpha_{1}$ and does not intersect $\alpha_{2}$.
\end{itemize}

\begin{figure}[h]

\includegraphics[width=5cm,clip]{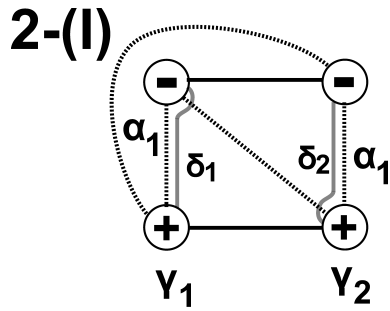} \\
\caption{$\delta = (\delta_{1},\delta_{2})$}\label{delta1}
\end{figure}

In this new diagram $(\Sigma,\delta,\gamma)$, $\alpha$ and $\beta$-curves become some framings of some knots in $U_{\delta} \cup U_{\gamma}$. Precisely, we can take three unknots $K_{1}$, $K_{2}$ and $C_{1}$. $K_{1}$ and $K_{2}$ are two unknots in $U_{\gamma}$ whose framings are $\beta$-curves. $C_{1}$ is an unknot in $U_{\delta}$ whose framing is $\alpha_{1}$-curve. Note that $\alpha_{2}$ can be written by $\delta$-curves as a homology in $\Sigma$, so there is no need to consider. These slopes can be determined as follows.

Let $r_{\alpha_{1}}$, $r_{\beta_{1}}$ and $r_{\beta_{2}}$ be the rational numbers representing the surgery framings of $\alpha_{1}$, $\beta_{1}$ and $\beta_{2}$ respectively. Precisely, these rational numbers are determined as follows. Put $r_{\alpha_{1}} = \rm{sgn}(r_{\alpha_{1}}) p_{1}/q_{1}$ and $r_{\beta_{i}} = \rm{sgn}(r_{\beta_{i}}) p'_{i}/q'_{i}$ for $i = 1,2$. Then,

\begin{itemize}
\item $\rm{sgn}(r_{\alpha_{1}}) = \rm{sgn}(r_{\beta_{1}}) = \rm{sgn}(r_{\beta_{2}}) = +1$,
\item $p_{1} = \#(\alpha_{1} \cap (\gamma_{1} \cup \gamma_{2})$,
\item $q_{1} = \#(\Gamma(+1,-2))$,
\item $p'_{i} = \#(\beta_{i} \cap \alpha_{2})$, for $i = 1,2$. 
\item $q'_{i} = \#(\beta_{i} \cap \gamma_{i})$, for $i = 1,2$. 
\end{itemize}

Note that $p_{1} > 2q_{1}$ and $p'_{i} > q'_{i}$ for $i = 1,2$.
Thus, $|r_{\alpha_{1}}| > 2$ and $|r_{\beta_{i}}| > 1$ for $i = 1,2$. (If $\gamma_{i} = \beta_{i}$, take $r_{\beta_{i}} = \infty$.)

\begin{figure}[h]

\includegraphics[width=10cm,clip]{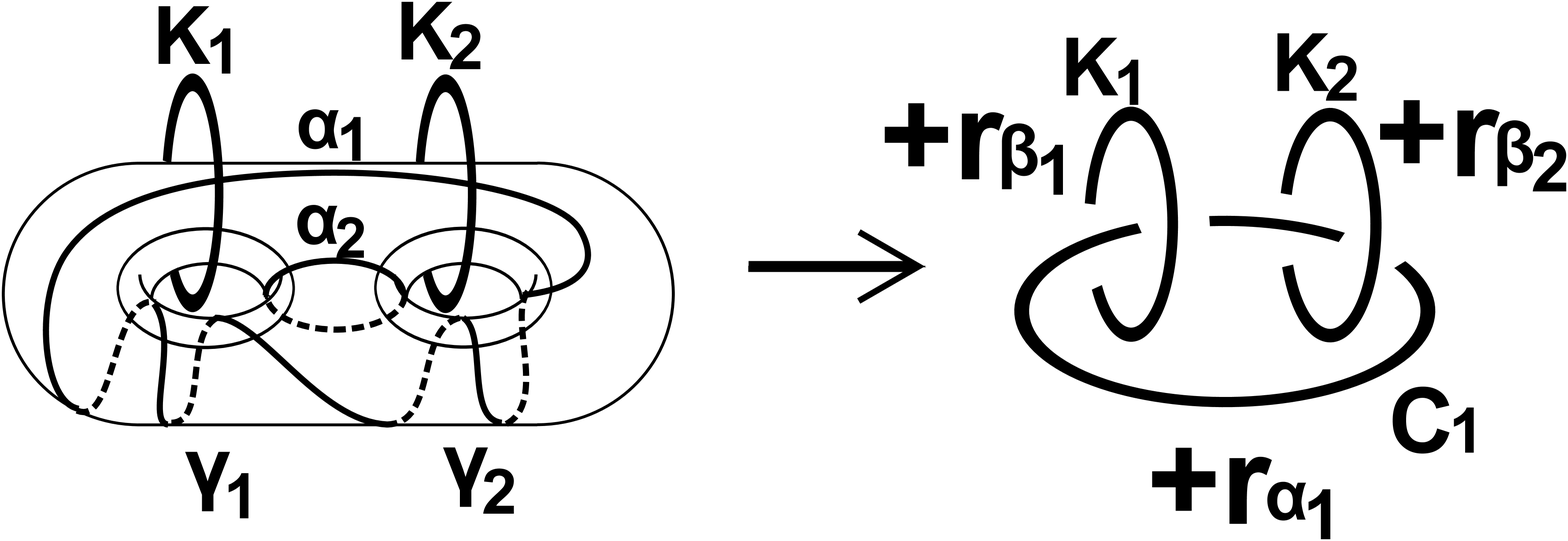} \\
\caption{surgery representation of $Y$}\label{k1k2-1}
\end{figure}

As a result, the three manifold $Y$ obtained from $(\Sigma, \alpha, \beta)$ can be represented as $Y = S^{3}(K_{1},K_{2},C_{1})$ (see Figure \ref{k1k2-1}).

\begin{figure}[h]

\includegraphics[width=6cm,clip]{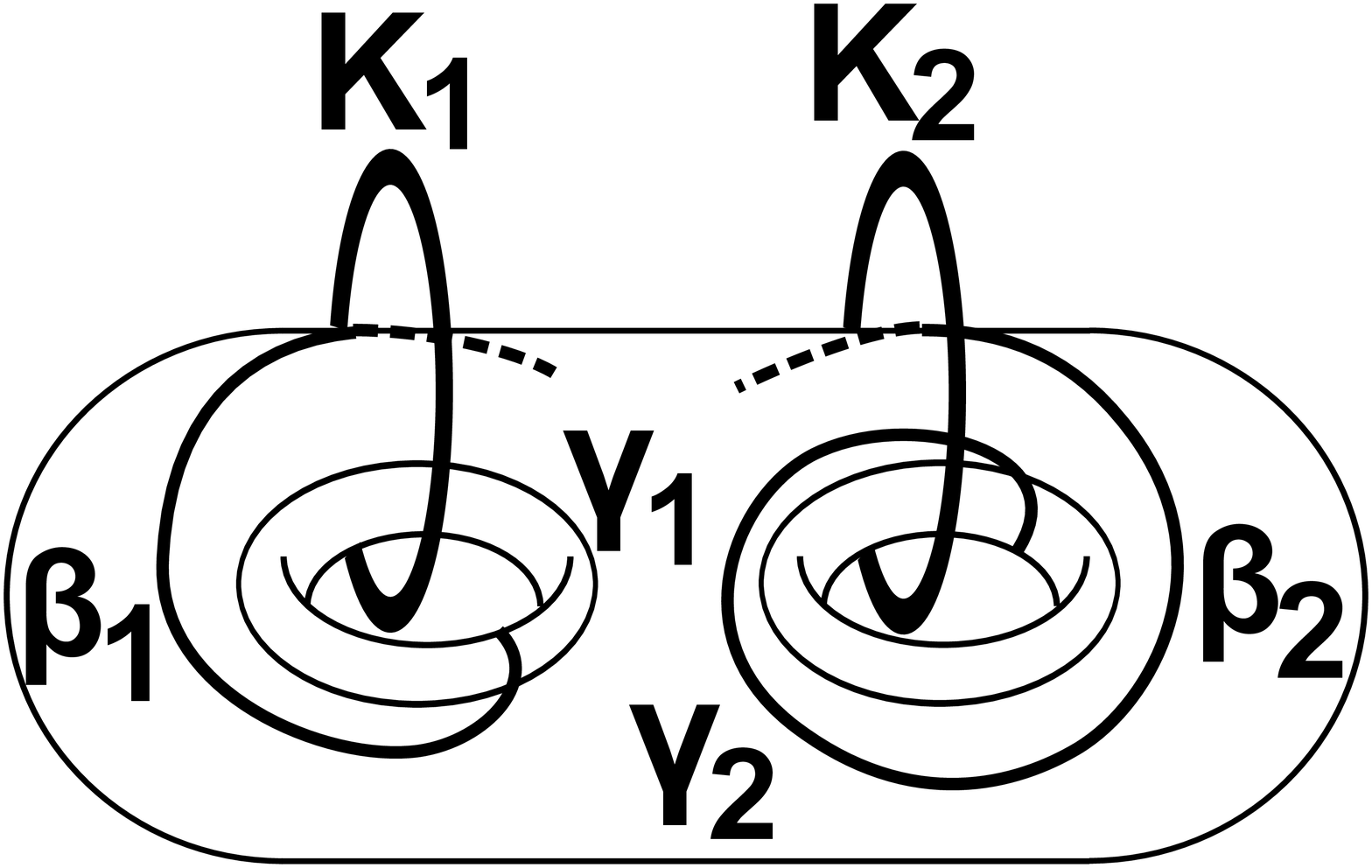} \\
\caption{}\label{k1k2}
\end{figure}

\begin{figure}[h]

\includegraphics[width=7cm,clip]{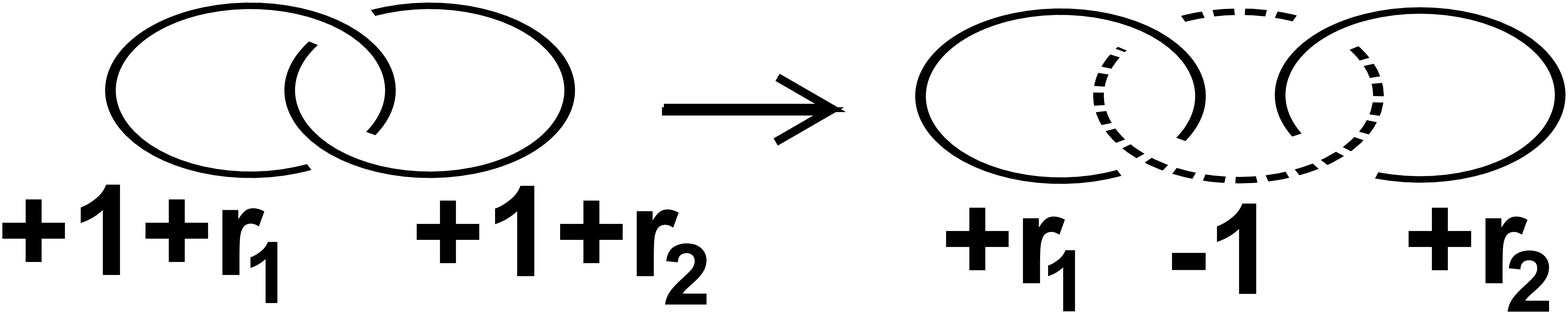} \\
\caption{}\label{slam1}
\end{figure}

It is easy to see that $Y$ belongs to $\mathcal{M}_{\mathcal{T}} = \mathcal{L}_{\overline{\rm{Brm}}}$. 
Actually, we can use the following Kirby calculus (see Figure \ref{slam1} and \cite{Gompf}). If two unknots with the linking number $\pm 1$ have rational framings $+1 +r_{1} \ge 1$ and $+1 +r_{2} \ge 1$, then we can perform the blow up operation so that the new link has alternating framings. Since $+ r_{\beta_{1}} > 1$, $+ r_{\beta_{2}} > 1$ and $+r_{\alpha_{1}} >2$, we get an alternatingly weighted link. Therefore, $Y$ is in $\mathcal{M}_{\mathcal{T}} = \mathcal{L}_{\overline{\rm{Brm}}}$.

\underline{type 2-(II)}

In this case, we can take another attaching circles $\delta = (\delta_{1},\delta_{2})$ in this diagram similarly such that
\begin{itemize}
\item $\#(\delta_{i} \cap \gamma_{j}) = \delta_{ij}$ for any $(i,j)$ (thus, $(\Sigma,\delta,\gamma)$ represents $S^3$),
\item $\delta_{1}$ intersects $\alpha_{1}$ and does not intersect $\alpha_{2}$.
\item $\delta_{2}$ intersects $\alpha_{2}$ and does not intersect $\alpha_{1}$.
\end{itemize}

\begin{figure}[h]

\includegraphics[width=5cm,clip]{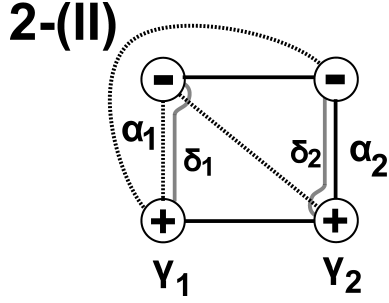} \\
\caption{$\delta = (\delta_{1},\delta_{2})$}\label{delta2}
\end{figure}

In this new diagram $(\Sigma,\delta,\gamma)$, we can take four unknots $K_{1}$, $K_{2}$, $C_{1}$ and $C_{2}$. $K_{1}$ and $K_{2}$ are similar in the above case. $C_{1}$ and $C_{2}$ are two unknots in $U_{\delta}$ whose framing is $\alpha_{1}$ and $\alpha_{2}$. These slopes can be determined as follows.

Let $r_{\alpha_{1}}$, $r_{\alpha_{2}}$, $r_{\beta_{1}}$ and $r_{\beta_{2}}$ be the rational numbers representing the surgery framings of $\alpha_{1}$, $\alpha_{2}$, $\beta_{1}$ and $\beta_{2}$ respectively. Precisely, put $r_{\alpha_{i}} = \rm{sgn}(r_{\alpha_{i}}) p_{i}/q_{i}$ and $r_{\beta_{i}} = \rm{sgn}(r_{\beta_{i}}) p'_{i}/q'_{i}$ for $i = 1,2$. Then,

\begin{itemize}
\item $\rm{sgn}(r_{\alpha_{1}}) = \rm{sgn}(r_{\beta_{1}}) = +1$,
\item $\rm{sgn}(r_{\alpha_{2}}) = \rm{sgn}(r_{\beta_{2}}) = -1$,
\item $p_{i} = \#(\alpha_{i} \cap \gamma_{i})$, for $i = 1,2$, 
\item $q_{1} = \#(\Gamma(+1,-2))$,
\item $q_{2} = \#(\Gamma(+2,+1))$.
\item $p'_{1} = \#(\beta_{1} \cap \alpha_{2})/\#(\Gamma(+2,+1))$,
\item $p'_{2} = \#(\beta_{2} \cap \alpha_{1})/\#(\Gamma(+2,-1))$,
\item $q'_{i} = \#(\beta_{i} \cap \gamma_{i})$, for $i = 1,2$. 
\end{itemize}

Note that $p_{i} > q_{i}$ and $p'_{i} > q'_{i}$ for $i = 1,2$.
Thus, $|r_{\alpha_{i}}| > 1$ and $|r_{\beta_{i}}| > 1$ for $i = 1,2$. (If $\gamma_{i} = \beta_{i}$, take $r_{\beta_{i}} = \infty$.)

\begin{figure}[h]

\includegraphics[width=11cm,clip]{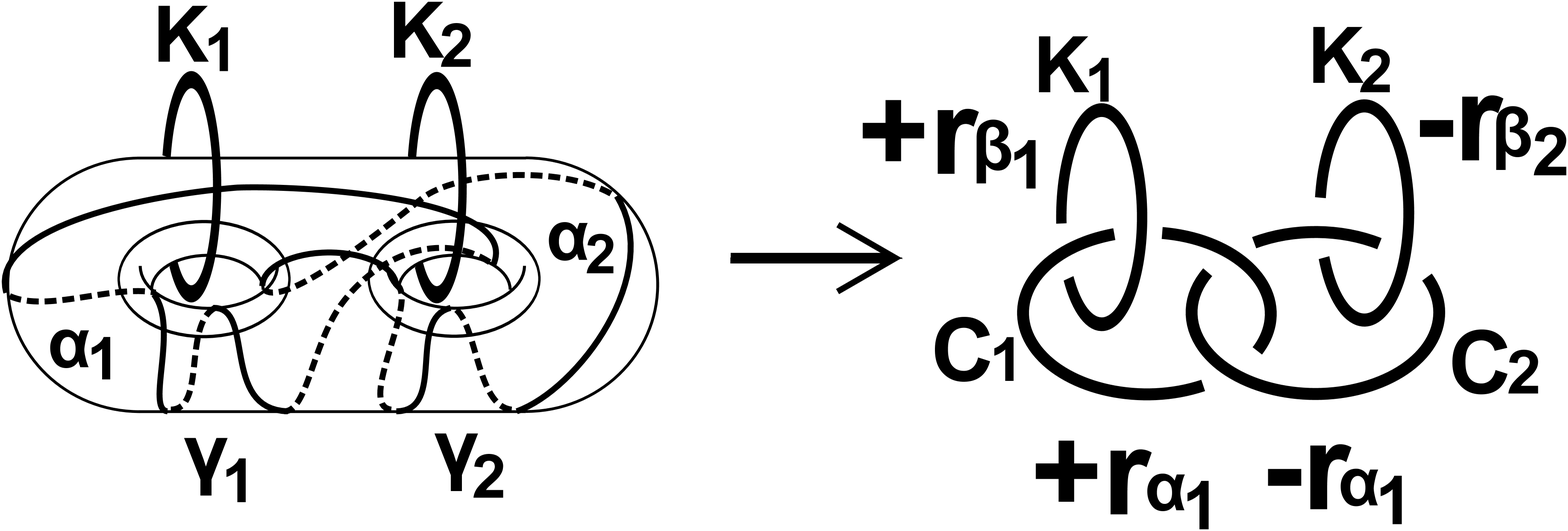} \\
\caption{surgery representation of $Y$}\label{k1k2-2}
\end{figure}

As a result, the three manifold $Y$ obtained from $(\Sigma, \alpha, \beta)$ can be represented as $Y = S^{3}(K_{1},K_{2},C_{1},C_{2})$ (see Figure \ref{k1k2-2}).

It is similar to prove that $Y$ belongs to $\mathcal{M}_{\mathcal{T}} = \mathcal{L}_{\overline{\rm{Brm}}}$. 

\subsection{Non-maximal cases for $g = 2$}\label{nonmax}

In this subsection, we study the case where the induced matrix $A_{(\Sigma,\alpha,\beta)}$ is effective, but not maximal. 

\begin{figure}[h]

\includegraphics[width=5cm,clip]{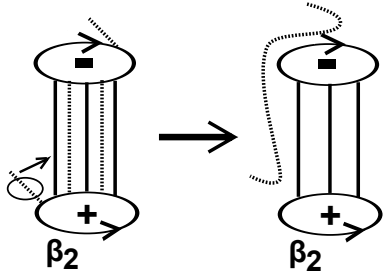} \\
\caption{not maximal case}\label{excep1}
\end{figure}

Let $(\Sigma,\alpha,\beta)$ be a strong diagram with genus two. If the induced matrix is not $A_{2}$, then
\[ A_{(\Sigma,\alpha,\beta)} \sim \left (\begin{array}{@{\,}cccc@{\,}}
+     & 0    \\
0     & +     \end{array} \right)
\text{ or } 
\left (\begin{array}{@{\,}cccc@{\,}}
+     & +    \\
0     & +     \end{array} \right).\]
The first matrix implies that $Y$ becomes a connected sum of Lens space.
If $A_{(\Sigma,\alpha,\beta)}$ is the second matrix, we find that $\Gamma(+2,-2) \neq \emptyset$, $\Gamma(+1,-2) \neq \emptyset$ and $\Gamma(+2,-1) \neq \emptyset$. So the neighborhood of $\beta_{2}$ looks as in Figure \ref{excep1}.
Let $x = (n_{1}, m_{1},\cdots,m_{k-1}, n_{k})$ be the sequence of integers representing $\Gamma(+2,-2)$, where $n$ means $\beta_{2} \cap \alpha_{1}$ and $m$ means $\beta_{2} \cap \alpha_{2}$. 
Note that $x$ can not be of the form $(m_{1},n_{1},\cdots,n_{k-1},m_{k})$.
Since $\Gamma(+2,+1) = \emptyset$, we find that $n_{i} = 1$ for all $i$.
Thus, we can transform the diagram by handle-slides $\alpha_{1} \leadsto \alpha_{2}$ (see Figure \ref{excep1}). 

This new diagram implies that $Y$ is a connected sum of lens spaces.

\section{Proof of Theorem \ref{classify}}

\subsection{Types of $A'_{3}$-strong diagrams for $g = 3$}

Let $(\Sigma,\alpha,\beta)$ be a strong diagram representing $Y$ with genus three. 
Suppose that the equivalence class of the induced matrix satisfies $[A'_{3}] \le [A_{(\Sigma,\alpha,\beta)}]$. We call such a diagram \textit{an } $A'_{3}$\textit{-strong diagram}. Recall that $[A'_{3}] \le [A_{3}]$.
There are just $22$ $\Gamma$-sets which may have some elements as follows.

$$\Gamma_{1}(+ \{1,2,3 \}, - \{1,2,3 \}),$$
$$\Gamma_{2}(-1,+1),
\Gamma_{2}(+ \{2,3\}, - \{2,3\}),
\Gamma_{2}(+ \{2,3\},+ 1) ,
\Gamma_{2}(- 1, - \{2,3\}),$$
$$\Gamma_{3}(-2,+2),
\Gamma_{3}(+ 3, - 3),
\Gamma_{3}(+ 3, +2),
\Gamma_{3}(- 2, - 3),$$ 
where $\{1,2,3\}$ means $1$ or $2$ or $3$.

\begin{prop}\label{A3}
Let $(\Sigma,\alpha,\beta)$ be an $A'_{3}$-strong diagram. 
Then,  $(\Sigma,\alpha,\beta)$ can be transformed by handle-slides and isotopies (if it is necessary) so that the new diagram is strong and $\Gamma(+ j, - j)$ has at least one element for each $j = 1,2,3$. 
%the one of the following conditions satisfies.
%\begin{itemize}
%\item the new diagram is also $A_{2}$-strong and $\Gamma(+ j, - j)$ has at least one element for each $j = 1,2,3$. 
%\item the new diagram is strong but not $A_{2}$-strong. 
%\end{itemize}
\end{prop}

\begin{prf}
We prove this proposition in three steps. Compare this proof with the proof of Proposition \ref{2ste}.
\begin{enumerate}
\item We can transform the diagram so that $\Gamma(+ 3, - 3) \neq \emptyset$.  
\item If $\Gamma(+ 3, - 3) \neq \emptyset$, we can transform the diagram so that $\Gamma(+ 2, - 2) \neq \emptyset$.
\item If $\Gamma(+ 3, - 3) \neq \emptyset$ and $\Gamma(+ 2, - 2) \neq \emptyset$, we can transform the diagram so that $\Gamma(+ 1, - 1) \neq \emptyset$.
\end{enumerate}
In each step, we must take a strong diagram.

\underline{Step 1}
This step is the same as step 1 of Proposition \ref{2ste}. 
Let $(\Sigma,\alpha,\beta)$ be an $A'_{3}$-strong diagram. Suppose $\Gamma(+ 3, - 3) = \emptyset$.
Since $\# (\alpha_{3} \cap \beta_{3}) \neq 0$, we find that $\Gamma(+ 3, + 2) \neq \emptyset$ and $\Gamma(- 2, - 3) \neq \emptyset$. Then, we can transform the diagram by a handle-slide $\beta_{3}^{+} \leadsto \beta_{2}^{+}$ similarly. 

\begin{figure}[h]

\includegraphics[width=7cm,clip]{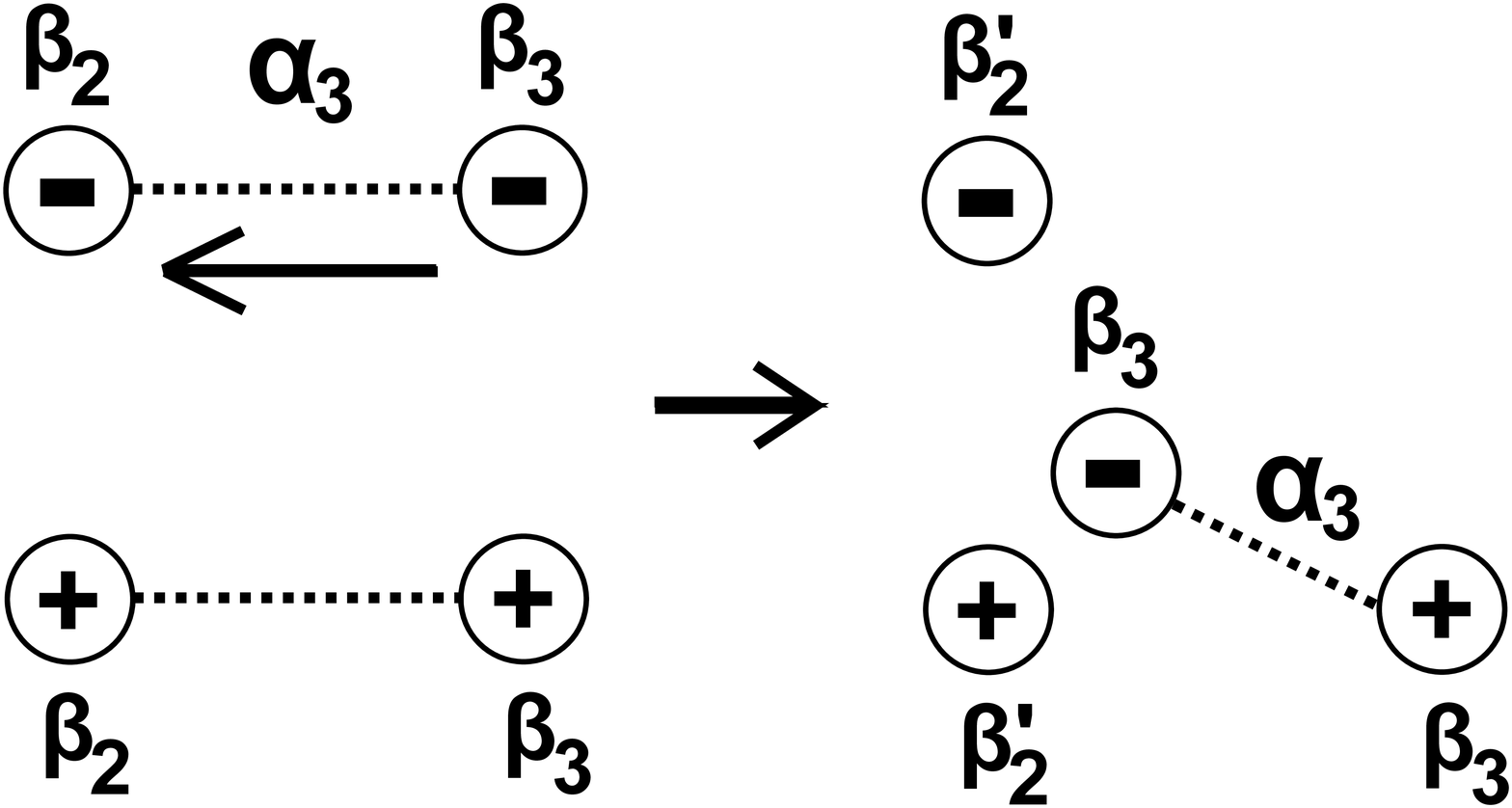} \\
\caption{}\label{selfcurve5}
\end{figure}

\underline{Step 2}
This step is also the same as step 2 in Proposition \ref{2ste}.
Let $(\Sigma,\alpha,\beta)$ be an $A'_{3}$-strong diagram with $\Gamma(+ 3, - 3) \neq \emptyset$.
Suppose that $\Gamma(+ 2, - 2) = \emptyset$.
Since $\# (\alpha_{3} \cap \beta_{2}) \neq 0$, we get that $\Gamma(+ 3, + 2) \neq \emptyset$ and  $\Gamma(- 2, - 3) \neq \emptyset$. Then, we can transform the diagram by handle-slides $\beta_{2}^{+} \leadsto \beta_{3}^{+}$ finitely many times (see Figure \ref{selfcurve6}). In finitely many steps, we will get an strong Heegaard diagram where $\Gamma(+ 2, - 2) \neq \emptyset$. We also get $\Gamma(+ 3, - 3) \neq \emptyset$ because the set $\alpha_{2} \cap \beta_{3}$ becomes non empty after these handle-slides.

\begin{figure}[h]

\includegraphics[width=7cm,clip]{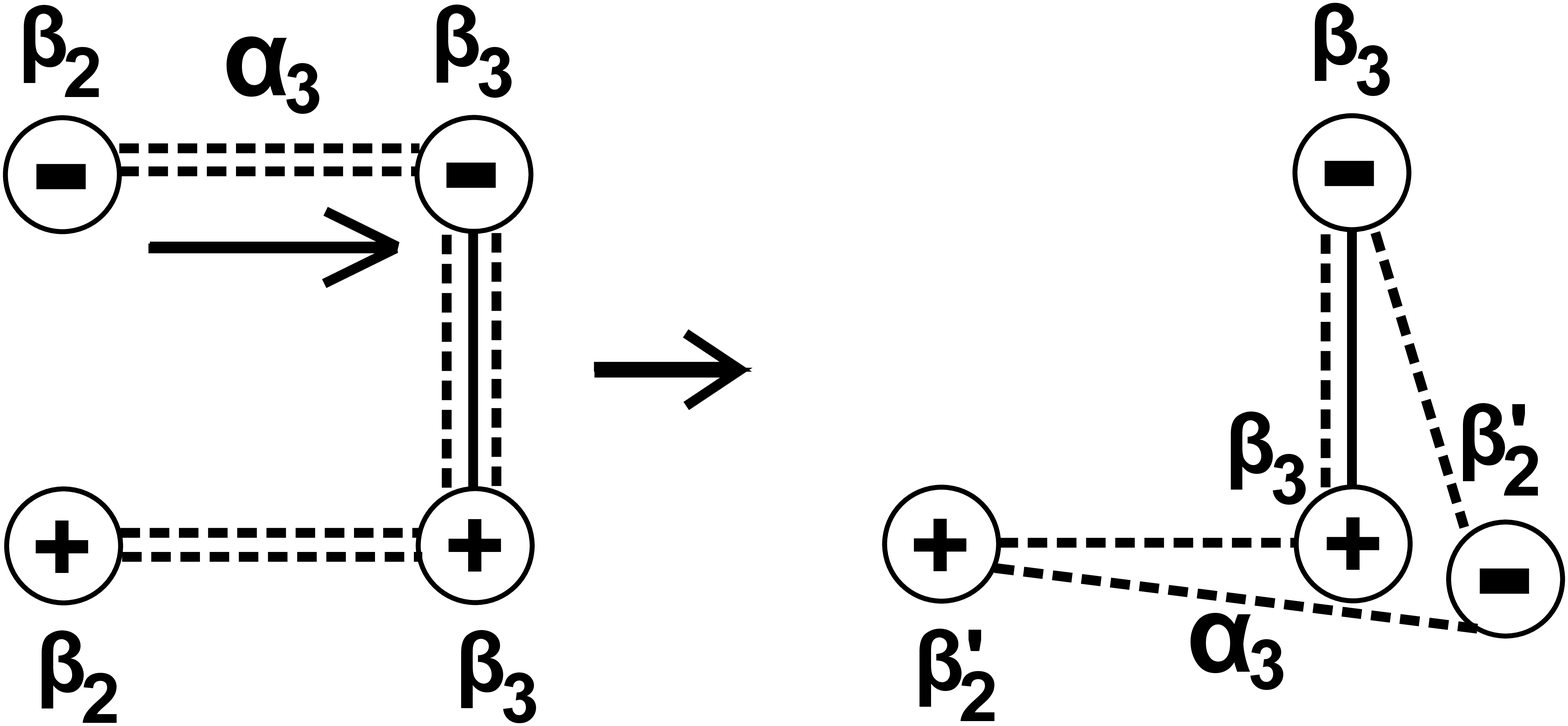} \\
\caption{}\label{selfcurve6}
\end{figure}

\underline{Step 3}
Let $(\Sigma,\alpha,\beta)$ be an $A'_{3}$-strong diagram with $\Gamma(+ 3, - 3) \neq \emptyset$, $\Gamma(+ 2, - 2) \neq \emptyset$. Suppose that $\Gamma(+ 1, - 1) = \emptyset$.
Since $\# (\alpha_{3} \cap \beta_{2}) \neq 0$ by the assumption, we can get Figure \ref{selfcurve1}. Moreover, any two arcs in $\Gamma(+ 3, - 3)$, $\Gamma(+ 3, + 2)$, $\Gamma(- 3, - 2)$ and $\Gamma(+ 2, - 2)$ are isotopic to each edge of the rectangle respectively. Since $\alpha_{1} \cap \beta_{1} \neq \emptyset$ and $\alpha_{2} \cap \beta_{1} \neq \emptyset$, we get that $\beta_{1}^{+}$ and $\beta_{1}^{-}$ are in or out of the rectangle. Actually, if $\beta_{1}^{+}$ is between two arcs in $\Gamma(-3,-2)$, we find that $\alpha_{1}$ can not intersect $\beta_{1}$. Moreover, if $\beta_{1}^{+}$ is between two arcs in $\Gamma(+3,-3)$, we can transform the diagram by handle-slides $\beta_{1}^{+} \leadsto \beta_{3}^{-}$ finitely many times so that $\beta_{1}^{+}$ is in or out of the rectangle (see Figure \ref{selfcurve2}). 

\begin{figure}[h]

\includegraphics[width=5cm,clip]{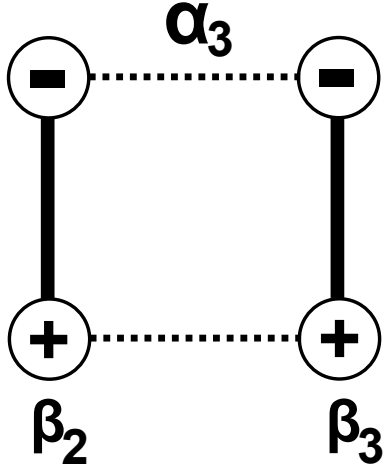} \\
\caption{$\Gamma(+ 3, - 3)$, $\Gamma(+ 3, + 2)$, $\Gamma(- 3, - 2)$ and $\Gamma(+ 2, - 2)$ make a rectangle} \label{selfcurve1}
\end{figure}
\begin{figure}[h]

\includegraphics[width=8cm,clip]{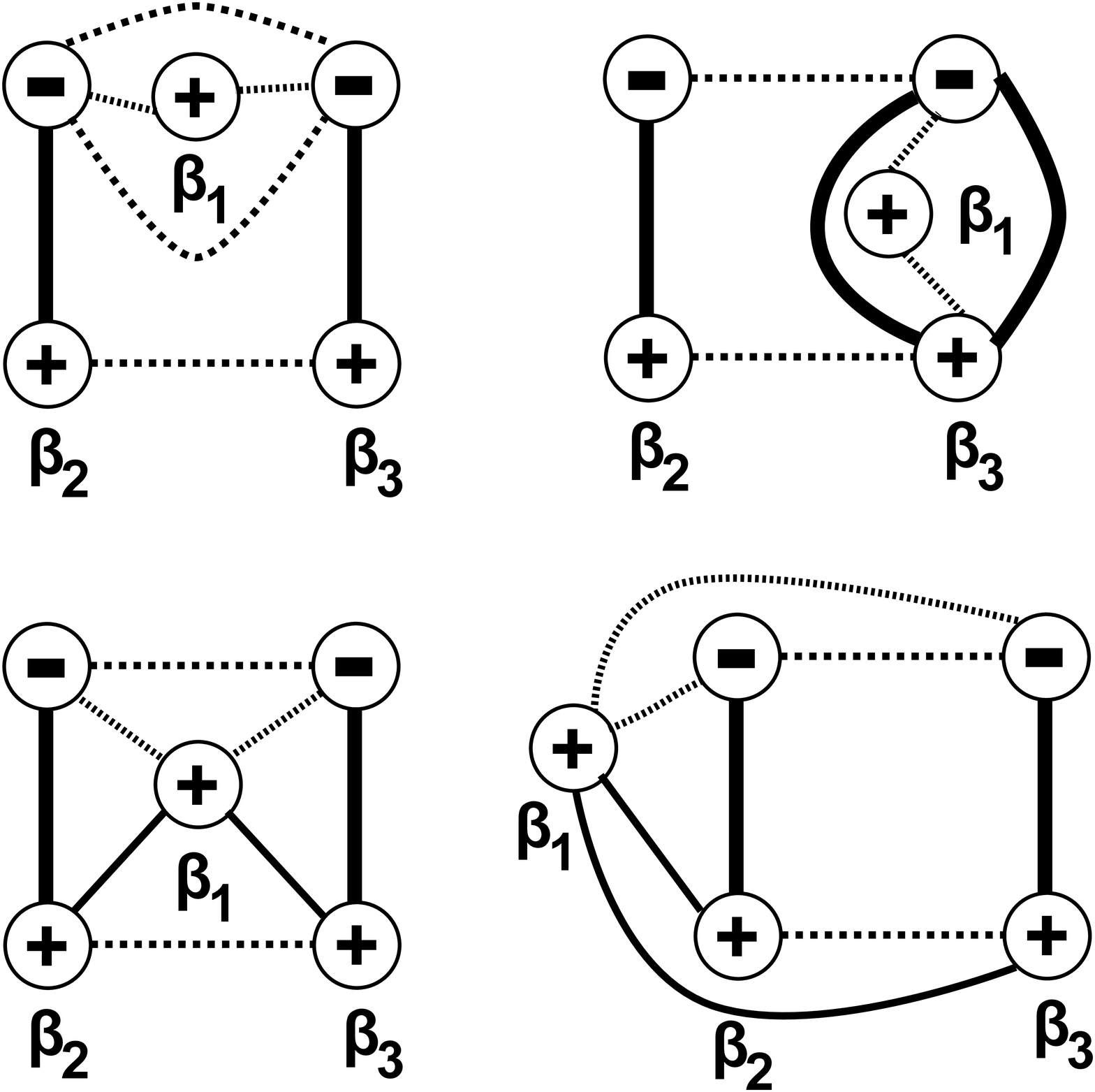} \\
\caption{positions of $\beta_{1}^{+}$}\label{selfcurve2}
\end{figure}

Assume $\beta_{1}^{+}$ is in the rectangle. There are four possible cases.

\begin{itemize}
\item If $\Gamma(+ 1, - 2) \neq \emptyset$, $\Gamma(+ 1, - 3) \neq \emptyset$, $\Gamma(+ 1, + 2) \neq \emptyset$, $\Gamma(+ 1, + 3) \neq \emptyset$ and $\Gamma(+ 3, - 2) \neq \emptyset$, then $\beta_{1}^{-}$ is in one of the three domains adjacent to $\beta_{3}^{-}$. But one of them is impossible because $\Gamma(+2, -1) \neq \emptyset$ (see Figure \ref{selfcurve3}). In the other two cases, we can transform the diagram by handle-slides $\beta_{1}^{-} \leadsto \beta_{3}^{-}$ finitely many times. Thus, we get $\Gamma(+ 1, - 1) \neq \emptyset$. Of course, the diagram is strong and $\Gamma(+ 3, - 3) \neq \emptyset$ and $\Gamma(+ 2, - 2) \neq \emptyset$.

\begin{figure}[h]

\includegraphics[width=5cm,clip]{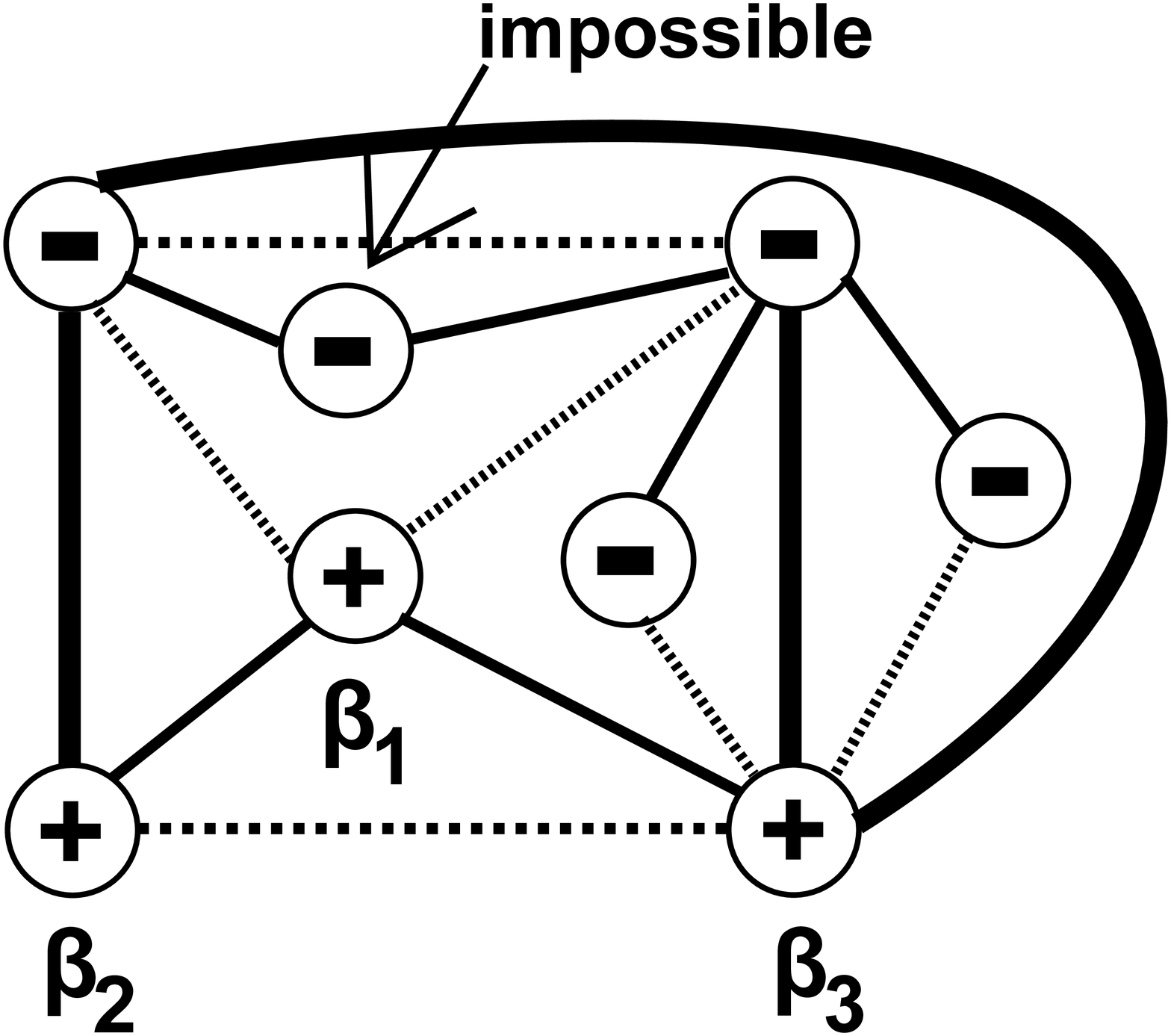} \\
\caption{}\label{selfcurve3}
\end{figure}

\item If $\Gamma(+ 1, - 2) \neq \emptyset$, $\Gamma(+ 1, - 3) \neq \emptyset$, $\Gamma(+ 1, + 2) \neq \emptyset$, $\Gamma(+ 1, + 3) \neq \emptyset$ and $\Gamma(- 3, + 2) \neq \emptyset$, we can get $\Gamma(+ 1, - 1) \neq \emptyset$ similarly.

\item If $\Gamma(+ 1, - 2) \neq \emptyset$, $\Gamma(+ 1, - 3) \neq \emptyset$, $\Gamma(+ 1, + 2) \neq \emptyset$, $\Gamma(+ 1, + 3) \neq \emptyset$, $\Gamma(+ 3, - 2) = \emptyset$ and $\Gamma(- 3, + 2) = \emptyset$, then $\beta_{1}^{-}$ is in one of the following two domains as in Figure \ref{selfcurve7} because $\# (\alpha_{1} \cap \beta_{3}^{+}) = \# (\alpha_{1} \cap \beta_{3}^{-})$. Thus, we can also transform the diagram by handle-slides $\beta_{1}^{-} \leadsto \beta_{3}^{-}$ finitely many times so that we get $\Gamma(+ 1, - 1) \neq \emptyset$.
\begin{figure}[h]

\includegraphics[width=5cm,clip]{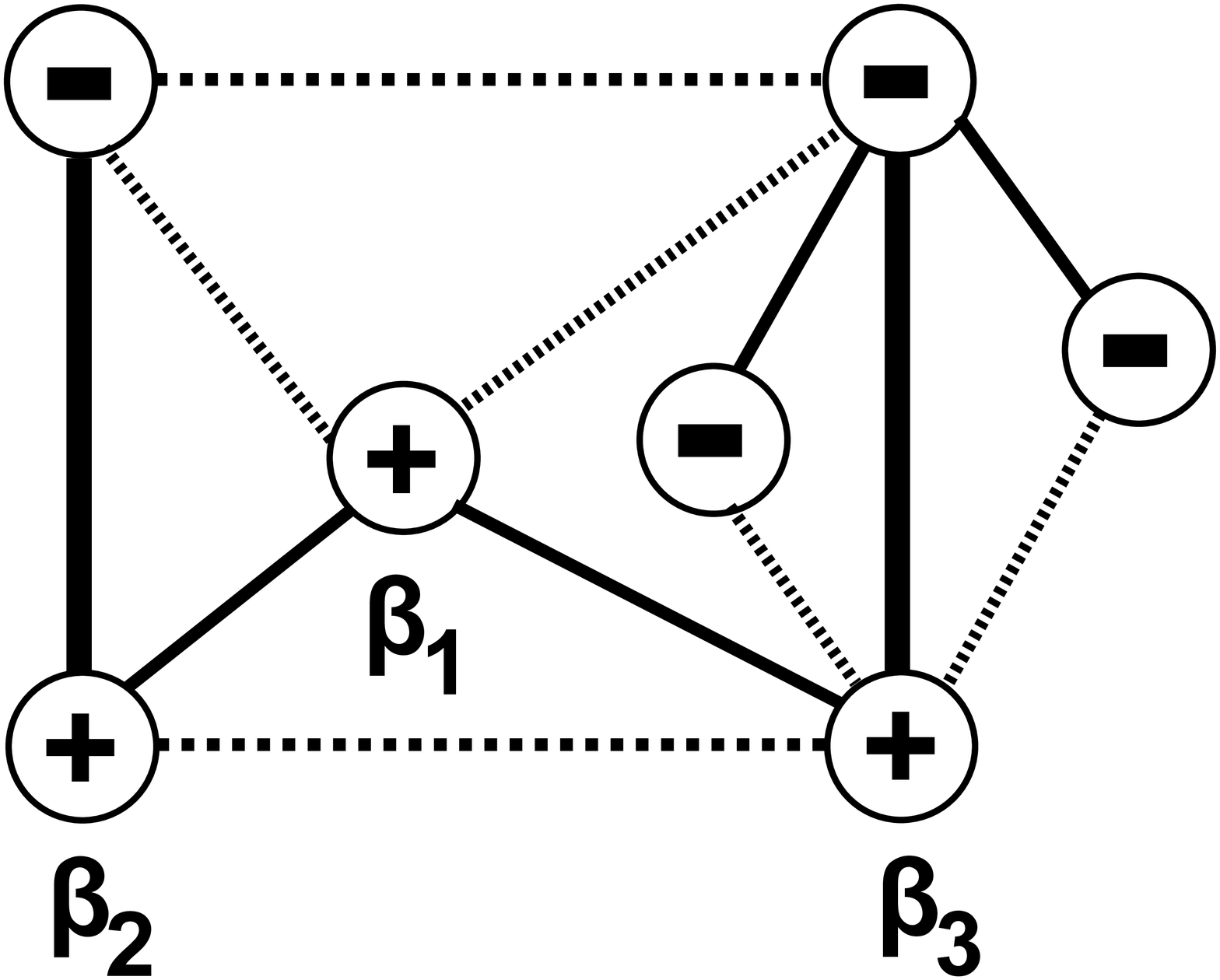} \\
\caption{}\label{selfcurve7}
\end{figure}
\item If one of $\Gamma(+ 1, - 2)$, $\Gamma(+ 1, - 3)$, $\Gamma(+ 1, + 2)$, $\Gamma(+ 1, + 3)$ is the empty set, then we can transform the diagram by handle-slides finitely many times so that we get $\Gamma(+ 1, - 1) \neq \emptyset$ as follows.

\begin{itemize}
\item $\Gamma(+ 1, + 3) = \emptyset \Rightarrow \beta_{1}^{+} \leadsto \beta_{2}^{+}$.  
\item $\Gamma(+ 1, - 3) = \emptyset \Rightarrow \beta_{1}^{+} \leadsto \beta_{2}^{-}$.
\item $\Gamma(+ 1, + 2) = \emptyset \Rightarrow \beta_{1}^{+} \leadsto \beta_{3}^{+}$.
\item $\Gamma(+ 1, - 2) = \emptyset \Rightarrow \beta_{1}^{+} \leadsto \beta_{3}^{-}$.
\end{itemize}
\end{itemize}
\qed
\end{prf}

\begin{figure}[h]

\includegraphics[width=12cm,clip]{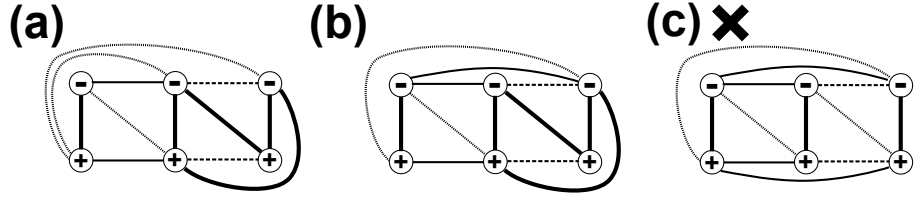} \\
\caption{three possoble types}\label{case2}
\end{figure}

\begin{prop} \label{A3next}
Let $(\Sigma,\alpha,\beta)$ be an $A'_{3}$-strong diagram. Suppose $\Gamma(+j,-j) \neq \emptyset$ for $j = 1,2,3$. Then, $(\Sigma,\alpha,\beta)$ can be transformed by handle-slides, isotopies, permutations of curves, changes of orientations(if it is necessary) so that the new strong diagram is of one of the following three types 3-(a), 3-(b) and 3-(c) (see Figure \ref{case2}).
\end{prop}

\begin{prf}
By Proposition \ref{A3}, we can assume $\Gamma(+ j, - j) \neq \emptyset$ for all $j$.
We put $\beta_{1}$ and $\beta_{2}$ as in Figure \label{2type1}.
Then, $\beta_{3}$ is in one of the three domains as in Figure \label{2type1}.

\begin{figure}[h]

\includegraphics[width=7cm,clip]{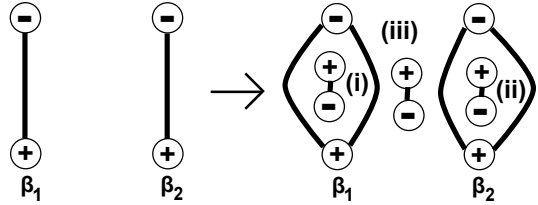} \\
\caption{three positions}\label{2type1}
\end{figure}

\begin{itemize}
\item The position (i) is impossible because $\Gamma(+ 3, + 2) \neq \emptyset$. 
\item If $\beta_{3}$ is in (ii), $\beta_{1}$ and $\beta_{2}$ looks as in Figure \ref{2type2}. Since $\Gamma(+1,-2) \neq \emptyset$ and $\Gamma(-1,-2) \neq \emptyset$, we can transform the diagram by handle-slides $(\beta_{1}^{+} \text{ and } \beta_{1}^{-}) \leadsto \beta_{2}^{-}$ finitely many times so that the new diagram is also $A'_{3}$-strong and $\beta_{3}$ is in the position (iii).

\begin{figure}[h]

\includegraphics[width=7cm,clip]{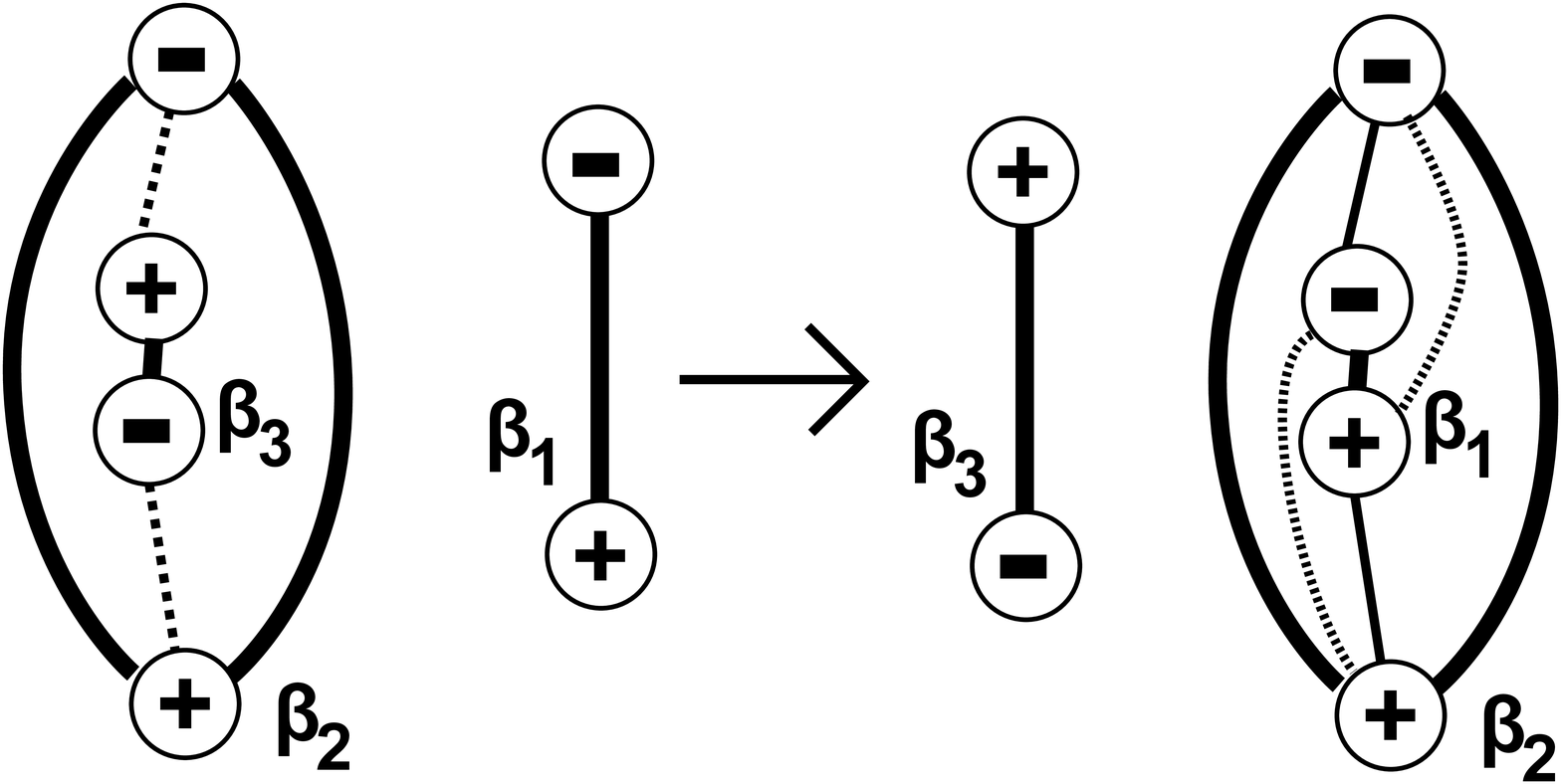} \\
\caption{}\label{2type2}
\end{figure}

\item If $\beta_{3}$ is in (iii), there exist at most two isotopy classes of arcs in $\Gamma(+ 3, -3)$ (see \ref{2type3}). If there exist two arcs which are not isotopic, then we can transform the diagram by handle-slides $(\beta_{2}^{+} \text{ and } \beta_{2}^{-}) \leadsto \beta_{3}^{-}$ finitely many times so that all arcs in $\Gamma(+ j, - j)$ are isotopic (see Figure \ref{2type3}).

\begin{figure}[h]

\includegraphics[width=6cm,clip]{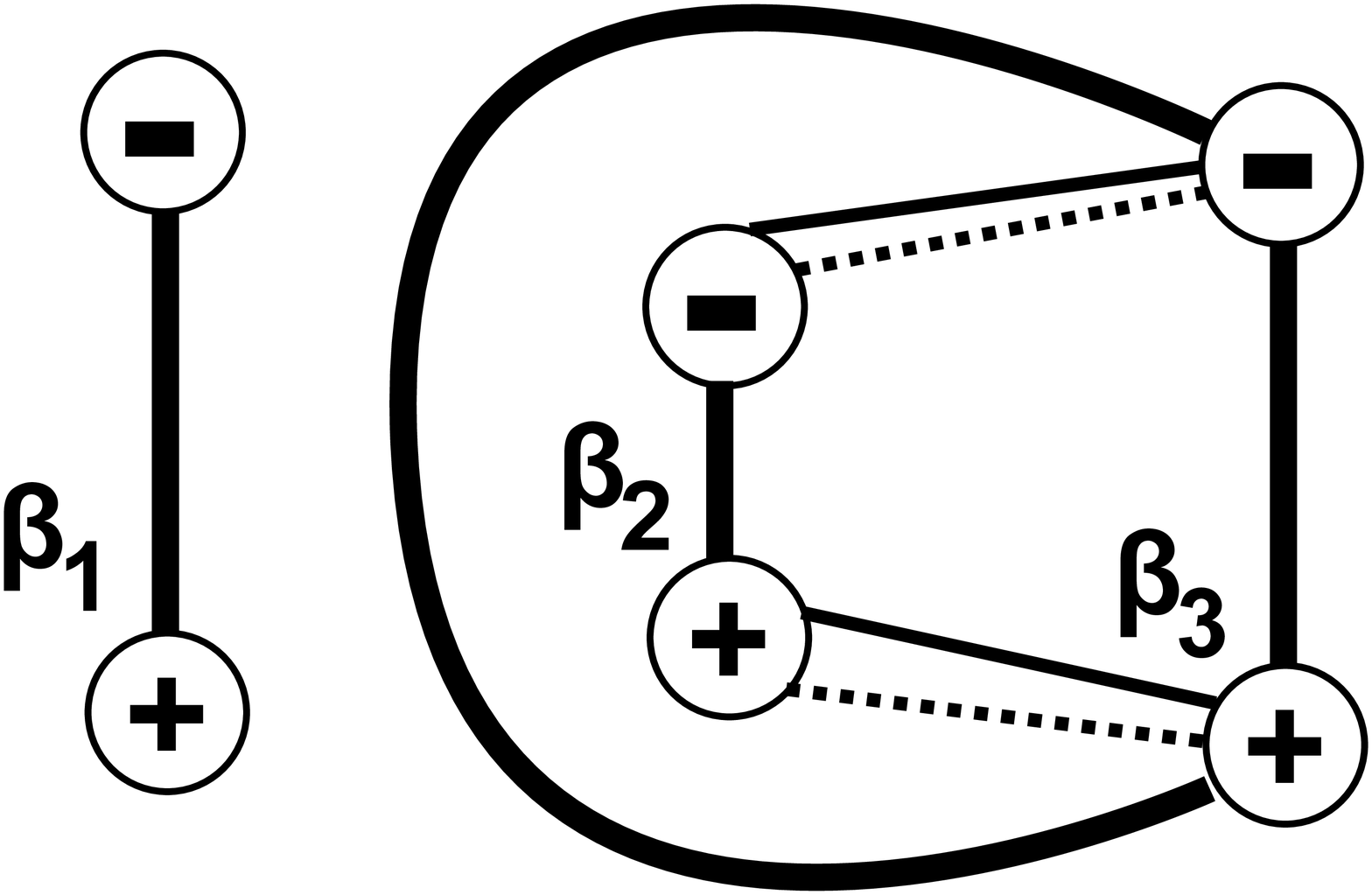} \\
\caption{}\label{2type3}
\end{figure}

\end{itemize}

Since $\Gamma(+ 3, + 2) \neq \emptyset$ and $\Gamma(- 3, - 2) \neq \emptyset$, we get two possible cases as in Figure \ref{2type4}. But they are equivalent under permutating $\beta_{2}$ and $\beta_{3}$ and reversing the orientation of $\alpha_{3}$.

\begin{figure}[h]

\includegraphics[width=12cm,clip]{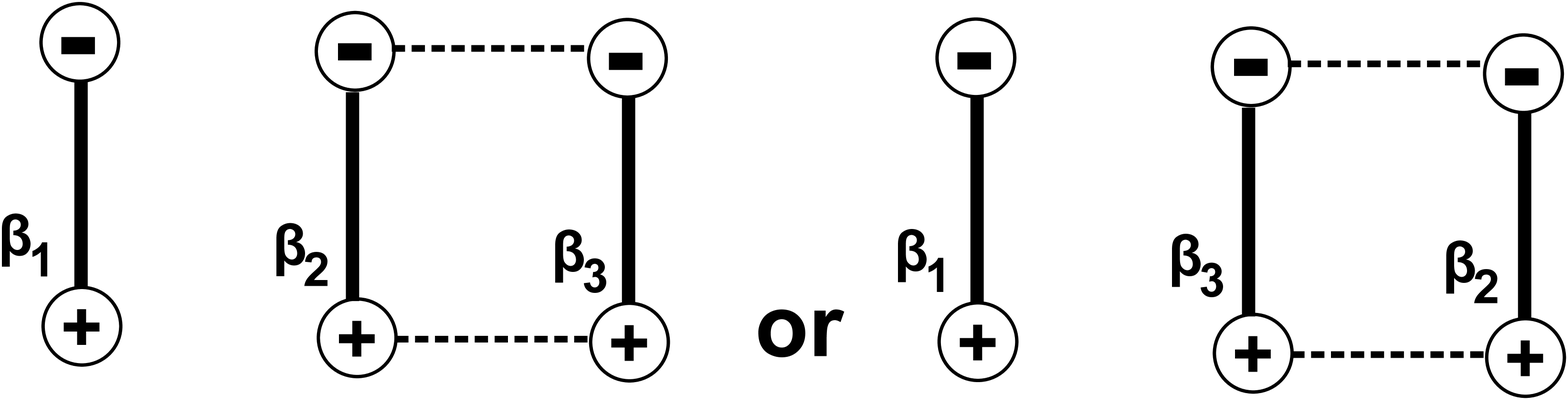} \\
\caption{}\label{2type4}
\end{figure}

Now we can describe all possible cases.
There are another $9$ cases, but we can transform these diagrams into one of the three types.
\qed
\end{prf}

\begin{figure}[h]

\includegraphics[width=12cm,clip]{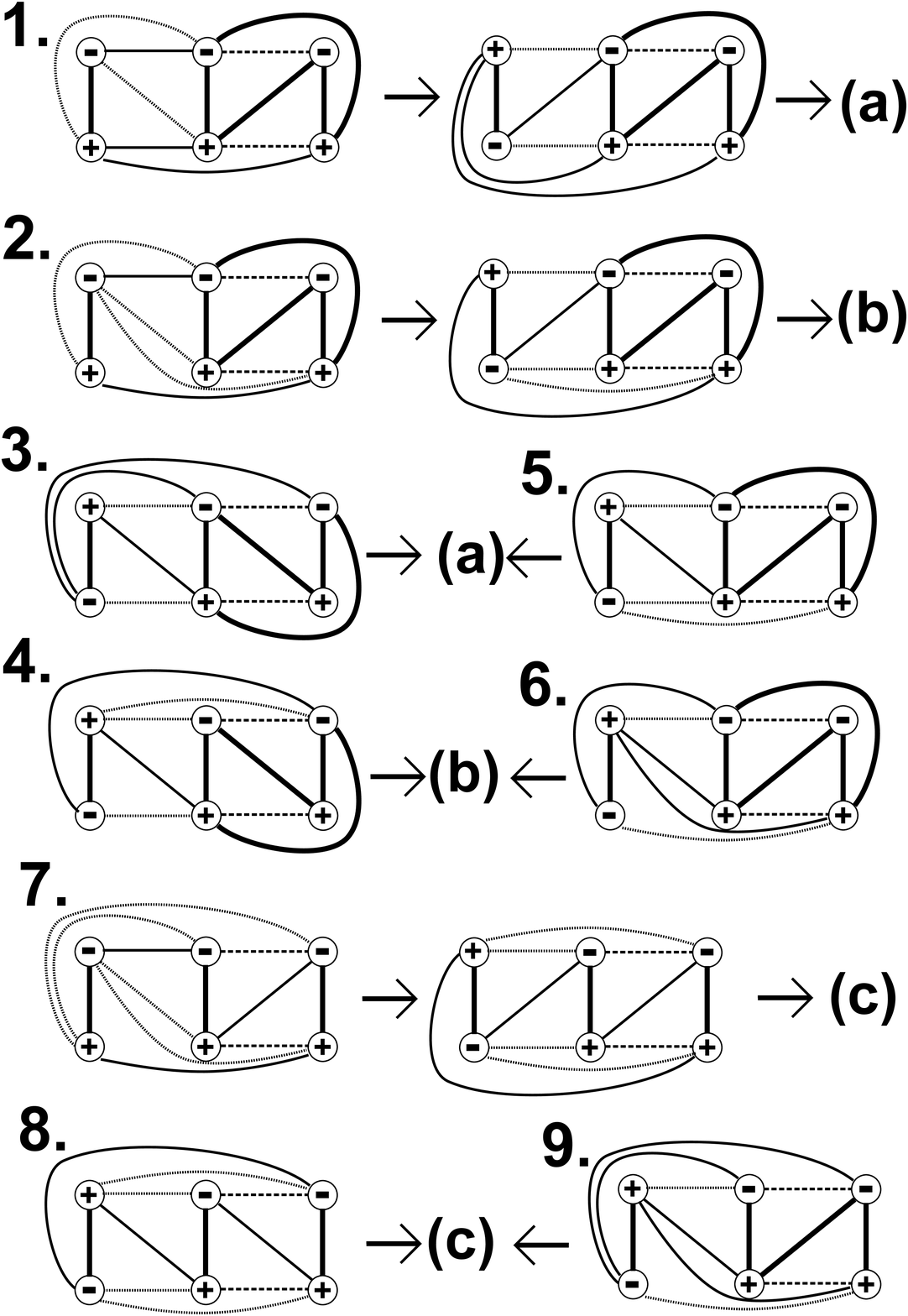} \\
\caption{}\label{case1}
\end{figure}

However, note that the type (c) never happens. Actually, the pattern of the intersection points at $\beta_{3}^{+}$ and $\beta_{3}^{-}$ never coincide, so we can not attach $\beta_{3}$-circles.

\subsection{Surgery representations for $g = 3$, $\beta$-curves}\label{surg3b}

Let $(\Sigma,\alpha,\beta)$ be an $A'_{3}$-strong diagram of type 3-(a) or 3-(b).

\underline{$\Gamma(+1,-1)$}:

We first consider the neighborhood of $\beta_{1}$. 
Let $n_{0} = \#\Gamma(+1,-2) + \#\Gamma(+1,-3)$ and $m_{0} = \#\Gamma(+1,+2) + \#\Gamma(+1,+3)$.
Then, the argument is the same as in the case when $g = 2$.

That is, we can define $x$ to be the sequence of the number of intersection points induced from $\Gamma(+1,-1)$ and one of the following two cases may happen, where $n$ comes from  $\beta_{1} \cap \alpha_{1}$ and $m$ comes from $\beta_{1} \cap \alpha_{2}$. 

\begin{itemize}
\item $x = (n_{1}, m_{1},\cdots,m_{k-1}, n_{k})$ where $m_{i} = m_{0}$ for all $i$ and $n_{1} = n_{k}$,
\item $x = (m_{1},n_{1},\cdots, n_{k-1},m_{k})$ where $n_{i} = n_{0}$ for all $i$ and $m_{1} = m_{k}$.
\end{itemize}

In each case, take a simple closed oriented curve $\gamma_{1}$ on $\Sigma$ so that the following conditions hold. If $\Gamma(+1, -1)$ consists of only $\alpha_{1}$-arcs (or only $\alpha_{2}$-arcs), then we take $\gamma_{1} = \beta_{1}$. 

\begin{itemize}
\item $\gamma_{1}$ intersects each arc in $\Gamma(+2,-1)$, $\Gamma(-1,-2)$ and $\Gamma(-1,-3)$ at a point. 
\item $\gamma_{1}$ intersects $\beta_{1}$ at some points, but does not intersect $\beta_{2}$ and $\beta_{3}$.
\item If $x = (n_{1}, m_{1},\cdots,m_{k-1}, n_{k})$, $\gamma_{1}$ intersects only $\alpha_{1}$-arcs in $\Gamma(+1,-1)$.
\item If $x = (m_{1},n_{1},\cdots, n_{k-1},m_{k})$, $\gamma_{1}$ intersects only $\alpha_{2}$-arcs in $\Gamma(+1,-1)$.
\end{itemize}

We consider the new $A'_{3}$-strong Heegaard diagram $(\Sigma, \alpha, (\gamma_{1}, \beta_{2},\beta_{3}))$. 
We describe $\beta_{1}$-curve more precisely later.

\begin{figure}[h]

\includegraphics[width=10cm,clip]{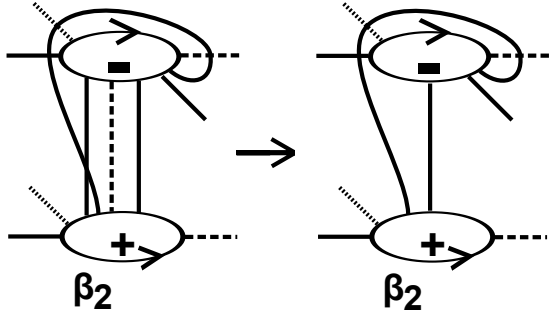} \\
\caption{$\gamma_{2}$-curve and new diagram $(\Sigma, \alpha, (\gamma_{1}, \gamma_{2},\beta_{3}))$}\label{nbd6}
\end{figure}

\underline{$\Gamma(+2,-2)$}:

Next, we consider the neighborhood of $\beta_{2}$-circles.
The neighborhood of $\beta_{2}$-circles looks as in Figure \ref{nbd6}.
Then, we can take $\gamma_{2}$ as follows (see Figure \ref{nbd6}). If there exists no $\alpha_{3}$-arc in $\Gamma(+2, -2)$, then we need not to take $\gamma_{2}$. 
\begin{itemize}
\item $\gamma_{2}$ intersects each arc in $\Gamma(-2,-3)$, $\Gamma(+1,-2)$ and $\Gamma(-1,-2)$ at a point. 
\item $\gamma_{1}$ intersects $\beta_{2}$ at some points, but does not intersect $\beta_{1}$ and $\beta_{3}$.
\item $\gamma_{1}$ does not intersect $\alpha_{3}$-arcs in $\Gamma(+2,-2)$.
\end{itemize}

We consider the new $A'_{3}$-strong Heegaard diagram $(\Sigma, \alpha, (\gamma_{1}, \gamma_{2},\beta_{3}))$. Then, the neighborhood of $\gamma_{2}$-circles looks as in Figure \ref{nbd6}. 
We describe $\beta_{2}$-curve more precisely later.

\begin{figure}[h]

\includegraphics[width=5cm,clip]{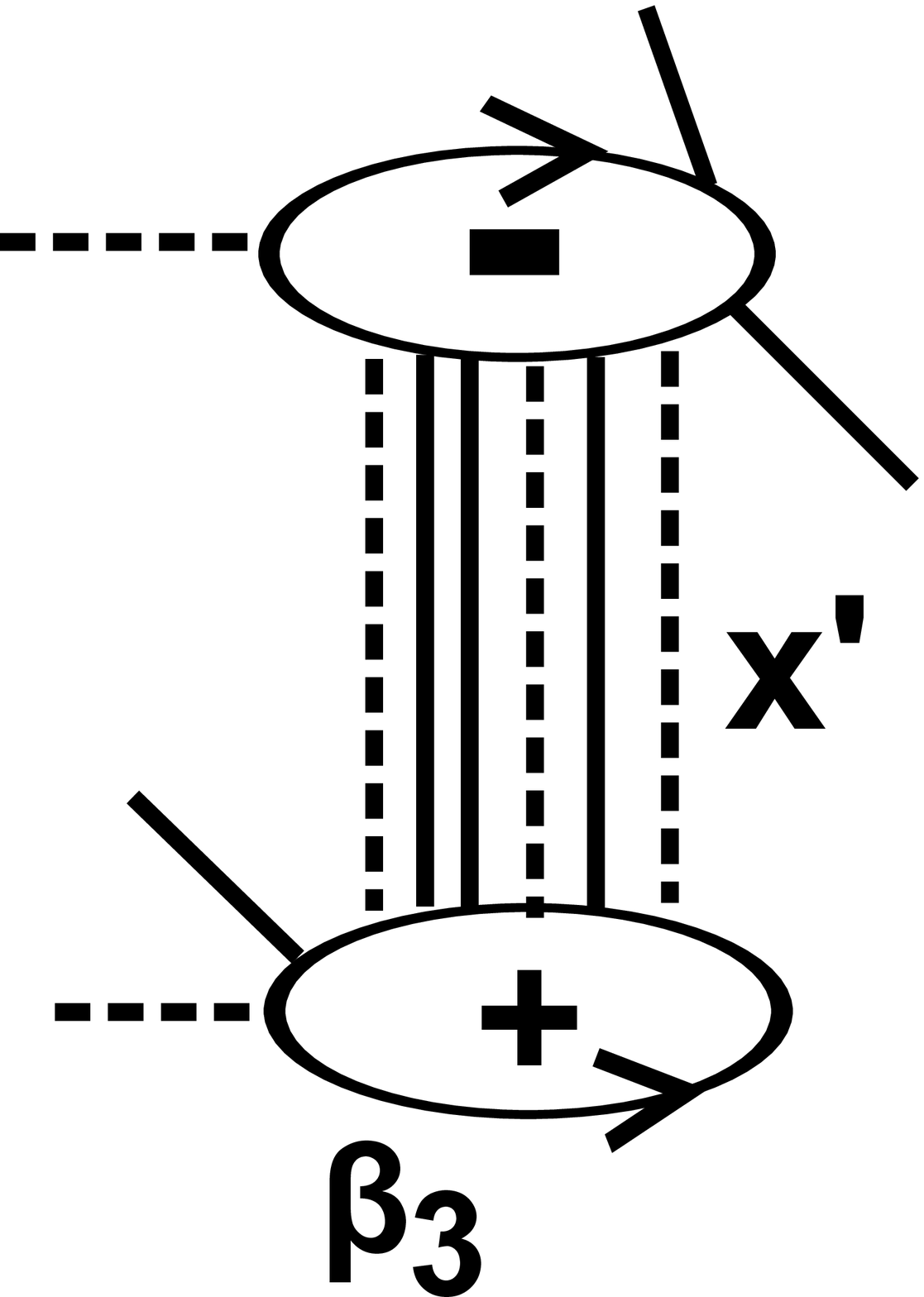} \\
\caption{}\label{nbd2}
\end{figure}

\underline{$\Gamma(+3,-3)$}:
Lastly, consider the neighborhood of $\beta_{3}$ looks as in Figure \ref{nbd2}.

Similarly, we can define $x'$ to be the the sequence of the number of intersection points induced from $\Gamma(+3,-3)$.
In this case, it is convinient not to distinguish $\alpha_{1}$ and $\alpha_{3}$. That is, there exist two cases,  where $n'$ comes from $\beta_{3} \cap \alpha_{3}$ and $m$ comes from $\beta_{3} \cap (\alpha_{1} \cup \alpha_{2})$. 

\begin{itemize}
\item $x' = (n'_{1}, m'_{1},\cdots,m'_{k-1}, n'_{k})$ where $m'_{i} = m'_{0}$ for all $i$ and $n'_{1} = n'_{k}$,
\item $x' = (m'_{1},n'_{1},\cdots, n'_{k-1},m'_{k})$ where $n'_{i} = n'_{0}$ for all $i$ and $m'_{1} = m'_{k}$.
\end{itemize}

\begin{figure}[h]

\includegraphics[width=7cm,clip]{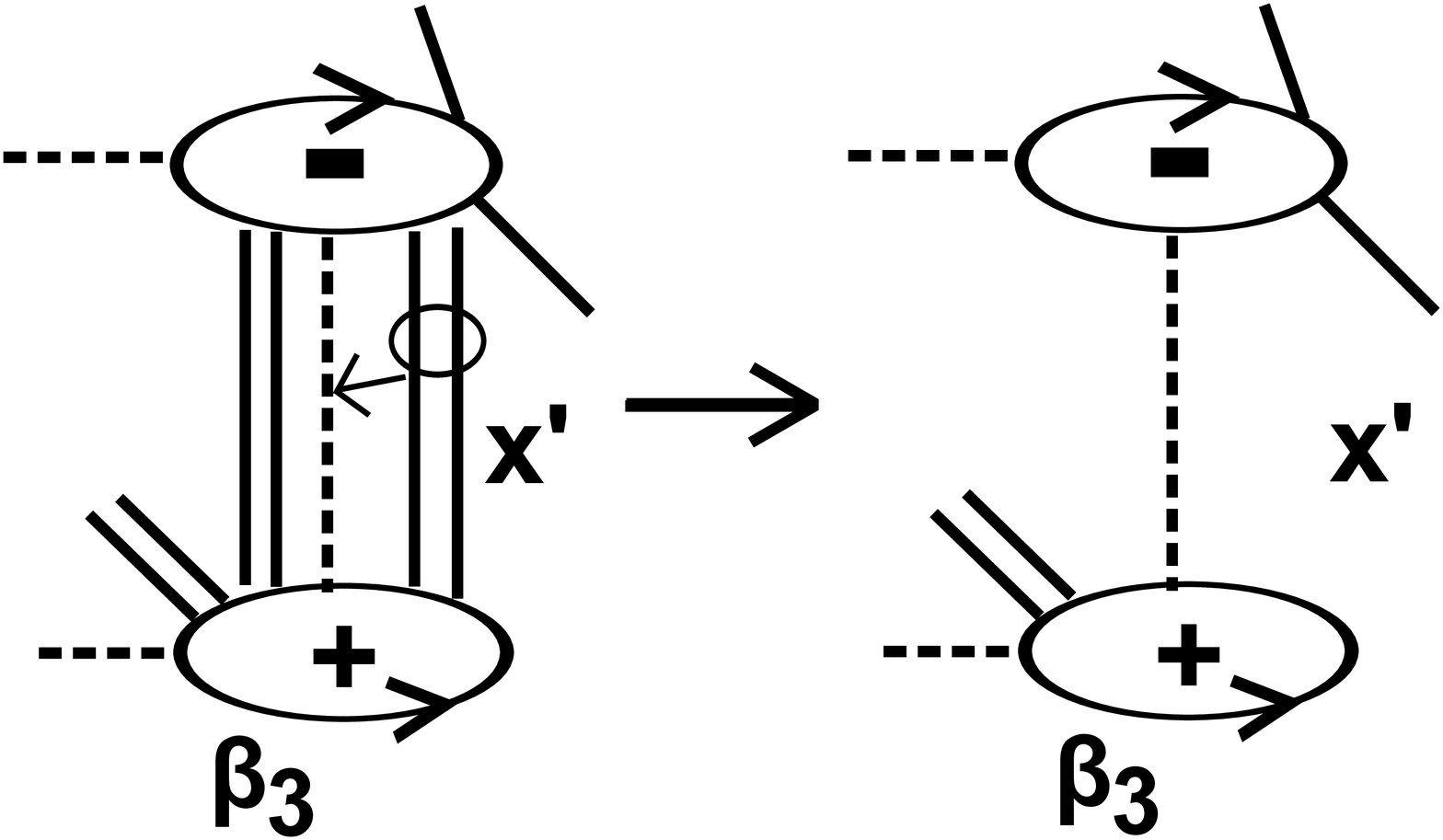} \\
\caption{}\label{nbd3}
\end{figure}

Moreover, if $x' = (m'_{1},n'_{1},\cdots, n'_{k-1},m'_{k})$, we find that $n'_{i} = 1$ for all $i$. Otherwise, $\alpha_{3}$-arcs can not be a closed curve. Then, We can transform the diagram by $m'_{1}$ handle-slides $+ (\alpha_{1} \text{ and } \alpha_{2}) \leadsto \alpha_{3}$ (see Figure \ref{nbd3}) so that the new diagram is also $A'_{3}$-strong diagram and $\Gamma_{i}(+3,-3) = \emptyset$ for $i = 1,2$. As a result, we return to the first case.

\begin{figure}[h]

\includegraphics[width=8cm,clip]{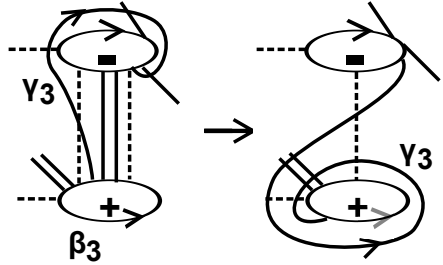} \\
\caption{$\gamma_{3}$ and new diagram $(\Sigma, \alpha, \gamma)$}\label{nbd7}
\end{figure}

If $x' = (n'_{1}, m'_{1},\cdots,m'_{k-1}, n'_{k})$, take a simple closed oriented curve $\gamma_{3}$ on $\Sigma$ similarly so that the following conditions hold (see Figure \ref{nbd7}). If $\Gamma(+3, -3)$ consists of only $\alpha_{3}$-arcs, then we take $\gamma_{3} = \beta_{3}$. 

\begin{itemize}
\item $\gamma_{3}$ intersects each arc in $\Gamma(-2,-3)$, $\Gamma(+1,-3)$ and $\Gamma(+2,-3) $at a point. 
\item $\gamma_{3}$ intersects $\beta_{3}$ at some points, but does not intersect $\beta_{1}$ and $\beta_{2}$.
\item $\gamma_{3}$ intersects only $\alpha_{3}$-arcs in $\Gamma(+2,-2)$.
\end{itemize}

We consider the new $A'_{3}$-strong Heegaard diagram $(\Sigma, \alpha, \gamma=(\gamma_{1}, \gamma_{2},\gamma_{3}))$. Then, the neighborhood of $\gamma_{2}$-circles looks as in Figure \ref{nbd7}. (In general, $\beta_{2}$ becomes arcs.) 
We describe $\beta_{3}$-curve more precisely later.

\begin{figure}[h]

\includegraphics[width=10cm,clip]{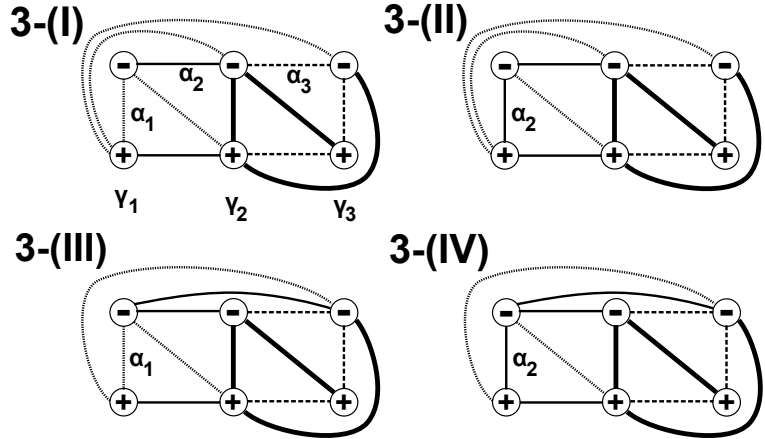} \\
\caption{possible four types}\label{4case1}
\end{figure}

Thus, there exist four possible types (see Figure \ref{4case1}).

\begin{figure}[h]

\includegraphics[width=7cm,clip]{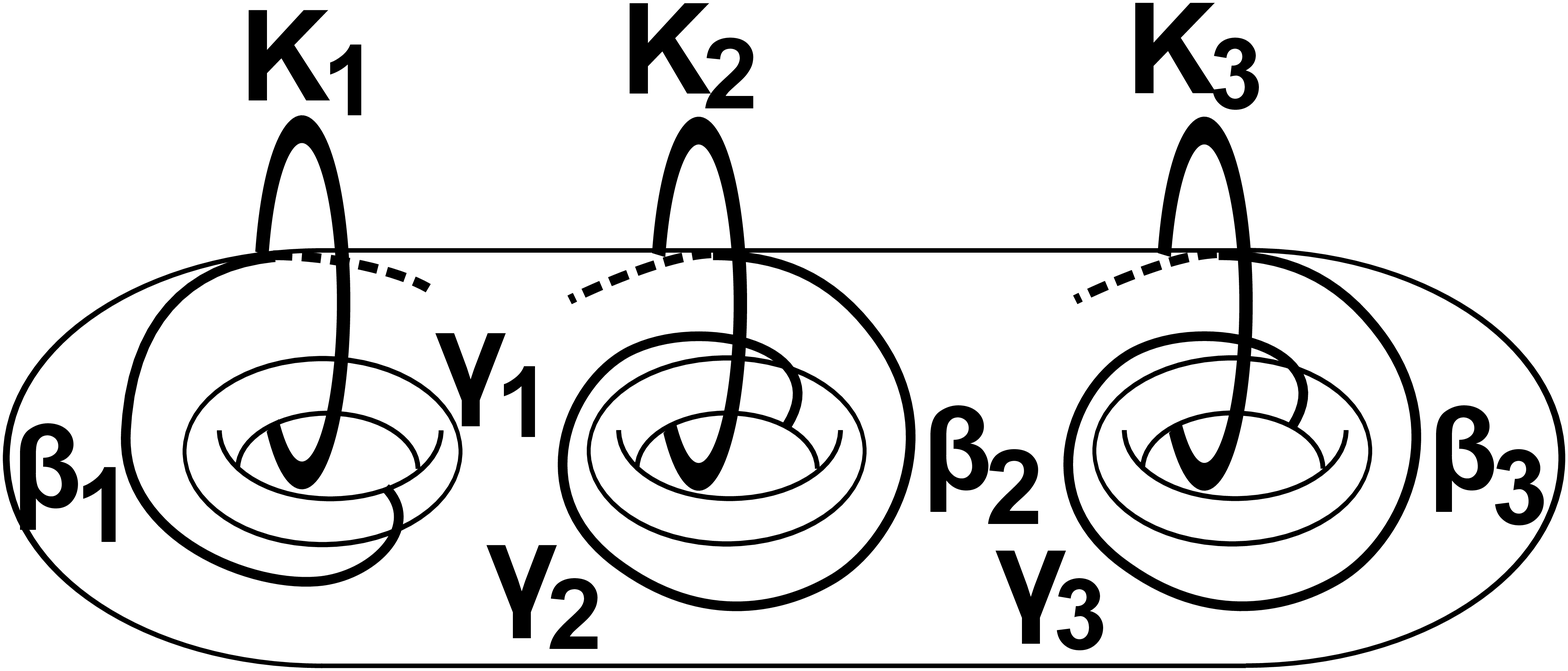} \\
\caption{}\label{k1k2k3}
\end{figure}

In each case, $\beta_{j}$-curves become some surgery framings of some unknots $K_{1}$, $K_{2}$ and $K_{3}$ by attaching $\gamma_{j}$-circles (see Figure \ref{k1k2k3}). 

\subsection{Surgery representations for $g = 3$, $\alpha$-curves}\label{surg3a}

\begin{figure}[h]

\includegraphics[width=7cm,clip]{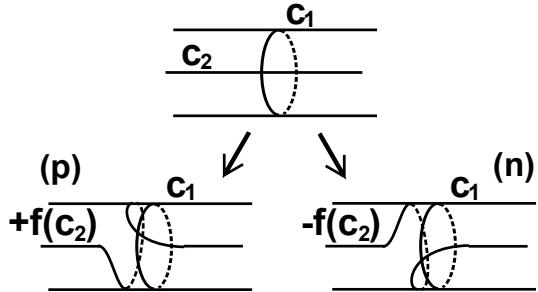} \\
\caption{positive/negative Dehm surgery}\label{dehn}
\end{figure}

Now, we recall positive (or negative) Dehn twists.
\begin{defn}
Let $\Sigma$ be a closed oriented genus $g$ surface and Let $c_{1}$ be a simple closed curve on $\Sigma$.
Then, a positive (or negative) Dehn twist is the self-homeomorphism $\pm f(c_{1})$ on $\Sigma$ defined as in Figure \ref{dehn}-(p) and (n) on a neighborhood of $c_{1}$, where a curve $c_{2}$ is mapped to $f(c_{2})$ by $f$. On the other hand, $\pm f(c_{1})$ is identity on $\Sigma \setminus nbd(c_{1})$.  
\end{defn}

\begin{figure}[h]

\includegraphics[width=10cm,clip]{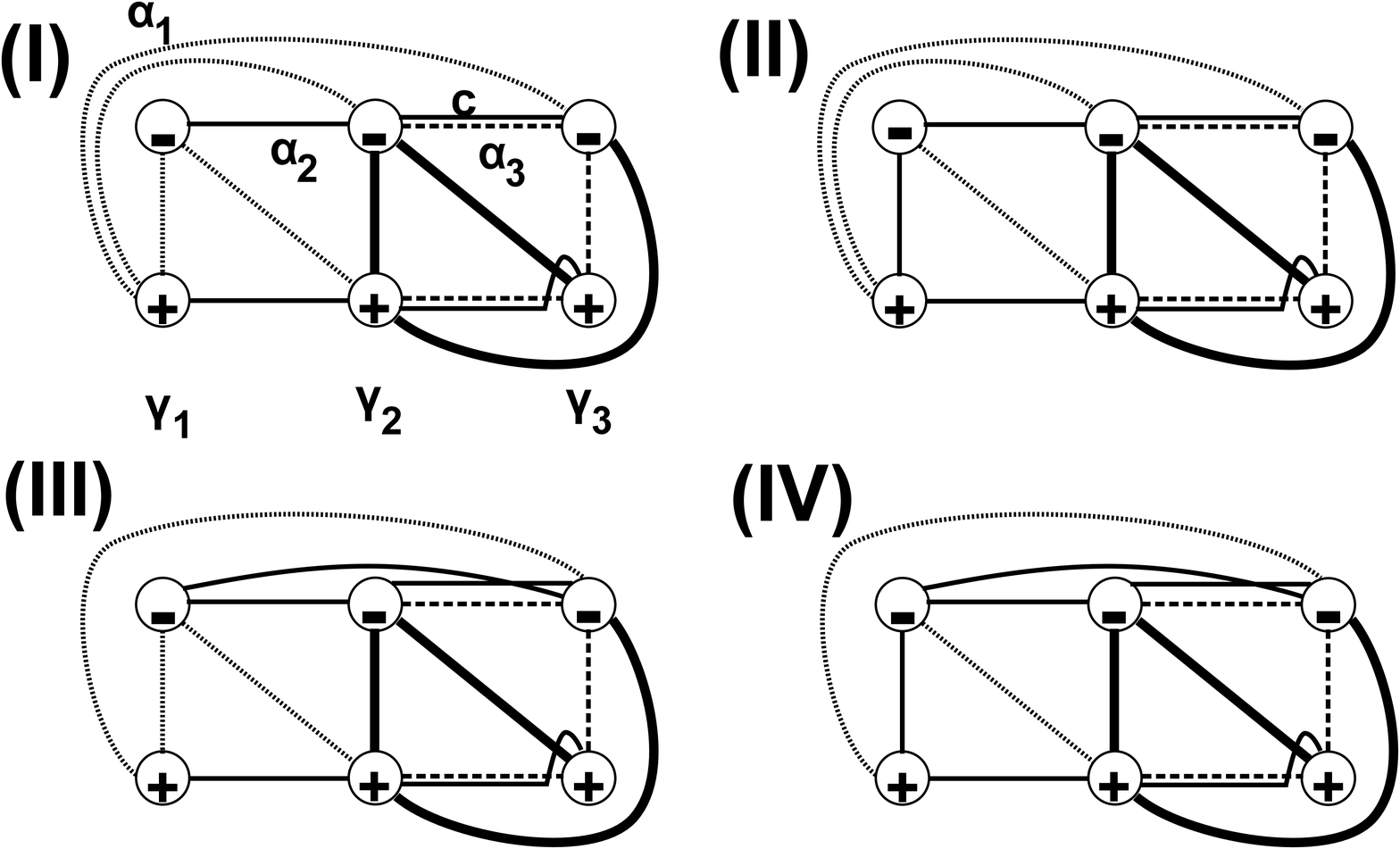} \\
\caption{}\label{4case2}
\end{figure}

\begin{figure}[h]

\includegraphics[width=6cm,clip]{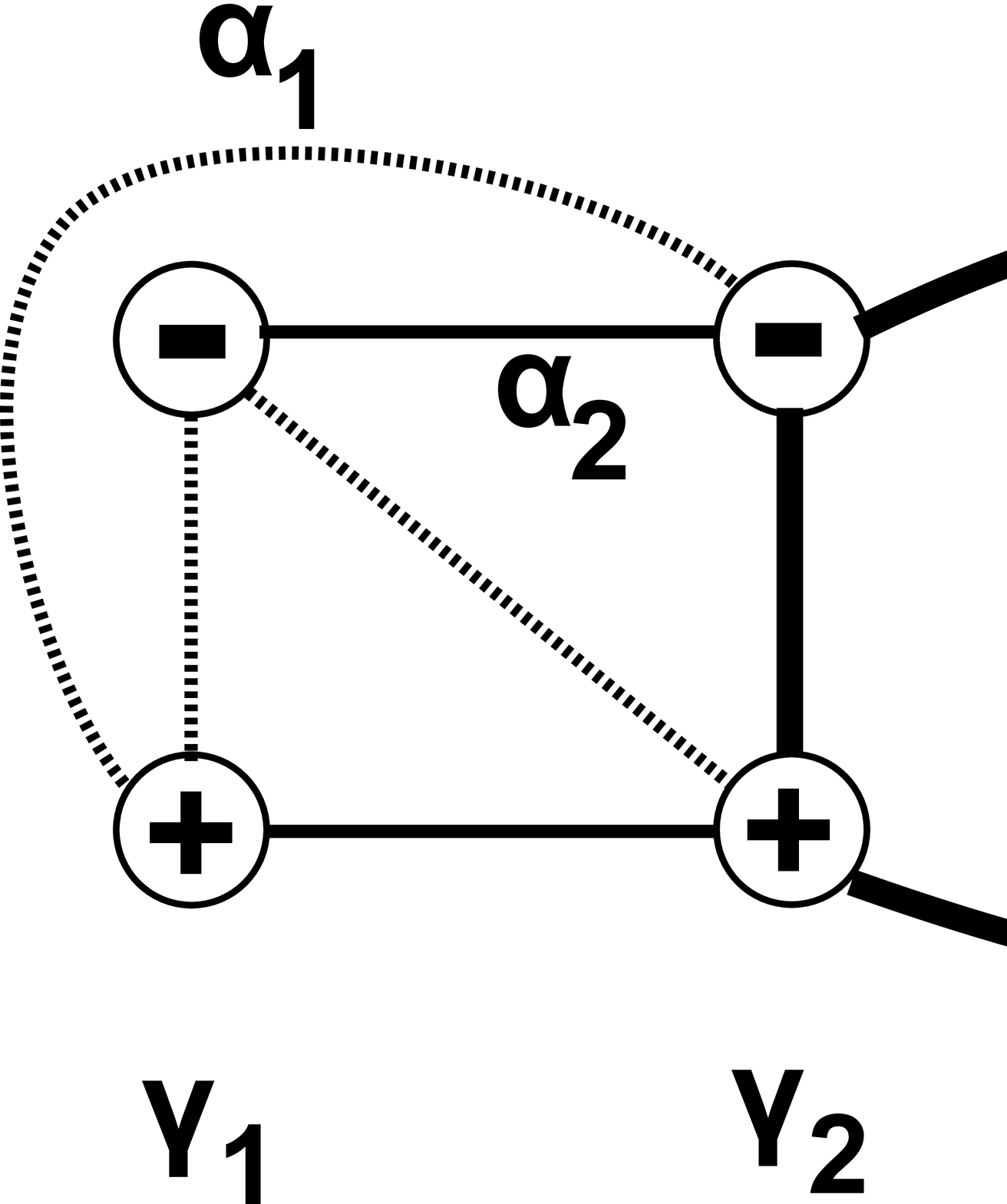} \\
\caption{}\label{4case3}
\end{figure}

\begin{figure}[h]

\includegraphics[width=6cm,clip]{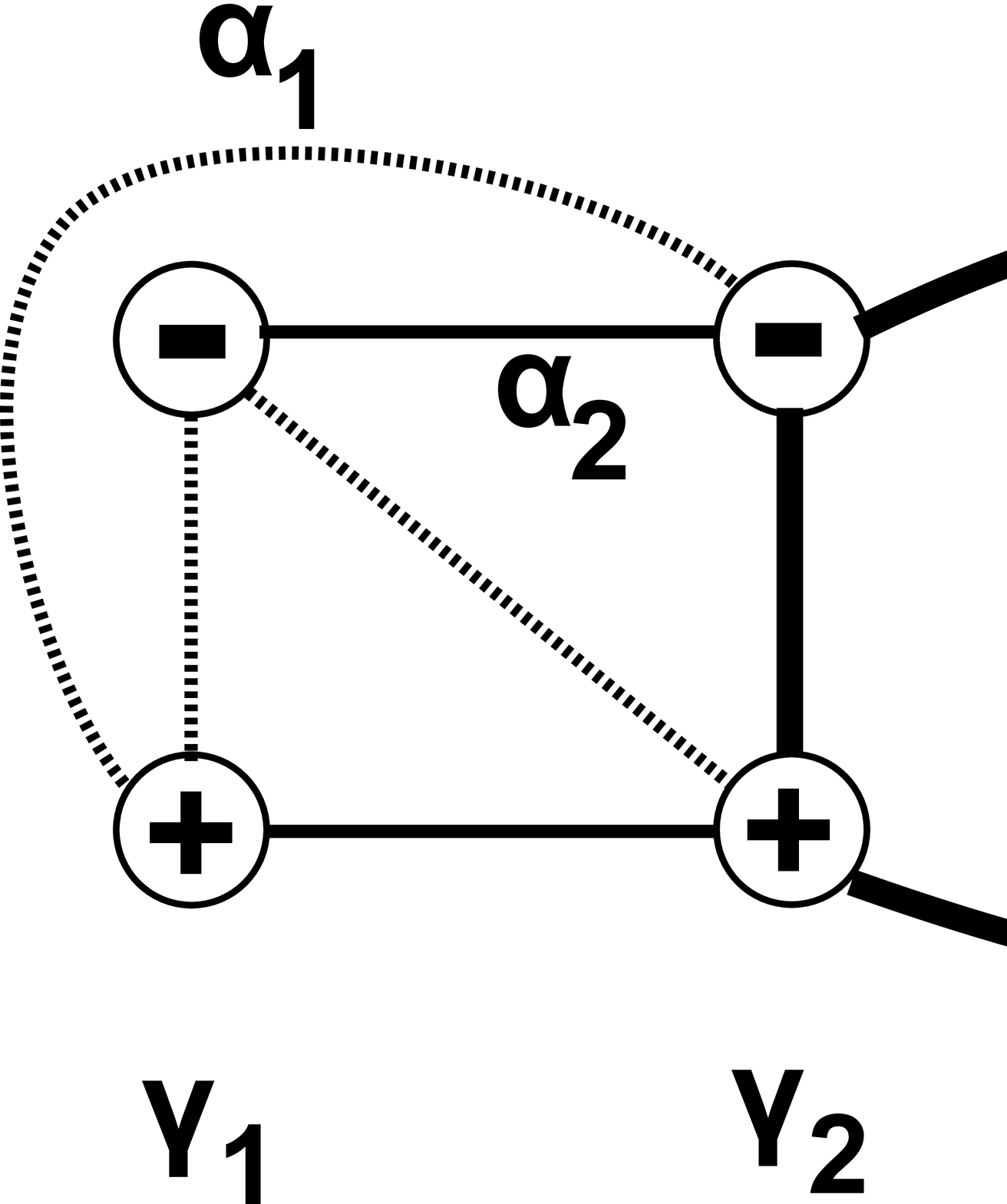} \\
\caption{}\label{4case4}
\end{figure}

Let $(\Sigma,\alpha,\beta)$ be an $A'_{3}$-strong diagram of type 3-(I), 3-(II), 3-(III) or 3-(IV).
Let $c$ be a simple closed curve on $\Sigma$ which intersects $\Gamma(+2,+3)$, $\Gamma(+3,-2)$, $\gamma_{2} $ and $\gamma_{3}$ as in Figure \ref{4case2}. We perform a positive Dehn twist along $c$. 
Then, $\alpha$-arcs are changed as in Figure \ref{4case3}. Note that $\gamma_{3}$ intersects only $\alpha_{3}$. Thus, $\alpha_{3}$-curves become a surgery framing of a unknot $C_{3}$ in $U_{\alpha}$. We describe the slope precisely later. Let $\alpha'_{3}$ be the simple closed curve which intersects only $\gamma_{3}$ at one point (see Figure \ref{4case4}).


\begin{figure}[h]

\includegraphics[width=10cm,clip]{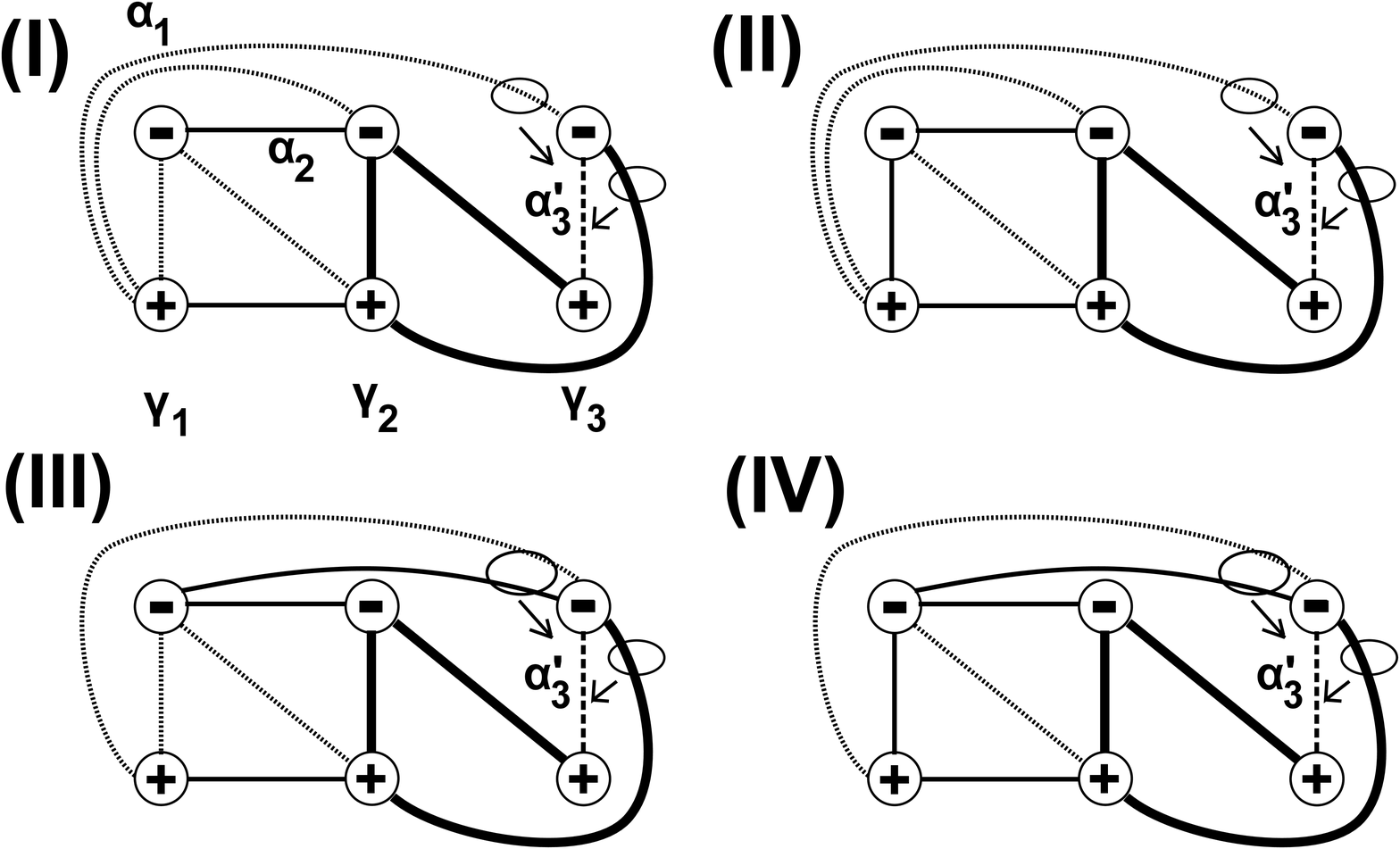} \\
\caption{}\label{4case2-}
\end{figure}

\begin{figure}[h]

\includegraphics[width=5cm,clip]{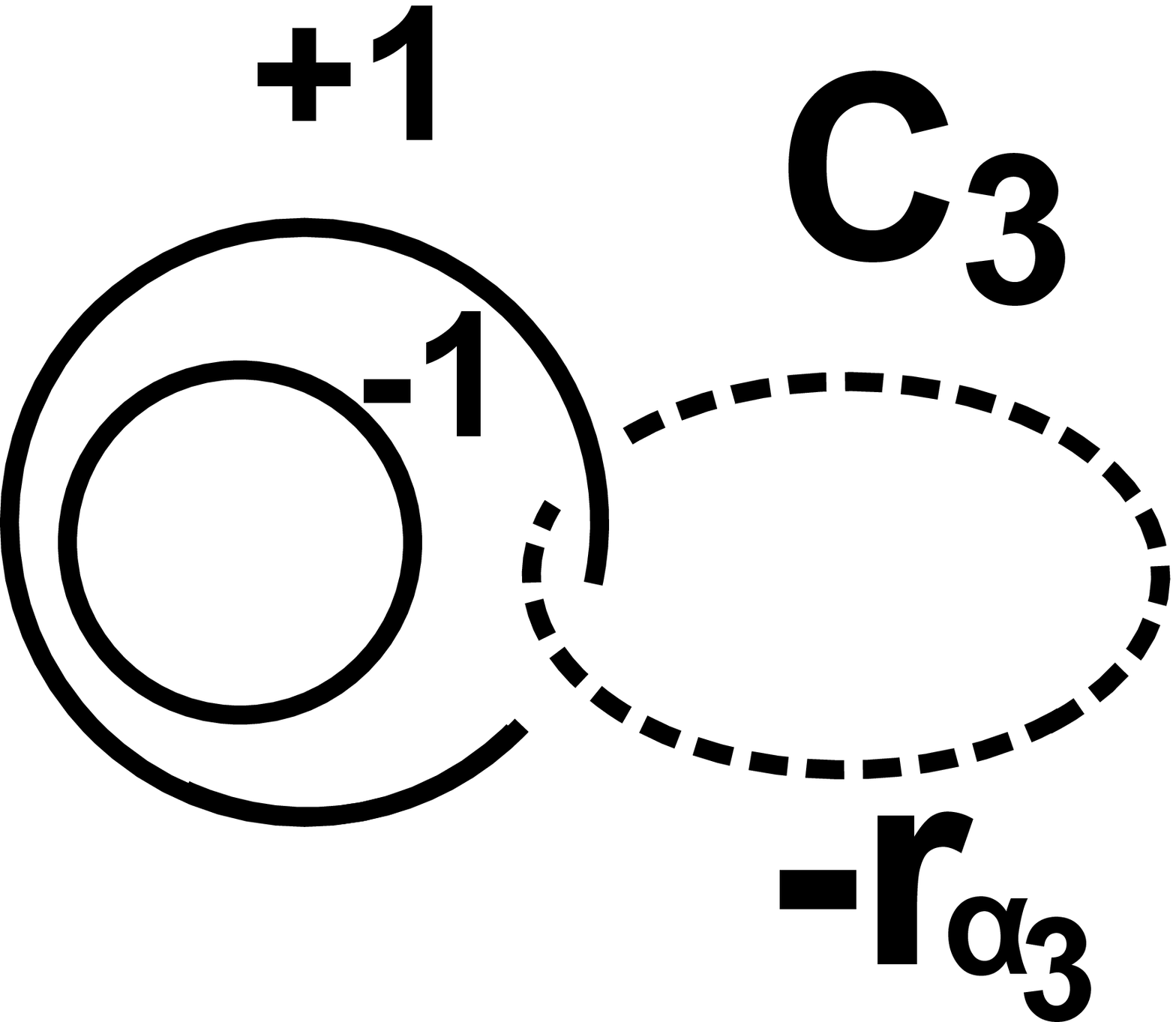} \\
\caption{}\label{type23-0}
\end{figure}

Next, we consider the simple closed curve $c$ as above again. We perform a negative Dehn twist along $c$ (see Figure \ref{4case2}). Since $\alpha'_{3}$ still intersects $\gamma_{3}$ at only one point, we can transform the diagram by a handle-slide $(\Gamma(+1,-3) \text{ and } \Gamma(+2,-3)) \leadsto \alpha_{3}$, where $(\Gamma(+1,-3) \leadsto \alpha_{3}$ means the collection of the handle slides $\gamma \leadsto \alpha_{3}$ and $\gamma \in \Gamma(+1,-3)$ (see Figure \ref{4case2-}).
These operations correspond to the following link (see Figure \ref{type23-0}).

\begin{figure}[h]

\includegraphics[width=10cm,clip]{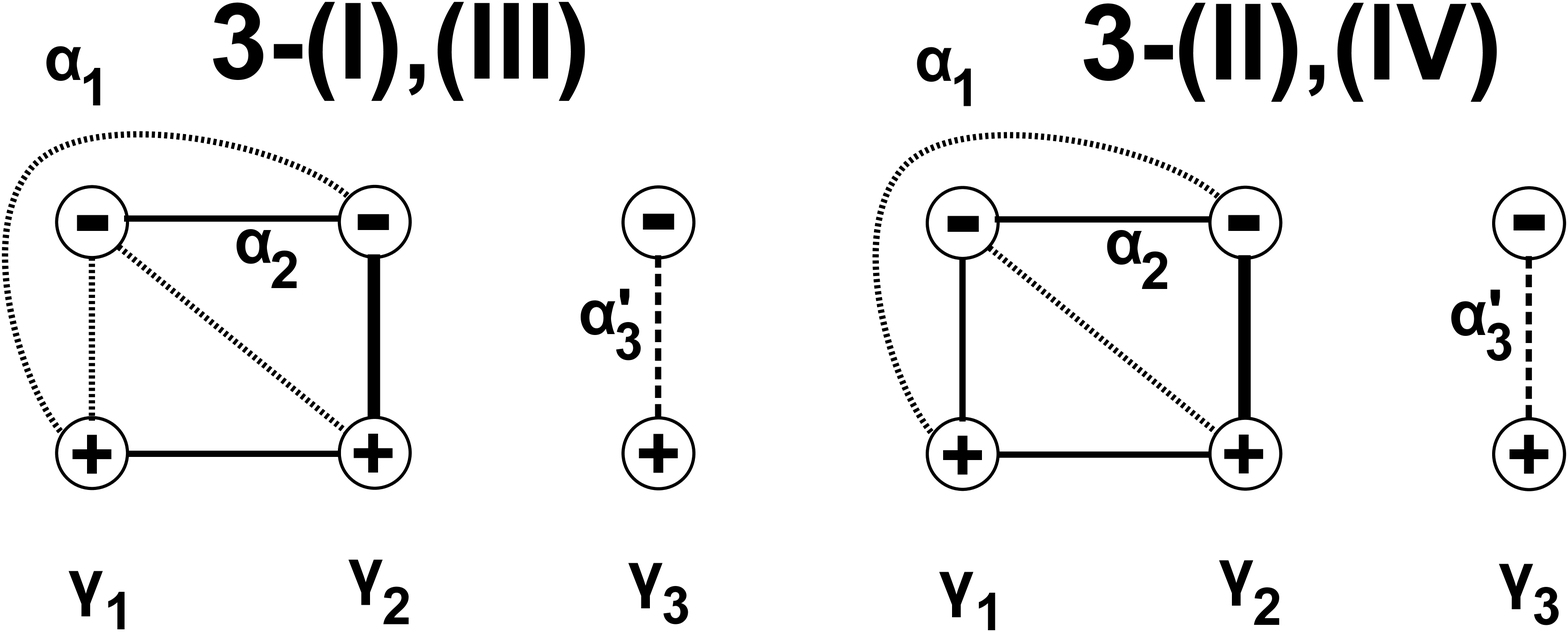} \\
\caption{}\label{4case2-new}
\end{figure}

The Figure \ref{4case2-new} implies that type 3-(I) and 3-(III) induces the same diagrams, and type 3-(II) and 3-(IV) induces the same diagrams. We call these diagrams 3-(I),(III) and 3-(II),(IV).
Let us consider the attaching circles $\alpha' = (\alpha_{1}, \alpha_{2},\alpha'_{3})$ and $(\gamma_{1}, \gamma_{2}, \gamma_{3})$ in Figure \ref{4case2-new}. 

Define another attaching circles $(\delta_{1},\delta_{2}, \delta_{3})$ in each diagram which satisfy the following conditions.
\begin{itemize}
\item $\#(\delta_{i} \cap \gamma_{j}) = \delta_{ij}$ for any $(i,j)$ (thus, $(\Sigma,\delta,\gamma)$ represents $S^3$),
\item If in the case of 3-(I),(III), $\delta_{1}$ intersects $\alpha_{1}$ and does not intersect $\alpha_{2}$.
\item If in the case of 3-(II),(IV), $\delta_{1}$ intersects $\alpha_{2}$ and does not intersect $\alpha_{1}$.
\item $\delta_{2}$ intersects $\alpha$-curves at $N$ points, where $N < \#(\Gamma(+2,-2))+ \#(\Gamma(+2,-1))+\#(\Gamma(+2,+1))$.
\item $\delta_{3} = \alpha'_{3}$
\end{itemize}

Let us denote $D_{1}$ and $D_{3}$ be properly embedded disks in $U_{\delta}$ such that $\partial{D_{1}} = \delta_{1}$ and $\partial{D_{3}} = \delta_{3}$.
If we cut the handlebody $U_{\delta}$ along $D_{1}$ and $D_{3}$, then $U_{\delta} \setminus (D_{1} \cup D_{3})$ becomes a solid torus. Let us denote the core of the solid torus by $C_{2}$. Then,
\begin{itemize}
\item $\alpha_{2}$ becomes the surgery framing of $C_{2}$ in the case of 3-(I),(III), and
\item $\alpha_{1}$ becomes the surgery framing of $C_{2}$ in the case of 3-(II),(IV).
\end{itemize}

On the other hand, it is easy to see that there exists a  $C_{1}$ in $U_{\delta}$ such that 

\begin{itemize}
\item $\alpha_{1}$ become the surgery framings of $C_{1}$ in the case of 3-(I),(III), and
\item $\alpha_{2}$ become the surgery framings of $C_{1}$ in the case of 3-(II),(IV).
\end{itemize}
Note that $C_{1}$ becomes a torus knot in $S^{3}$ in genaral. We describe these slopes later. 

\begin{figure}[h]

\includegraphics[width=7cm,clip]{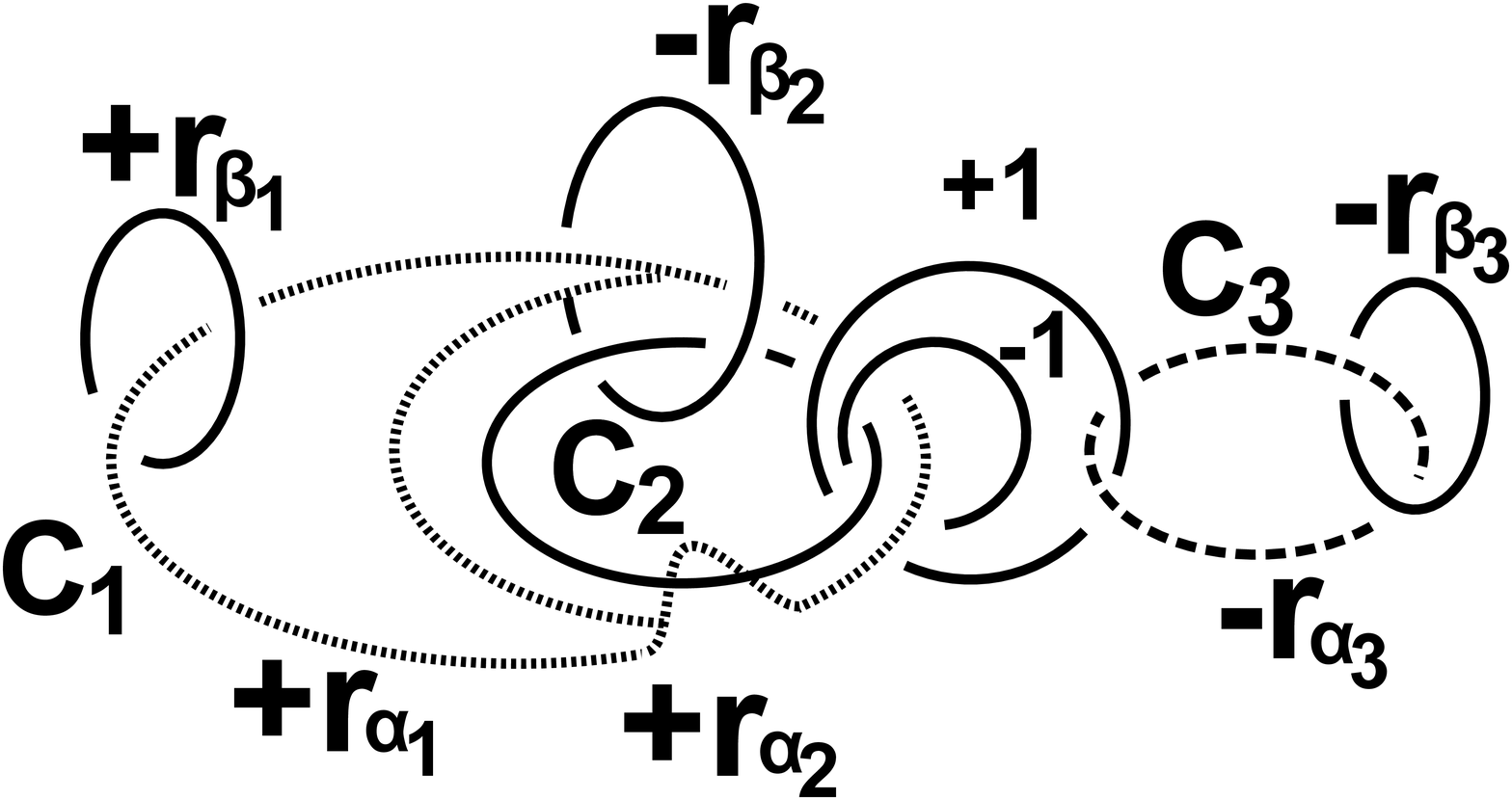} \\
\caption{3-(I),(III)}\label{type23-1}
\end{figure}

\begin{figure}[h]

\includegraphics[width=7cm,clip]{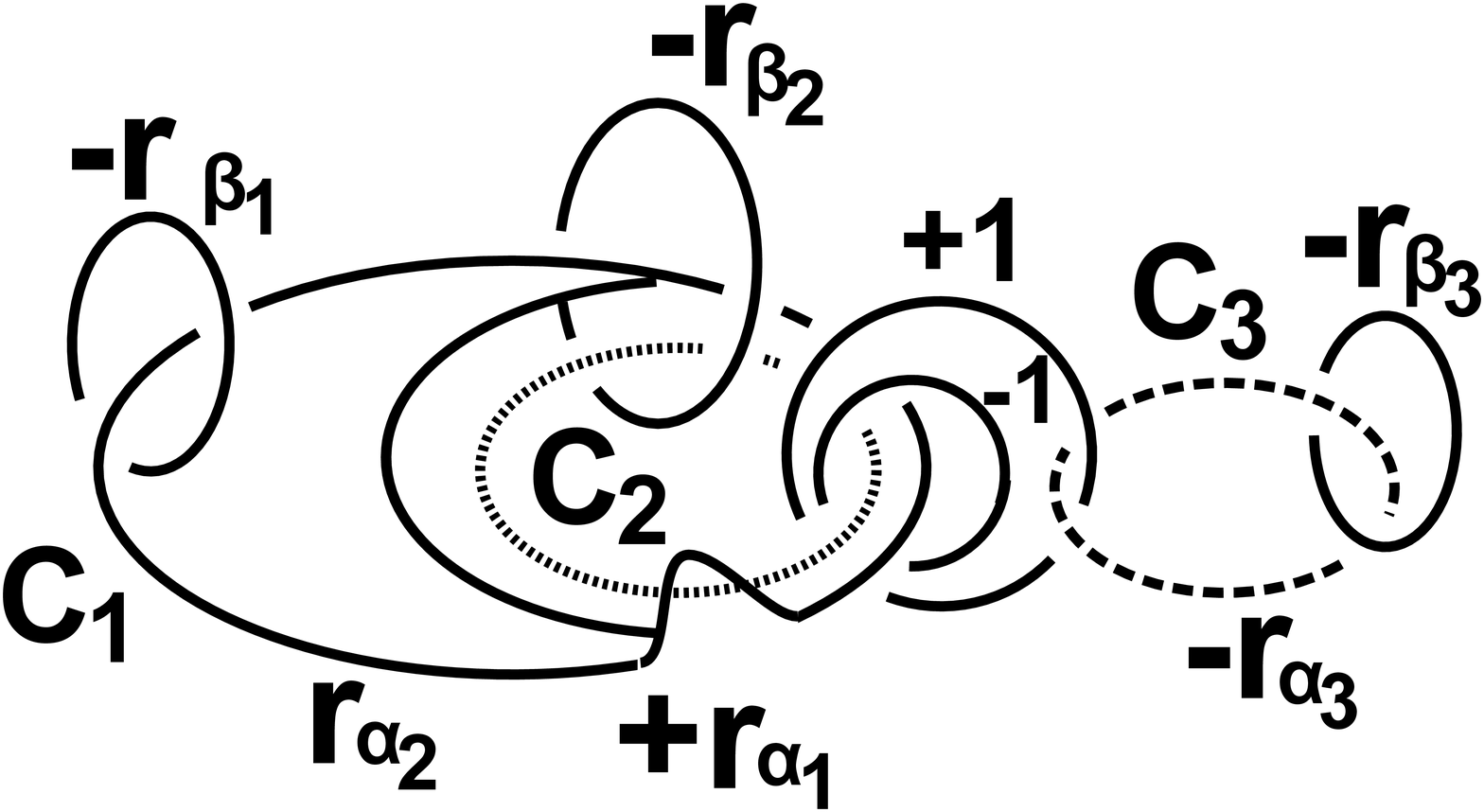} \\
\caption{3-(II),(IV)}\label{type14-1}
\end{figure}

As a result, a Heegaard diagrams of each type is represented by a surgery of $S^3$ along some link (see Figure \ref{type23-1} and \ref{type14-1}). Let us denote the framed link induced from 3-(I),(III) by $L(13)$ and the framed link induced from 3-(II),(IV) by $L(24)$. Moreover, let us denote the manifolds induced from these links $L(13)$ and $L(24)$ by $M(13)$ and $M(24)$.
Now we denote the framing of these links shortly. Let $r_{\alpha_{i}}$ be the framing of $C_{j}$ corresponding to $\alpha_{i}$ for any $(i,j)$. Let $r_{\beta_{j}}$ be the framing of $K_{j}$ corresponding to $\beta_{j}$ for any $j$.

\subsection{Determination of manifolds for $g = 3$}\label{determ}

In this subsection, we finish to prove Theorem \ref{classify}. Recall there are two cases to be considered.

\underline{3-(I),(III)}
Let $(\Sigma, \alpha, \beta)$ be an $A'_{3}$-strong diagram of type 3-(I),(III) (see Figure \ref{4case2-new}).
Recall that $r_{\alpha_{i}}$ is the framing of $C_{i}$ corresponding to $\alpha_{i}$ for any $i$, and $r_{\beta_{j}}$ is the framing of $K_{j}$ corresponding to $\beta_{j}$ for any $j$.

First, we write these slopes concretely.
\begin{itemize}
\item $r_{\alpha_{1}} = + p_{1}/q_{1} +p_{2} + q_{2} -1$, where 
	\begin{itemize}
	\item $p_{1} = \#(\alpha_{1} \cap \gamma_{1})$, 
	\item $q_{1} = \#\Gamma(+1, -2)$,
	\item $p_{2} = \#(\alpha_{1} \cap \gamma_{2})/q_{1}$. 
	\item $q_{2} = \#(\alpha_{1} \cap \delta_{2})$.
	\end{itemize}
\item $r_{\alpha_{2}} = +p'_{2}/q'_{2}$, where
	\begin{itemize}
	\item $p'_{2} = \#(\alpha_{2} \cap \gamma_{2})$,
	\item $q'_{2} = \#(\alpha_{2} \cap \delta_{2})$.
	\end{itemize}
\item $r_{\alpha_{3}} = - p_{3}/q_{3}$, where
	\begin{itemize}
	\item $p_{3} = \#(\alpha_{3} \cap \gamma_{3})$,
	\item $q_{3} = \#(\alpha_{3} \cap \delta_{3})$.
	\end{itemize}
\item $r_{\beta_{1}} = + \#(\beta_{1} \cap \delta_{1})/\#(\beta_{1} \cap \gamma_{1})$.
\item $r_{\beta_{2}} = - \#(\beta_{2} \cap \delta_{2})/\#(\beta_{1} \cap \gamma_{1})$.
\item $r_{\beta_{3}} = - \#(\beta_{3} \cap \delta_{3})/\#(\beta_{1} \cap \gamma_{1})$.
\end{itemize}

Then, $C_{1}$ becomes the $(p_{2},q_{2})$-torus knot and it is linking with $C_{2}$, where the framing of the $(p_{2},q_{2})$-torus becomes the integer $p_{2} + q_{2} -1$.
We also find that, by easy obsevations,
\begin{itemize}
\item $+p_{1}/q_{1} > +1$,
\item $r_{\alpha_{2}}> +1$,
\item $r_{\beta_{1}} > +1$,
\item $r_{\beta_{2}} <-1$,
\item $r_{\beta_{3}} <-1$.
\end{itemize}

Next, we describe the relation between $p_{2}/q_{2}$ and $p'_{2}/q'_{2}$ precisely.
We can represent $p_{2}/q_{2}$ as the continuous fraction expansion as follows.

\begin{equation}
\frac{p_{2}}{q_{2}}  =  k_{1} + \cfrac{1}{ k_{2} + 
						\cfrac{1}{ \cdots +
						\cfrac{1}{ k_{n-1} +
						\cfrac{1}{ k_{n}}}}}, \text{ where } k_{i} \ge 1 \text{ and } k_{n} \ge 2.
\end{equation}

By using these integers, we put a new rational number $R(p_{2},q_{2},p'_{2},q'_{2})$ as follows.

\begin{equation}
R(p_{2},q_{2},p'_{2},q'_{2}) = -( k_{n} + \cfrac{1}{ k_{n-1} + 
						\cfrac{1}{ \cdots +
						\cfrac{1}{ k_{2} +
						\cfrac{1}{ k_{1} -
						\cfrac{p'_{2}}{q'_{2}}}}}}).
\end{equation}

Then, we can prove the following claim.

\begin{clm}\label{clm3}
Let $p_{2}/q_{2}$, $p'_{2}/q'_{2})$ and $R(p_{2},q_{2},p'_{2},q'_{2})$ be the rational numbers defined as above. Then, we get that
\begin{itemize}
\item $R(p_{2},q_{2},p'_{2},q'_{2}) > 0$ if $n$ is odd, and
\item $-1 < R(p_{2},q_{2},p'_{2},q'_{2}) < 0$ if $n$ is even.
\end{itemize}
\end{clm}

\begin{prf}
Originaly, these two rational numbers $p_{2}/q_{2}$ and $p'_{2}/q'_{2})$ come from the slopes of $\alpha_{1}$ and $\alpha_{2}$ near $\gamma_{2}$. Thus, we study about the neighborhood of $\beta_{2}$ more precisely.

Recall the neighborhood of $\gamma_{2}$ looks as in Figure \ref{Rnbd}, where we can define $x = (m_{1},n_{1},m_{2},\cdots, n_{\tilde{k}-1},m_{\tilde{k}})$ as the sequence of the number of intersection points induced from $\Gamma(+2,-2)$. $n$ comes from  $\gamma_{2} \cap \alpha_{1}$ and $m$ comes from $\gamma_{2} \cap \alpha_{2}$. 

\begin{figure}[h]

\includegraphics[width=7cm,clip]{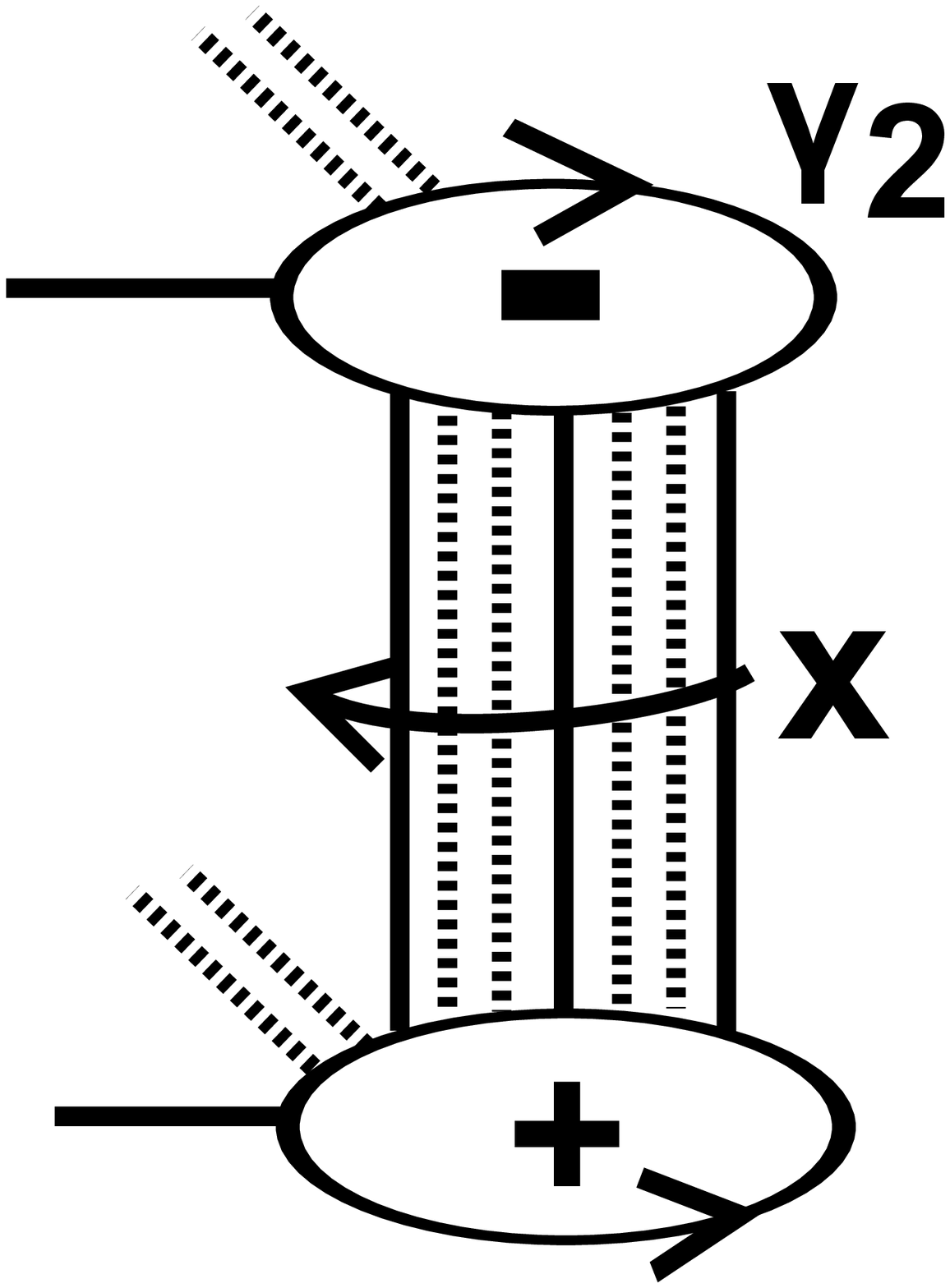} \\
\caption{}\label{Rnbd}
\end{figure}

Note that it is enough to consider the case when $m_{1} = m_{\tilde{k}}$ is not zero. Otherwise, $\#\Gamma(+2,+1) = 1$ and we can change the diagram by handle-slides so that $n_{1} = n_{\tilde{k}-1}$ is equal to zero in the new diagram (see Figure \ref{nbd3}). 

We put $d_{1} = \# \Gamma(+2,+1)$ and $d_{2} = m_{1}= m_{\tilde{k}}$. Moreover, let $1 \le \bar{p}_{2} \le p_{2}$ and $1\le \bar{q}_{2} \le q_{2}$ be positive coprime integers so that $\bar{p}_{2} q_{2} = \bar{q}_{2} p_{2} + 1$.

Then, $p'_{2} / q'_{2}$ can be written by these integers as follows. 
$$p'_{2} / q'_{2} = \frac{\bar{p}_{2}(d_{1}+d_{2})+(p_{2}-\bar{p}_{2})d_{2}}{\bar{q}_{2}(d_{1}+d_{2})+(q_{2}-\bar{q}_{2})d_{2}} = \frac{\bar{p}_{2}+p_{2}z}{\bar{q}_{2}+q_{2}z}, \text{ where } z = d_{2}/d_{1}.$$

Actually, $\bar{p}_{2}$ represents the number of the integer $d_{1}+d_{2}$ among the sequence $(m_{1},m_{2},\cdots,m_{\tilde{k}})$ and $\bar{q}_{2}$ represents the number of the integer $d_{1}+d_{2}$ among the sequence $(m_{1},m_{2},\cdots,m_{q_{2}})$. It is easy to see that these integers satisfies $\bar{p}_{2} q_{2} = \bar{q}_{2} p_{2} + 1$. Then, we find that $p'_{2} / q'_{2}$ is written as above.

Since $0 < z < +\infty$ we get that $p_{2} / q_{2}<p'_{2} / q'_{2}<\bar{p}_{2} / \bar{q}_{2}$. 

Moreover, we can prove that $\bar{p}_{2} / \bar{q}_{2}$ can be written precisely as follows.
$$ \bar{p}_{2} / \bar{q}_{2} =  k_{1} + \cfrac{1}{ k_{2} + 
						\cfrac{1}{ \cdots +
						\cfrac{1}{ k_{n-1}}}} , \text{ if } n \text{ is odd.}$$
						
$$ \bar{p}_{2} / \bar{q}_{2} =  k_{1} + \cfrac{1}{ k_{2} + 
						\cfrac{1}{ \cdots +
						\cfrac{1}{ k_{n} - 1}}} ,\text{ if } n \text{ is even.}$$
This equation comes from the Euclidean algorithm. 
Thus, we also find that the length of the continuous fraction expansion of $p'_{2} / q'_{2}$ is $n$ or grater than $n$. In particular, $R(p_{2},q_{2},p'_{2},q'_{2})$ is well-defined.

Finally, we conclude that 
$$R(p_{2},q_{2},p'_{2},q'_{2}) = - ( k_{n} + \cfrac{1}{ k_{n-1} + 
						\cfrac{1}{ \cdots +
						\cfrac{1}{ k_{2} +
						\cfrac{1}{ k_{1} - p'_{2} / q'_{2}}}}} )   > 0 ,$$ if $n$ is odd, and
$$0 > R(p_{2},q_{2},p'_{2},q'_{2}) > - ( k_{n} + \cfrac{1}{ k_{n-1} + 
						\cfrac{1}{ \cdots +
						\cfrac{1}{ k_{2} +
						\cfrac{1}{ k_{1} - \bar{p}_{2} / \bar{q}_{2}}}}} ) = -1,$$ if $n$ is even. 
\qed
\end{prf}

We consider the framed link $L(13)$ again. Our goal is to prove that $M(13)$ is in $\mathcal{L}_{\overline{\rm{Brm}}}$.

\begin{figure}[h]

\includegraphics[width=7cm,clip]{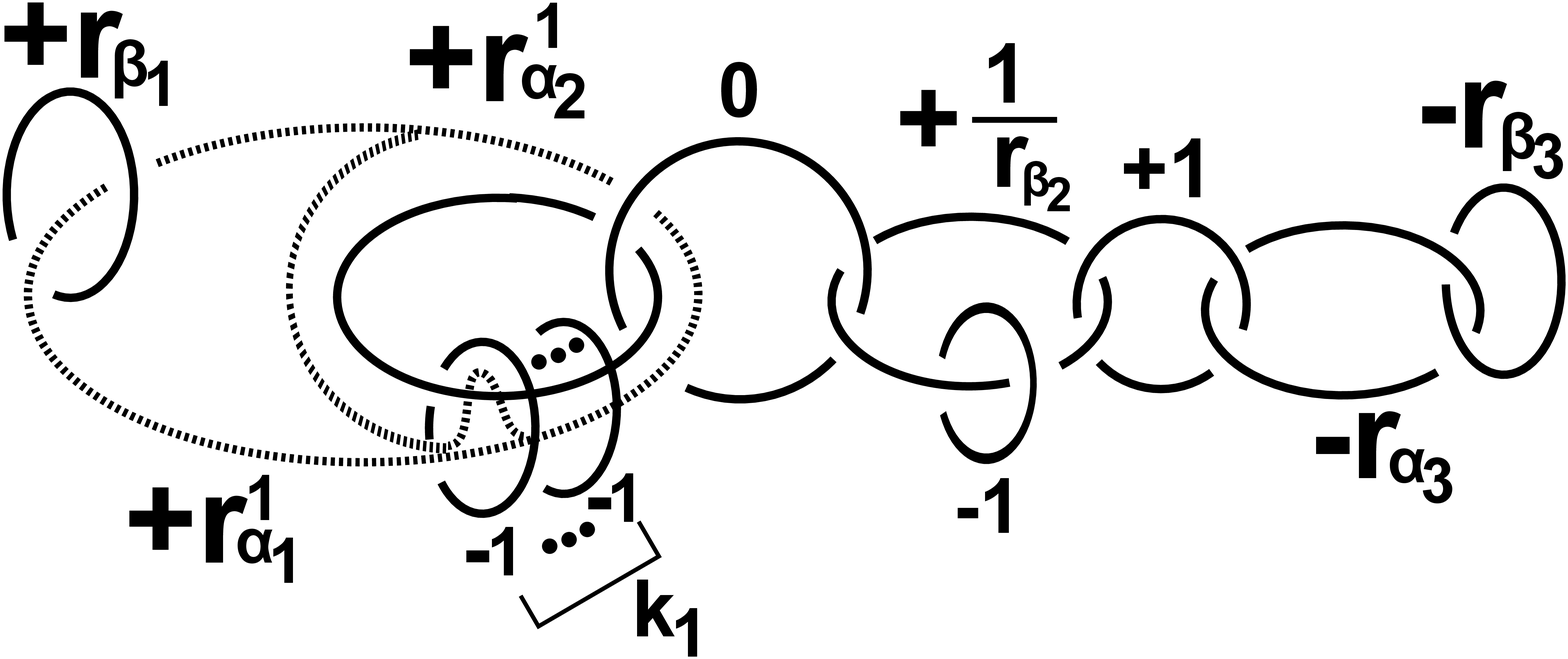} \\
\caption{}\label{type23-2}
\end{figure}

We perform the blow up operations finitely many times so that the new framed link consists of only unknots.
First, add $k_{1}$ unknots with framing $-1$ near $C_{2}$ and $C_{1}$ as in Figure \ref{type23-2}. After that, the slopes of $C_{1}$ and $C_{2}$ are changed as follows. Denote the new slopes by $r^{1}_{\alpha_{1}}$ and $r^{1}_{\alpha_{2}}$.

\begin{itemize}
\item $r_{\alpha_{1}} \leadsto r^{1}_{\alpha_{1}} = r_{\alpha_{1}} - k_{1}q^{0}_{2}$, where $q^{0}_{2} = q_{2}$,
\item $r_{\alpha_{2}} \leadsto r^{1}_{\alpha_{2}} = r_{\alpha_{2}} - k_{1}$,
\item the new knot $C^{1}_{1}$ induced from $C_{1}$ is linking with $C_{2}$ with the slope $p^{1}_{2}/q^{1}_{2}$, where $p^{1}_{2}$ and $q^{1}_{2}$ are positive coprime integer such that $p_{2}/q_{2} = k_{1}+p^{1}_{2}/q^{1}_{2}$.
\end{itemize}

\begin{figure}[h]

\includegraphics[width=7cm,clip]{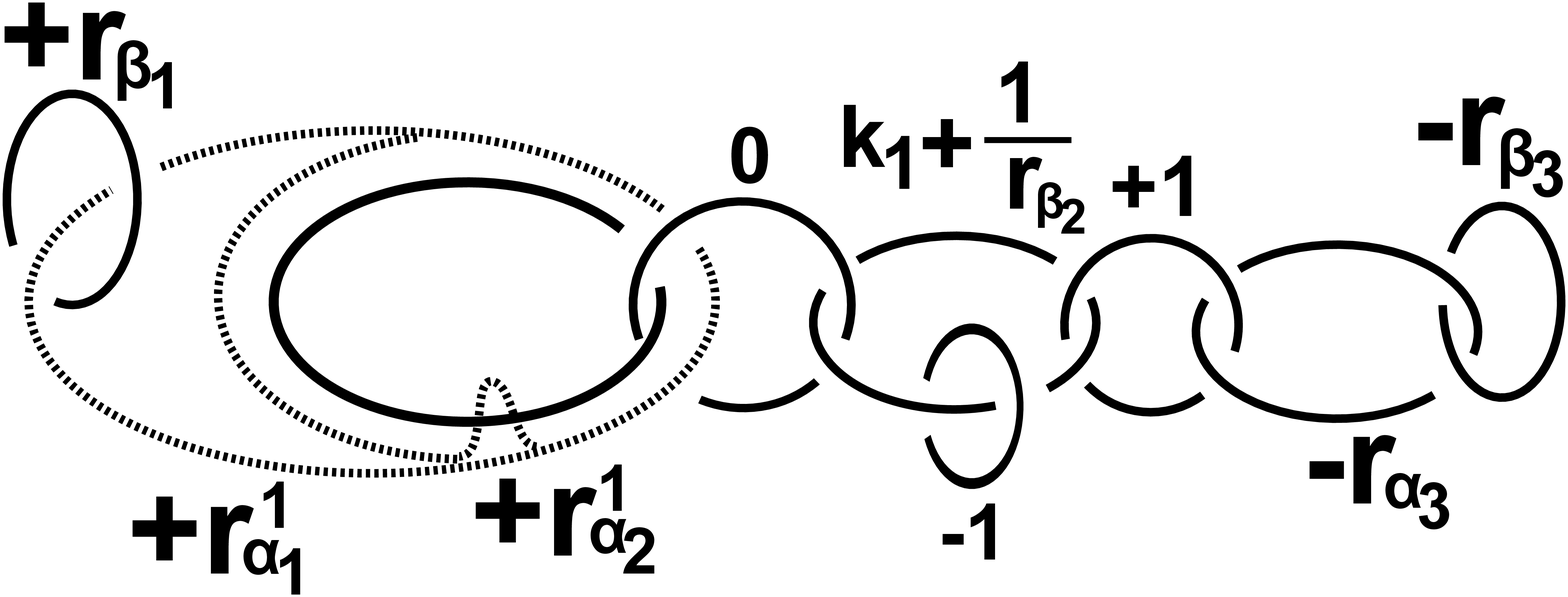} \\
\caption{}\label{type23-2-}
\end{figure}

This new link can be transformed into the following link (see Figure \ref{type23-2-}).

\begin{figure}[h]

\includegraphics[width=7cm,clip]{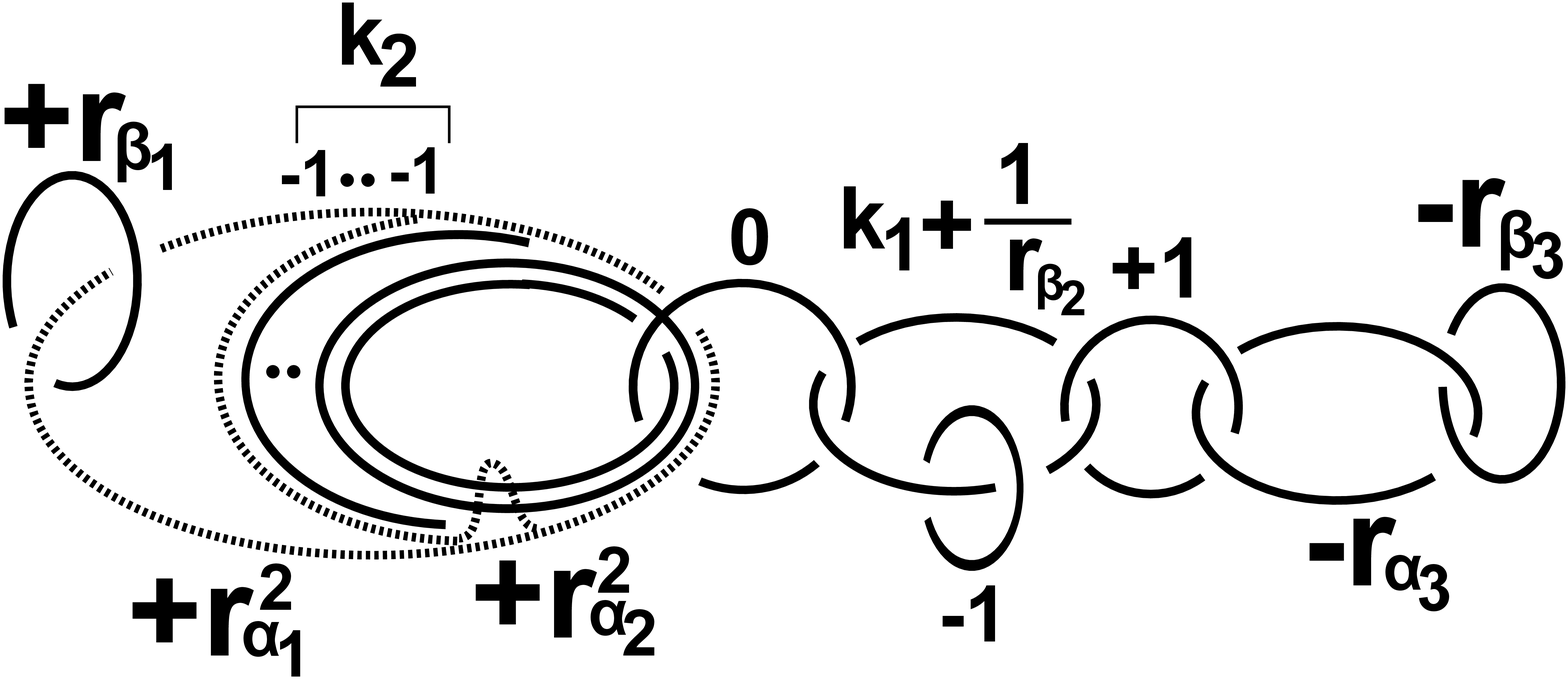} \\
\caption{}\label{type23-3}
\end{figure}

Next, add $k_{2}$ unknots with framing $-1$ isotopic to $C_{2}$ as in Figure \ref{type23-3}. After some Kirby calculus, the slopes of $C_{1}$ and $C_{2}$ are changed as follows.
\begin{itemize}
\item $r^{1}_{\alpha_{1}} \leadsto r^{2}_{\alpha_{1}} = r^{1}_{\alpha_{1}} - k_{2}p^{1}_{2}$,
\item $r^{1}_{\alpha_{2}} \leadsto r^{2}_{\alpha_{2}} = -  \cfrac{1}{ k_{2} + 
						\cfrac{1}{- r^{1}_{\alpha_{2}}}}$,
\item the new $C^{2}_{1}$ is linking with $C_{2}$ with the slope $p^{2}_{2}/q^{2}_{2}$, where $p^{2}_{2}$ and $q^{2}_{2}$ are positive coprime integer such that $p^{1}_{2}/q^{1}_{2} = k_{2}+p^{2}_{2}/q^{2}_{2}$.
\end{itemize}

\begin{figure}[h]

\includegraphics[width=7cm,clip]{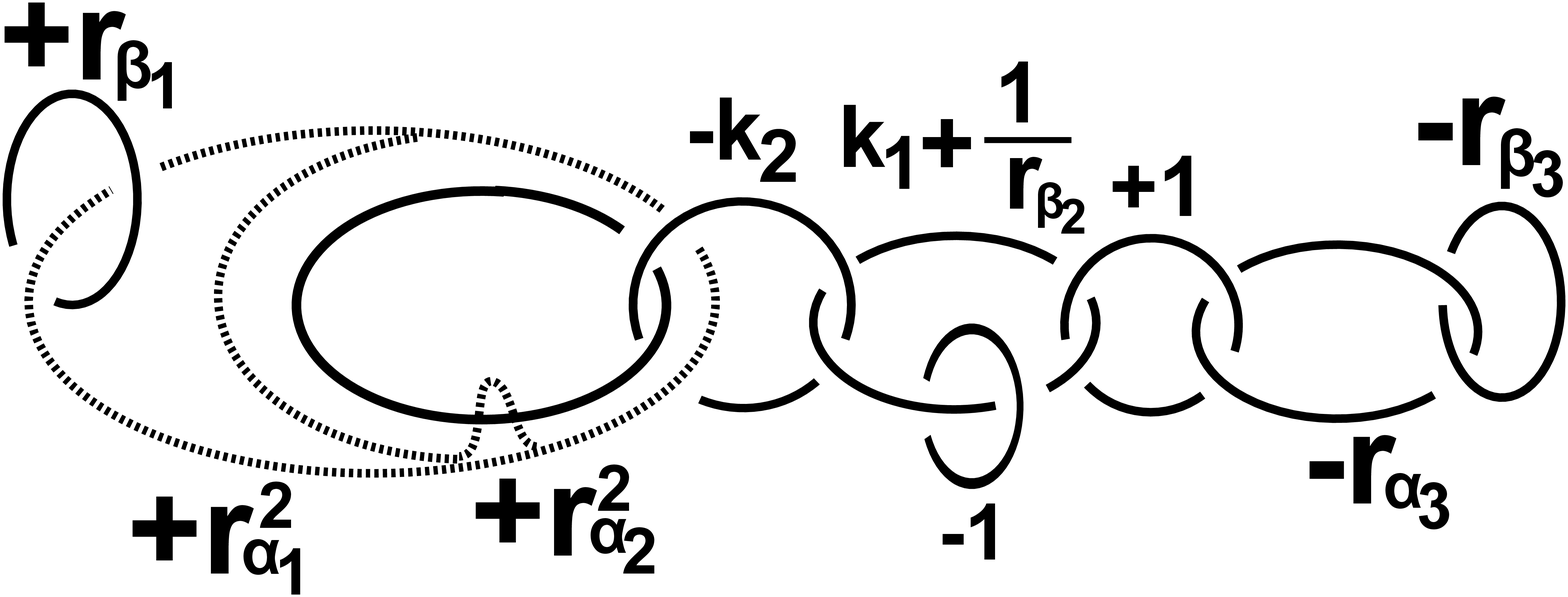} \\
\caption{}\label{type23-3-}
\end{figure}

\begin{figure}[h]

\includegraphics[width=7cm,clip]{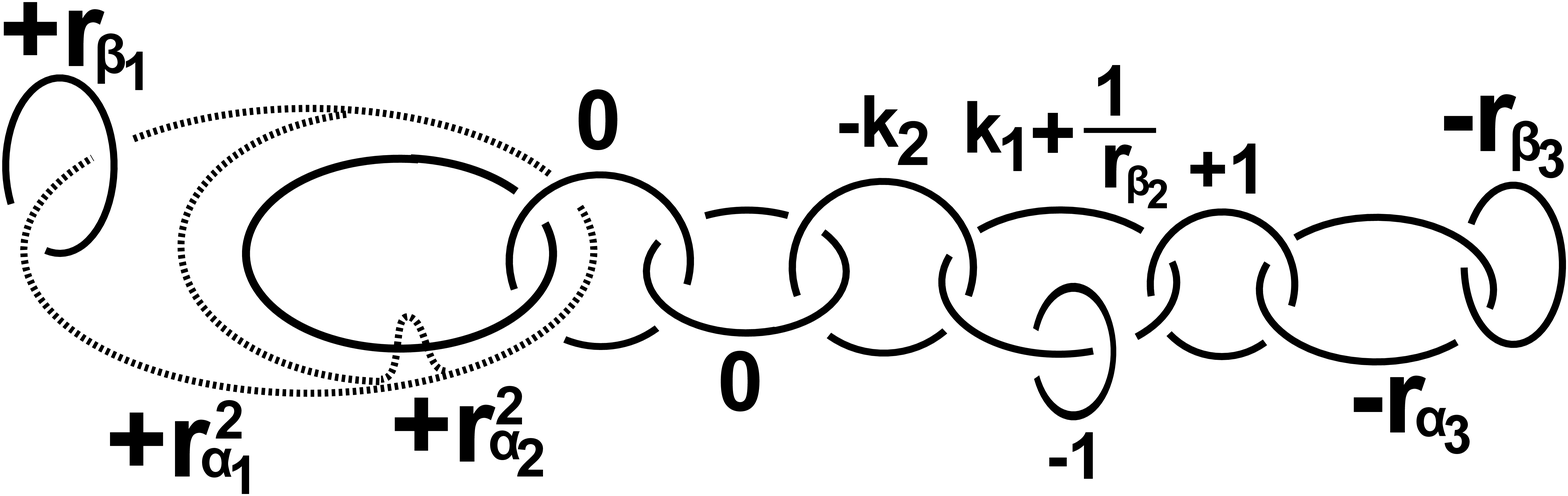} \\
\caption{}\label{type23-3+}
\end{figure}

This new link can be transformed into the following link (see Figure \ref{type23-3+}).

We find that $C^{2}_{1}$ is linking with $C_{2}$ with slope $p^{2}_{2}/q^{2}_{2}$. Thus, we can transform this link again by adding new $k_{3}$ unknots as above. 

\begin{figure}[h]

\includegraphics[width=10cm,clip]{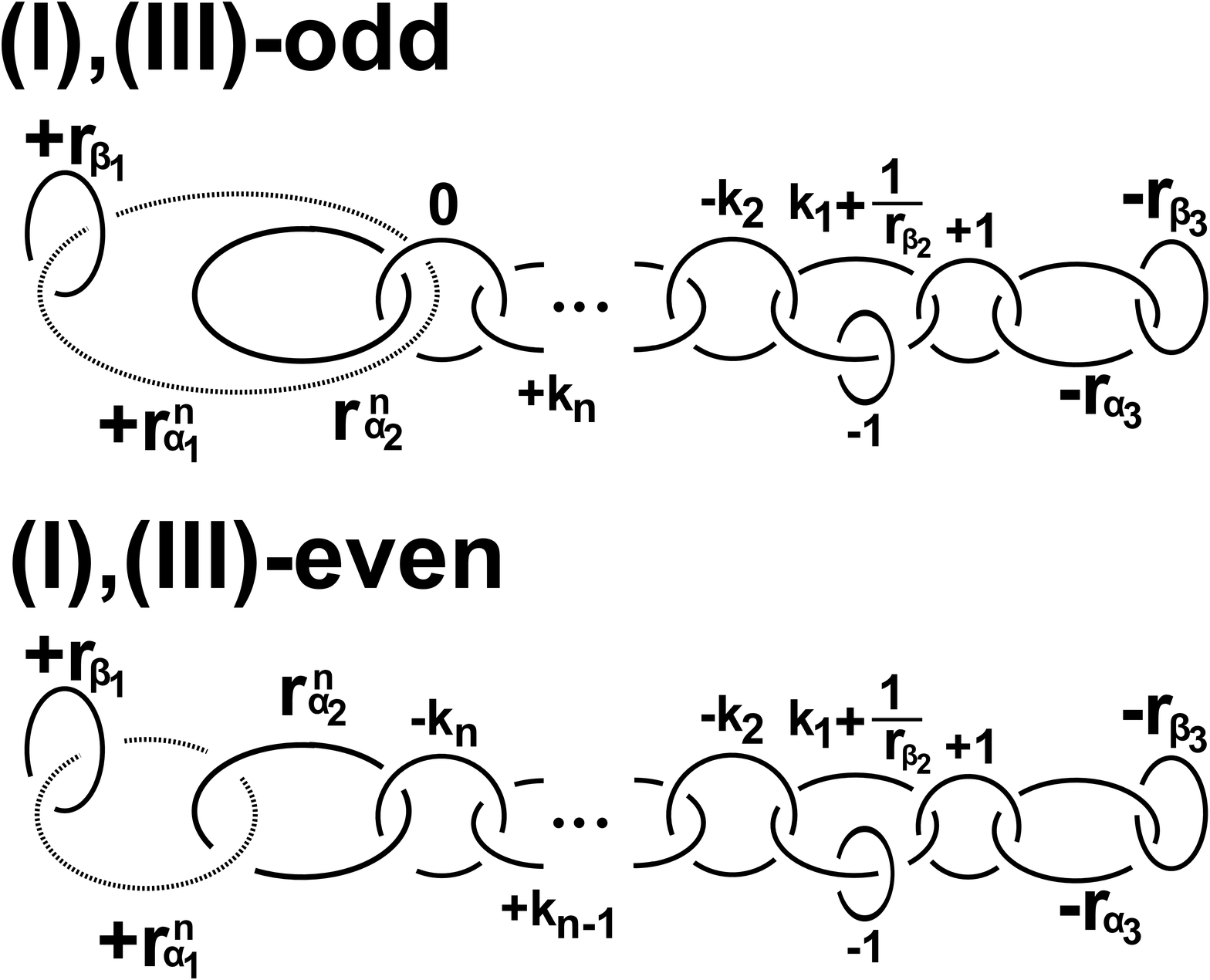} \\
\caption{}\label{type23-f}
\end{figure}

In finitely many steps, we finally get the following framed links (see Figure \ref{type23-f}).

\begin{itemize}
\item $r^{n}_{\alpha_{1}} = r_{\alpha_{1}} - k_{1}q_{2} - k_{2}p^{1}_{2}-k_{3}q^{2}_{2}-k_{4}p^{3}_{2} \cdots = r_{\alpha_{1}} - \sum_{i \text{:odd}} k_{i}q^{i-1}_{2} - \sum_{i \text{:even}} k_{i}p^{i-1}_{2} = r_{\alpha_{1}}-(p+q-1) = +p_{1}/q_{1}$,
\item $r^{n}_{\alpha_{2}} = - ( k_{n} + \cfrac{1}{ k_{n-1} + 
						\cfrac{1}{ \cdots +
						\cfrac{1}{ k_{2} +
						\cfrac{1}{ k_{1} - r_{\alpha_{2}}}}}} ) = R(p_{2},q_{2},p'_{2},q'_{2})$
\item the new $C^{n}_{1}$ is the unknot shown in Figure \ref{type23-f}, where $p^{i}_{2}$ and $q^{i}_{2}$ are positive coprime integer such that $p^{i-1}_{2}/q^{i-1}_{2} = k_{i}+p^{i}_{2}/q^{i}_{2}$.
\end{itemize}

\begin{figure}[h]

\includegraphics[width=7cm,clip]{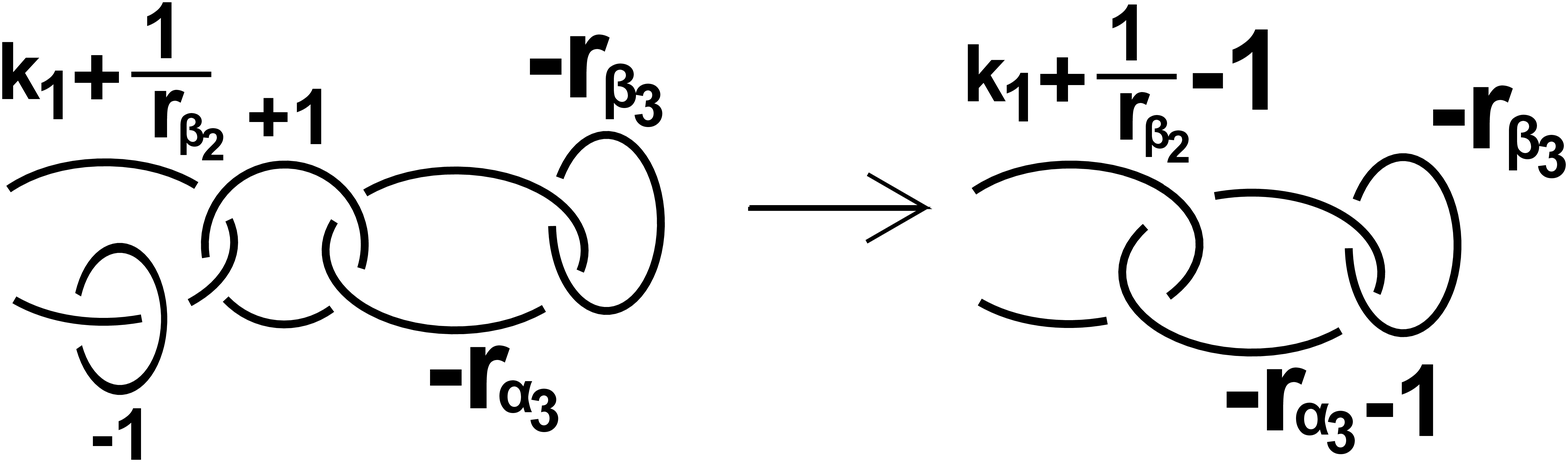} \\
\caption{}\label{type1234-1}
\end{figure}

In each case, the right part of the link can be changed to have alternating weights because $k_{1} + 1/r_{\beta_{2} -1 \ge 0}$, $-r_{\alpha_{3} -1 < 0}$ and $-r_{\beta_{3}}$ (see Figure \ref{type1234-1}).

To prove that the left part of these link have also alternating weights, we can just use claim \ref{clm3}.
Thus, $L(13)$ can be represented by alternatingly weighted unknots, That is, $M(13)$ belongs to $\mathcal{L}_{\overline{\rm{Brm}}}$.

\underline{3-(II),(IV)}

Let $(\Sigma, \alpha, \beta)$ be an $A'_{3}$-strong diagram of type 3-(II),(IV) (see Figure \ref{4case2-new}).
Recall that $r_{\alpha_{1}}$, $r_{\alpha_{2}}$ and $r_{\alpha_{3}}$ are the framings of $C_{2}$, $C_{1}$ and $C_{3}$ corresponding to $r_{\alpha_{1}}$, $r_{\alpha_{2}}$ and $r_{\alpha_{3}}$ respectively, and $r_{\beta_{j}}$ is the framing of $K_{j}$ corresponding to $\beta_{j}$ for any $j$.

First, we also write these slopes concretely.
\begin{itemize}
\item $r_{\alpha_{2}} = - p_{1}/q_{1} +p_{2} + q_{2} -1$, where 
	\begin{itemize}
	\item $p_{1} = \#(\alpha_{2} \cap \gamma_{1})$, 
	\item $q_{1} = \#\Gamma(-1, -2)$,
	\item $p_{2} = \#(\alpha_{2} \cap \gamma_{2})/q_{1}$. 
	\item $q_{2} = \#(\alpha_{2} \cap \delta_{2})$.
	\end{itemize}
\item $r_{\alpha_{1}} = +p'_{2}/q'_{2}$, where
	\begin{itemize}
	\item $p'_{2} = \#(\alpha_{1} \cap \gamma_{2})$,
	\item $q'_{2} = \#(\alpha_{1} \cap \delta_{2})$.
	\end{itemize}
\item $r_{\alpha_{3}} = - p_{3}/q_{3}$, where
	\begin{itemize}
	\item $p_{3} = \#(\alpha_{3} \cap \gamma_{3})$,
	\item $q_{3} = \#(\alpha_{3} \cap \delta_{3})$.
	\end{itemize}
\item $r_{\beta_{1}} = - \#(\beta_{1} \cap \delta_{1})/\#(\beta_{1} \cap \gamma_{1})$.
\item $r_{\beta_{2}} = - \#(\beta_{2} \cap \delta_{2})/\#(\beta_{1} \cap \gamma_{1})$.
\item $r_{\beta_{3}} = - \#(\beta_{3} \cap \delta_{3})/\#(\beta_{1} \cap \gamma_{1})$.
\end{itemize}

Then, $C_{1}$ becomes the $(p_{2},q_{2})$-torus knot and it is linking with $C_{2}$.

\begin{itemize}
\item $-p_{1}/q_{1} < -1$,
\item $r_{\alpha_{1}}> +1$,
\item $r_{\beta_{1}} < -1$,
\item $r_{\beta_{2}} <-1$,
\item $r_{\beta_{3}} <-1$.
\end{itemize}

Next, we describe the relation between $p_{2}/q_{2}$ and $p'_{2}/q'_{2}$ precisely.
We can represent $p_{2}/q_{2}$ as the continuous fraction expansion similarly.

\begin{equation}
\frac{p_{2}}{q_{2}}  =  k_{1} + \cfrac{1}{ k_{2} + 
						\cfrac{1}{ \cdots +
						\cfrac{1}{ k_{n-1} +
						\cfrac{1}{ k_{n}}}}}, \text{ where } k_{i} \ge 1 \text{ and } k_{n} \ge 2.
\end{equation}

We set $R(p_{2},q_{2},p'_{2},q'_{2})$ similarly.

\begin{equation}
R(p_{2},q_{2},p'_{2},q'_{2}) = -( k_{n} + \cfrac{1}{ k_{n-1} + 
						\cfrac{1}{ \cdots +
						\cfrac{1}{ k_{2} +
						\cfrac{1}{ k_{1} -
						\cfrac{p'_{2}}{q'_{2}}}}}}).
\end{equation}

Then, we can prove the following claim.

\begin{clm}\label{clm4}
Let $p_{2}/q_{2}$, $p'_{2}/q'_{2})$ and $R(p_{2},q_{2},p'_{2},q'_{2})$ be the rational numbers defined as above. Then, we get that
\begin{itemize}
\item $-1 < R(p_{2},q_{2},p'_{2},q'_{2}) < 0$ if $n$ is odd, and
\item $R(p_{2},q_{2},p'_{2},q'_{2}) > 0$ if $n$ is even.
\end{itemize}
\end{clm}

\begin{prf}
We can prove this claim similarly to the proof of Claim \ref{clm3}.

Actually, if we define $x = (n_{1},m_{1},n_{2},\cdots, m_{\tilde{k}-1},n_{\tilde{k}})$ as the sequence of the number of intersection points induced from $\Gamma(+2,-2)$, then we can assume $n_{1} = n_{\tilde{k}}$ is not zero. 

Let $d_{1} = \# \Gamma(+2,-1)$ and $d_{2} = n_{1}= n_{\tilde{k}}$. 
We take  positive coprime integers  $1 \le \bar{p}_{2} \le p_{2}$ and $1\le \bar{q}_{2} \le q_{2}$ so that $\bar{p}_{2} q_{2} = \bar{q}_{2} p_{2} - 1$.

Then, $p'_{2} / q'_{2}$ can be written by these integers as follows. 
$$p'_{2} / q'_{2} = \frac{\bar{p}_{2}(d_{1}+d_{2})+(p_{2}-\bar{p}_{2})d_{2}}{\bar{q}_{2}(d_{1}+d_{2})+(q_{2}-\bar{q}_{2})d_{2}} = \frac{\bar{p}_{2}+p_{2}z}{\bar{q}_{2}+q_{2}z}, \text{ where } z = d_{2}/d_{1}$$.

Actually, $\bar{p}_{2}$ represents the number of the integer $d_{1}+d_{2}$ among the sequence $(n_{1},n_{2},\cdots,n_{\tilde{k}})$ and $\bar{q}_{2}$ represents the number of the integer $d_{1}+d_{2}$ among the sequence $(n_{1},n_{2},\cdots,n_{q_{2}})$. It is easy to see that these integers satisfies $\bar{p}_{2} q_{2} = \bar{q}_{2} p_{2} - 1$. Then, we find that $p'_{2} / q'_{2}$ is written as above.

Since $0 < z < +\infty$ we get that $p_{2} / q_{2}<p'_{2} / q'_{2}<\bar{p}_{2} / \bar{q}_{2}$. 

Moreover, we can prove that $\bar{p}_{2} / \bar{q}_{2}$ can be written precisely as follows.
$$ \bar{p}_{2} / \bar{q}_{2} =  k_{1} + \cfrac{1}{ k_{2} + 
						\cfrac{1}{ \cdots +
						\cfrac{1}{ k_{n} - 1}}} , \text{ if } n \text{ is odd.}$$
						
$$ \bar{p}_{2} / \bar{q}_{2} =  k_{1} + \cfrac{1}{ k_{2} + 
						\cfrac{1}{ \cdots +
						\cfrac{1}{ k_{n-1}}}} ,\text{ if } n \text{ is even.}$$
This equation comes from the Euclidean algorithm. 
Thus, we also find that the length of the continuous fraction expansion of $p'_{2} / q'_{2}$ is $n$ or grater than $n$. In particular, $R(p_{2},q_{2},p'_{2},q'_{2})$ is well-defined.

Finally, we conclude that 
$$0 > R(p_{2},q_{2},p'_{2},q'_{2}) > - ( k_{n} + \cfrac{1}{ k_{n-1} + 
						\cfrac{1}{ \cdots +
						\cfrac{1}{ k_{2} +
						\cfrac{1}{ k_{1} - \bar{p}_{2} / \bar{q}_{2}}}}} ) = -1 ,$$ if $n$ is odd, and
$$R(p_{2},q_{2},p'_{2},q'_{2}) = - ( k_{n} + \cfrac{1}{ k_{n-1} + 
						\cfrac{1}{ \cdots +
						\cfrac{1}{ k_{2} +
						\cfrac{1}{ k_{1} - p'_{2} / q'_{2}}}}} ) > 0,$$ if $n$ is even. 
\qed
\end{prf}

We consider the framed link $L(24)$ again. Our goal is to prove that $M(24)$ is in $\mathcal{L}_{\overline{\rm{Brm}}}$.

We perform the blow up operations finitely many times so that the new framed link consists of only unknots.
First, add $k_{1}$ unknots with framing $-1$ near $C_{2}$ and $C_{1}$. After that, the slopes of $C_{1}$ and $C_{2}$ are changed as follows. Denote the new slopes by $r^{1}_{\alpha_{1}}$ and $r^{1}_{\alpha_{2}}$.

\begin{itemize}
\item $r_{\alpha_{2}} \leadsto r^{1}_{\alpha_{2}} = r_{\alpha_{2}} - k_{1}q^{0}_{2}$, where $q^{0}_{2} = q_{2}$,
\item $r_{\alpha_{1}} \leadsto r^{1}_{\alpha_{1}} = r_{\alpha_{1}} - k_{1}$,
\item the new knot $C^{1}_{1}$ induced from $C_{1}$ is linking with $C_{2}$ with the slope $p^{1}_{2}/q^{1}_{2}$, where $p^{1}_{2}$ and $q^{1}_{2}$ are positive coprime integer such that $p_{2}/q_{2} = k_{1}+p^{1}_{2}/q^{1}_{2}$.
\end{itemize}

Next, add $k_{2}$ unknots with framing $-1$ isotopic to $C_{2}$. After some Kirby calculus, the slopes of $C_{1}$ and $C_{2}$ are changed as follows.
\begin{itemize}
\item $r^{1}_{\alpha_{2}} \leadsto r^{2}_{\alpha_{2}} = r^{1}_{\alpha_{2}} - k_{2}p^{1}_{2}$,
\item $r^{1}_{\alpha_{1}} \leadsto r^{2}_{\alpha_{1}} = -  \cfrac{1}{ k_{2} + 
						\cfrac{1}{- r^{1}_{\alpha_{1}}}}$,
\item the new $C^{2}_{1}$ is linking with $C_{2}$ with the slope $p^{2}_{2}/q^{2}_{2}$, where $p^{2}_{2}$ and $q^{2}_{2}$ are positive coprime integer such that $p^{1}_{2}/q^{1}_{2} = k_{2}+p^{2}_{2}/q^{2}_{2}$.
\end{itemize}

We find that $C^{2}_{1}$ is linking with $C_{2}$ with slope $p^{2}_{2}/q^{2}_{2}$. Thus, we can transform this link again by adding new $k_{3}$ unknots as above. 

\begin{figure}[h]

\includegraphics[width=7cm,clip]{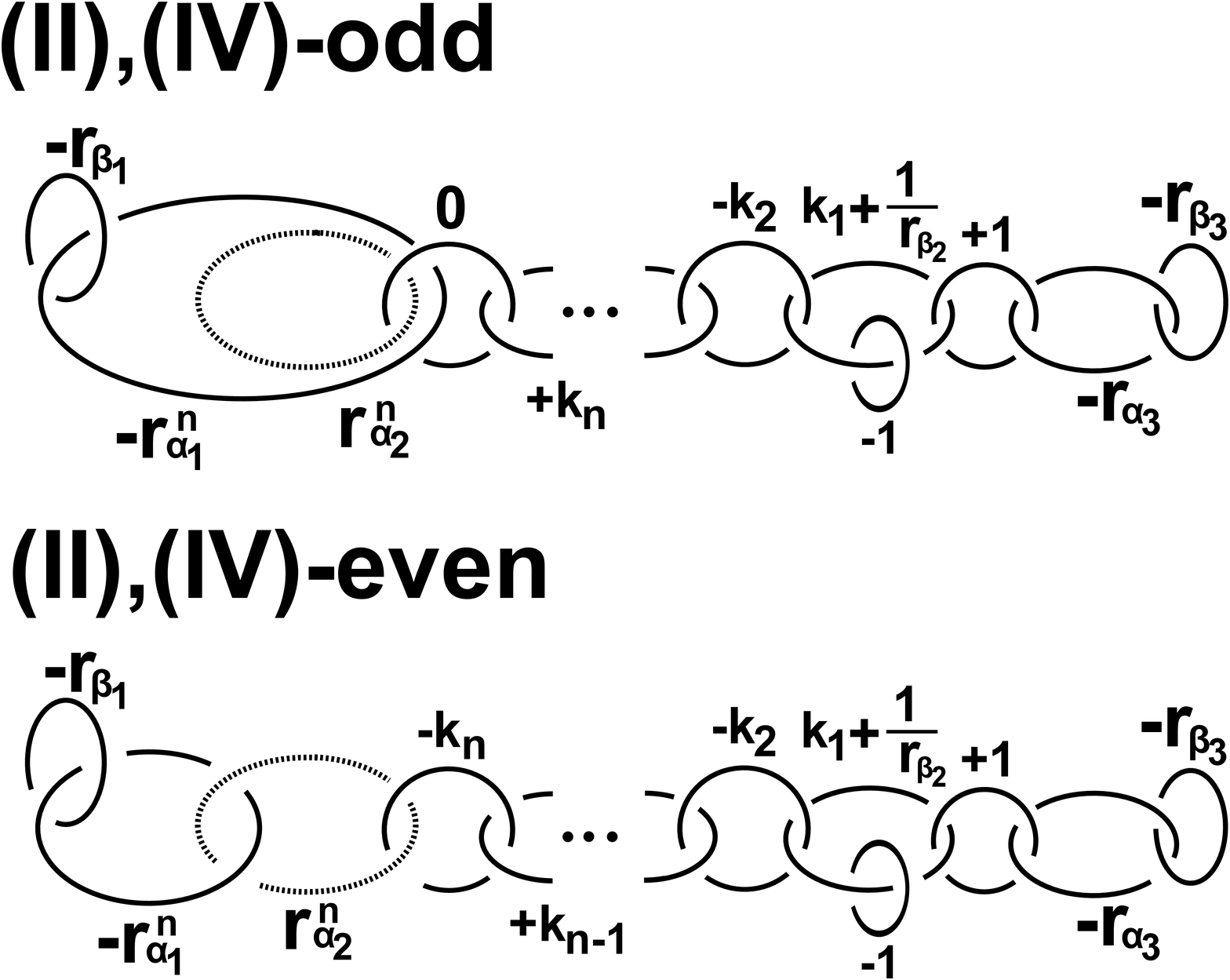} \\
\caption{}\label{type14-f}
\end{figure}

In finitely many steps, we finally get the following framed links (see Figure \ref{type14-f}).

\begin{itemize}
\item $r^{n}_{\alpha_{2}} = r_{\alpha_{2}} - k_{1}q_{2} - k_{2}p^{1}_{2}-k_{3}q^{2}_{2}-k_{4}p^{3}_{2} \cdots = r_{\alpha_{2}} - \sum_{i \text{:odd}} k_{i}q^{i-1}_{2} - \sum_{i \text{:even}} k_{i}p^{i-1}_{2} = r_{\alpha_{2}}-(p+q-1) = -p_{1}/q_{1}$,
\item $r^{n}_{\alpha_{1}} = - ( k_{n} + \cfrac{1}{ k_{n-1} + 
						\cfrac{1}{ \cdots +
						\cfrac{1}{ k_{2} +
						\cfrac{1}{ k_{1} - r_{\alpha_{1}}}}}} ) = R(p_{2},q_{2},p'_{2},q'_{2})$
\item the new $C^{n}_{1}$ is the unknot shown in Figure \ref{type14-f}, where $p^{i}_{2}$ and $q^{i}_{2}$ are positive coprime integer such that $p^{i-1}_{2}/q^{i-1}_{2} = k_{i}+p^{i}_{2}/q^{i}_{2}$.
\end{itemize}

In each case, the right part of the link can be changed to have alternating weights similarly (see Figure \ref{type1234-1}).

To prove that the left part of these link have also alternating weights, we can just use Claim \ref{clm4}.
Thus, $L(24)$ can be represented by alternatingly weighted unknots, That is, $M(24)$ belongs to $\mathcal{L}_{\overline{\rm{Brm}}}$.

\subsection{Not-$A'_{3}$-strong cases for $g = 3$}\label{notm3}

In this subsection, we study the case where the diagram is strong, but not $A'_{3}$-strong. 
Since $ME_{3} = \{[A_{3}]\}$, we get $A_{(\Sigma,\alpha,\beta)} \le A'_{3}$. Thus, it is enough to consider the case when $x_{8} = 0$, where we put

\[ A_{(\Sigma,\alpha,\beta)} = \left (\begin{array}{@{\,}cccc@{\,}}
+     & x_{2} & x_{3} \\
x_{4} & +     & x_{6} \\
0     & x_{8} & +     \end{array} \right). \]
That is, in the following argument, we do not use any conditions on $x_{2}$, $x_{3}$, $x_{4}$ and $x_{6}$.

Since $x_{8} = 0$, we get $\Gamma_{3}(+3 , +2) = \emptyset$.
Thus, the same argument as subsection \ref{nonmax} can be applied to determine this manifold as follows.

Let $x = (n_{1}, m_{1},\cdots,m_{k-1}, n_{k})$ be the sequence of integers representing $\Gamma(+3,-3)$, where $n$ means $\beta_{3} \cap (\alpha_{1}\cup \alpha_{2})$ and $m$ means $\beta_{3} \cap \alpha_{3}$. 
If $n_{1} = n_{k} = 0$, we have $n_{i} = 1$ for $1<i<k$ and we can transform the diagram by handle-slides so that $m_{1} = m_{k-1} = 0$ (see Figure \ref{nbd3}).

If $n_{1}$ and $n_{k}$ are not zero, we also have $n_{i} = 1$ for all $i$ and we can also transform the diagram by handle-slides so that $m_{i} = 0$ for all $i$ (see Figure \ref{excep1}).

Therefore, the Heegaard diagram has a lens space component. The remaining manifold has a strong Heegaard diagram with Heegaard genus two. 

\begin{prf7}
Let $(\Sigma,\alpha,\beta)$ be a strong Heegaard diagram representing $Y$ with genus three.
If the induced matrix $A_{(\Sigma,\alpha,\beta)}$ is equivalent to $A'_{3}$, we can apply Proposition \ref{A3} and Proposition \ref{A3next}. Thus, subsection \ref{surg3b}, \ref{surg3a} and \ref{determ} tell us that $Y$ is in $\mathcal{L}_{\overline{\rm{Brm}}}$. 
If $A_{(\Sigma,\alpha,\beta)}$ is not equivalent to $A'_{3}$, we return to the genus two case.
Finally, the genus two case are proved in section \ref{g2case}.
\qed
\end{prf7}

\section*{acknowledgement}
I would like to express my deepest gratitude to Prof. Kohno who provided helpful comments and suggestions. I would also like to express my gratitude to my family for their moral support and warm encouragements.

\end{document}